 \input amstex
 \input epsf

 \documentstyle{amsppt}

 \def\al{\alpha}
 \def\th{\theta}
 \def\be{\beta}
 \def\ep{\varepsilon}
 \def\de{{\delta}}
 \def\ka{\kappa}
 \def\la{\lambda}
 \def\La{\Lambda}
 \def\ga{\gamma}
 \def\Ga{\Gamma}
 \def\si{\sigma}
 \def\Si{\Sigma}
 \def\de{{\delta}}
 \def\De{\Delta}
 \def\th{\theta}
 \def\Th{\Theta}
 \def\ze{\zeta}
 
 \def\Om{\Omega}
 
 \def\Up{\Upsilon}

 
 \def\cB{\Cal B}
 \def\cC{\Cal C}
 \def\cD{\Cal D}
 \def\cE{\Cal E}
 
 \def\cF{\Cal F}
 \def\cH{\Cal H}
 
\def\cJ{\Cal J}
 
 \def\cL{\Cal L}
 \def\cM{\Cal M}

 \def\cV{\Cal V}
 \def\cW{\Cal W}

 \def\cZ{\Cal Z}
 \def\cL{\Cal L}
 \def\cP{\Cal P}
 \def\cQ{\Cal Q}
 \def\cS{\Cal S}

 \def\limt{\lim_{t\to\infty}}
 \def\limn{\lim_{n\to\infty}}

 \def\wt{\widetilde}
 \def\wh{\widehat}
 \redefine\l{\ell}
 
 \def\qqquad{\quad \quad \quad}
 \def\sobre#1#2{\lower 1ex \hbox{ $#1 \atop #2 $ } }
 \def\bajo#1#2{\raise 1ex \hbox{ $#1 \atop #2 $ } }

 \magnification=\magstep1
 \vcorrection{.5truecm}
 \vsize 19cm

\global\newcount\numsec\global\newcount\numfor
\global\newcount\numfig
\global\newcount\numcon
\global\newcount\numkon
\gdef\profonditastruttura{\dp\strutbox}
\def\senondefinito#1{\expandafter\ifx\csname#1\endcsname\relax}
\def\SIA #1,#2,#3 {\senondefinito{#1#2}
\expandafter\xdef\csname #1#2\endcsname{#3}\else
\write16{???? ma #1,#2 e' gia' stato definito !!!!} \fi}

\def\etichetta(#1){(\veroparagrafo.\veraformula)
\SIA e,#1,(\veroparagrafo.\veraformula)
 \global\advance\numfor by 1
 \write16{ EQ \equ(#1) == #1  }}

\def\letichetta(#1){\veroparagrafo.\veraformula
\SIA e,#1,{\veroparagrafo.\veraformula}
\global\advance\numfor by 1
 \write16{ Sta \equ(#1) == #1}}

\def\tetichetta(#1){\veroparagrafo.\veraformula 
\SIA e,#1,{(\veroparagrafo.\veraformula)}
\global\advance\numfor by 1
 \write16{ tag \equ(#1) == #1}}

\def \FU(#1)#2{\SIA fu,#1,#2 }

\def\etichettaa(#1){(A\veroparagrafo.\veraformula)
 \SIA e,#1,(A.\veroparagrafo.\veraformula)
\global\advance\numfor by 1
 \write16{ EQ \equ(#1) == #1  }}

\def\getichetta(#1){Fig. \verafigura
 \SIA e,#1,{\verafigura}
 \global\advance\numfig by 1
 \write16{ Fig. \equ(#1) ha simbolo  #1  }}

\def\cetichetta(#1){C_\veracon
 \SIA c,#1,{\veracon}
 \global\advance\numcon by 1
 \write16{ Const \cequ(#1) =  #1 }}

\def\ketichetta(#1){K_{\verakon}
 \SIA k,#1,{\verakon}
 \global\advance\numkon by 1
 \write16{ Const \kequ(#1) =  #1 }}

\newdimen\gwidth
\def\BOZZA{
\def\alato(##1){
 {\vtop to \profonditastruttura{\baselineskip
 \profonditastruttura\vss
 \rlap{\kern-\hsize\kern-1.2truecm{$\scriptstyle##1$}}}}}
\def\galato(##1){ \gwidth=\hsize \divide\gwidth by 2
 {\vtop to \profonditastruttura{\baselineskip
 \profonditastruttura\vss
 \rlap{\kern-\gwidth\kern-1.2truecm{$\scriptstyle##1$}}}}}
}
\def\alato(#1){}
\def\galato(#1){}
\def\calato(#1){}
\def\kalato(#1){}
\def\veroparagrafo{\number\numsec}
\def\veraformula{\number\numfor}
\def\verafigura{\number\numfig}
\def\veracon{\number\numcon}
\def\verakon{\number\numkon}
\def\geq(#1){\getichetta(#1)\galato(#1)}
\def\Eq(#1){\eqno{\etichetta(#1)\alato(#1)}}
\def\Ceq(#1){\cetichetta(#1)\calato(#1)} 
\def\Keq(#1){\ketichetta(#1)\kalato(#1)} 
\def\teq(#1){\tag{\tetichetta(#1)\hskip-1.6truemm\alato(#1)}}  
\def\eq(#1){\etichetta(#1)\alato(#1)}
\def\Eqa(#1){\eqno{\etichettaa(#1)\alato(#1)}}
\def\eqa(#1){\etichettaa(#1)\alato(#1)}
\def\eqv(#1){\senondefinito{fu#1}$\clubsuit$#1\else\csname fu#1\endcsname\fi}
\def\equ(#1){\senondefinito{e#1}\eqv(#1)\else\csname e#1\endcsname\fi}
\def\cequ(#1){C_\senondefinito{c#1}\eqv(#1)\else\csname c#1\endcsname\fi}
\def\kequ(#1){K_\senondefinito{k#1}\eqv(#1)\else\csname k#1\endcsname\fi}
\def\Lemma(#1){Lemma \letichetta(#1)\hskip-1.6truemm}
\def\Theorem(#1){{Theorem \letichetta(#1)}\hskip-1.6truemm}
\def\Proposition(#1){{Proposition \letichetta(#1)}\hskip-1.6truemm}
\def\Corollary(#1){{Corollary \letichetta(#1)}\hskip-1.6truemm.}
\def\Remark(#1){{\noindent{\bf Remark \letichetta(#1)\hskip-1.6truemm.}}}

\def\include#1{
\openin13=#1.aux \ifeof13 \relax \else
\input #1.aux \closein13 \fi}
\openin14=\jobname.aux \ifeof14 \relax \else
\input \jobname.aux \closein14 \fi

 \TagsOnRight
 \NoBlackBoxes
 \document
 \topmatter
 \title
 A phase transition in a model for the spread of an infection  
 \endtitle
 \author 
 Harry Kesten and Vladas Sidoravicius
 \endauthor
 \leftheadtext{Harry Kesten and Vladas Sidoravicius}
 \rightheadtext{Phase transition for an infection}
 \abstract
 We show that a certain model for the spread of an infection has a
 phase transition in the recuperation rate. The model is as follows:
 There are particles or individuals of type $A$ and type $B$, interpreted
 as healthy and infected, respectively. All particles perform
 independent, continuous time, simple random walks on $\Bbb Z^d$ with
 the same jump rate $D$. The only interaction between the particles is
 that at the moment when a $B$-particle jumps to a site which contains
 an $A$-particle, or vice versa, the $A$-particle turns into a
 $B$-particle. All $B$-particles recuperate (that is, turn back into
 $A$-particles) independently of each other at a rate $\la$. 
  We assume that we start the system with $N_A(x,0-)$
 $A$-particles at $x$, and that the $N_A(x,0-), \, x \in \Bbb Z^d$, 
 are i.i.d., mean $\mu_A$ Poisson random variables. In addition we
 start with one additional $B$-particle at the origin. We show that
 there is a critical recuperation rate $\la_c > 0$ such that the 
 $B$-particles survive (globally) with positive probability if 
$\la < \la_c$ and
 die out with probability 1 if $\la > \la_c$. 
 \endabstract

 \address
 Harry Kesten, 
 Department of Mathematics,
 Malott Hall,
 Cornell University,
 Ithaca NY 14853, USA
 \endaddress
 \email
 kesten\@math.cornell.edu
 \endemail

 \address
 Vladas Sidoravicius,
 IMPA,
 Estr. Dona Castorina 110,
 Rio de Janeiro,
 Brasil,
 \endaddress
 \email
 vladas\@impa.br
 \endemail

 \endtopmatter

 \subhead
 1. Introduction
 \endsubhead
 \numsec=1
 \numfor=1
 In \cite{KSb},\cite {KSc} we investigated the model discussed in the
 abstract, but without recuperation, that is, with $\la = 0$ only.
 We heard of the present version from Ronald Meester and we also
 learned from him the conjecture that there would be a phase transition
 in $\la$, as is now confirmed by our principal theorem here. 

 Before formally stating our theorem we 
 make some comments about the precise formulation
 of the model, and introduce some notation.
First we define for $\eta = A$ or $B$
 $$
 N_\eta(x,t) = \text{number of $\eta$-particles at the space-time point
 $(x,t)$}.
 $$
 Throughout we write $\bold 0$ for the origin. As stated in the
 abstract, we put $N_A(x,0-)$ $A$-particles at $x$ just
 before we start. We then introduce a $B$-particle at the origin and
 turn some of the particles at the origin instantaneously to $B$-particles, so
 that at time 0 we start with $N_A(x,0)= N_A(x,0-)$ $A$-particles at 
 $x \ne \bold 0$ and $N_B(\bold 0,0) \in [1, N_A(\bold 0, 0-) +1]$
 $B$-particles at $\bold 0$. However, at any time $t > 0$ 
 an $A$-particle can 
 turn into a $B$-particle only if the $A$-particle itself jumps at $t$
 or if some $B$-particle jumps to the position of the $A$-particle at
 time $t$.
 Thus, we are not saying that an $A$-particle turns  into a
 $B$-particle whenever it coincides with a $B$-particle. We adopted
 the rule that a jump is required for the following reason. If we did
 not make this requirement, then $B$-particles could effectively not
 recover at a space-time point $(x,t)$ with several $B$-particles present.
 Indeed, if one of them tried to turn back into an $A$-particle at time
 $t$, it would immediately become of type $B$ again because it coincided with
 another $B$-particle. This creates some sort of singularity in the
 model which we are unable to handle at the moment (see, however,
 Remark 3 below). This is the reason
 for the requirement of a jump for a change from type $A$ to type
 $B$ at all strictly positive times $t$. Only at $t=0$ did we change
 some $A$-particles at $\bold 0$ to $B$-particles because they coincided
 with a $B$-particle (even though no jump occurred). 
The choice of the set of $A$-particles at
 $\bold 0$ which is turned into $B$-particles at time 0 
 will not influence our arguments. 
 Note that because of the jump requirement there may be particles of
 both types at a single space-time point.

 We have not attempted to give a formal proof of the existence of our
 process here as a strong Markov process on a suitable probability
 space. We did carry out such a proof for the model without
 recuperation in \cite{KSc}, and this indicates that such an existence
 proof for the present model is  probably non-trivial, and in any case
 rather tedious. Probably one can build on the proof for the case
 without recuperation, because there are fewer $B$-particles in the
 model with recuperation than in the one without recuperation, as
 shown in Corollary 3 below. We merely mention that in \cite {KSc} our basic
 probability space for the process without recuperation was a subset of
 the collection of right continuous paths with left limits 
 from $[0,\infty)$ into 
 $$
 \Si:= \prod_{k \ge 1}\big( (\Bbb Z^d \cup \partial_k) 
 \times \{A,B\}\big).
 \teq(1.1)
 $$
 The $\partial_k$ are cemetery points which we can ignore here, since
 the process is defined such that it almost surely does not reach any
 of these points. The initial particles are ordered in some way as 
 $\rho_1,\rho_2, \dots$. A typical point of $\Si$ is written as $\si =
 \big(\si'(k), \si''(k)\big)_{k \ge 1}$. $t \mapsto \big(\si'_t(k), 
\si''_t(k)\big)$
 is a path from $[0,\infty)$ into $\Bbb Z^d \times \{A,B\}$. The value
 of this path at time $t$ represents the position and type of $\rho_k$
 at time $t$. We often write $\pi(t, \rho_k)$ and $\eta(t, \rho_k)$ 
 for the position and type of $\rho_k$ at time $t$.
 Thus we have attached to each particle $\rho$ 
 a  path $t \mapsto \pi(t, \rho)$. The quantity 
 $\{\pi_A(t, \rho):= \pi(t,\rho)-\pi(0,\rho)\}_{t \ge 0}$ 
gives the displacement at time $t$ of $\rho$ from its
 starting point. The paths $\pi_A(\cdot, \rho)$
for the different $\rho$ are all taken  
 as independent copies of a continuous time simple random walk
 $\{S_t\}_{t \ge 0}$ with jump rate $D$ and starting point $S_0 = \bold
 0$. 
The type of $\rho_k$ at time $t$ is a complicated function of the
 initial types and the restrictions to $[0,t]$ of all the paths
 $\pi_A(\cdot, \rho)$. More
 details of dependence of the types as functions of the paths can be
 found in Section 2 of \cite {KSc}. 

 In the case were recuperation is allowed, as in the present article,
 we further attach to each particle $\rho$ a sequence of
 potential recuperation times $r(1, \rho) < r(2, \rho) <\dots $. The
 $r(i,\rho)$ are the jump times of a rate $\la$ Poisson process, and
 these processes are all independent of each other for different $\rho$
 and independent of the $\pi(\cdot,\rho)$. 
If $\rho$ is of type $B$ at a time $t$,
 then its type will turn back to $A$ at the first $r(i,\rho) \ge t$.
 A great advantage of the assumption that the random walks are
 independent of the types is that the $\pi(\cdot,\rho)$ and the
 $r(i,\rho)$ can be determined once and for all at time 0. The actual
 evolution of the type of each particle over time is then a complicated
 function of all the paths and recuperation times for all particles.
 We shall make a few more comments about this function in the beginning
 of Section 2. 
 We point out that another reason for our restriction to 
 the case of equal random walks for the different types is that 
 the basic monotonicity properties of the next section
 may fail if the random walks are different for the different types. 

 We say that the infection {\it survives} if 
 $$
 P\{\text{there are some $B$-particles at all times}\} > 0.
 \teq(1.1a)
 $$
 Since there cannot be any $B$-particles after time $t$ if there are no
 $B$-particles at $t$, it follows that \equ(1.1a) is equivalent to
 $$
 \limt P\{\text{there are some $B$-particles at time }t\} > 0.
 \teq(1.2)
 $$
 One may even replace $\limt$ by $\liminf_{t \to \infty}$ in
 \equ(1.2). Note that the survival in \equ(1.1a) or \equ(1.2) is only  
 global survival. Local survival in its strongest form would say that
 $$
 \liminf_{t \to \infty}P\{N_B(\bold 0,t) > 0\} > 0.
 \teq(1.2aa) 
 $$
 A weaker form of local survival would be that
 $$
 P\{N_B(\bold 0,t)>0 \text{ for arbitrily large }t \} > 0.
 \teq(1.2ab) 
 $$
 Clearly \equ(1.2aa) implies \equ(1.2ab), and this, in turn implies \equ(1.1a).
 We do not know how to prove that either of the forms \equ(1.2aa) or 
  \equ(1.2ab) of local survival holds if $\la$ is small enough.

 The infection is said to {\it die out} or to {\it become
 extinct} if it does not survive, i.e., if
 $$
 P\{\text{there is some (random) $t$ such that there are no
 $B$-particles after $t$}\}=1.
 \teq(1.3)
 $$  

 Here is our principal result.
 \proclaim{Theorem 1} There exists a $0 < \la_c < \infty$
 such that the infection survives if $\la < \la_c$ and dies out if $\la
 > \la_c$.
 \endproclaim

 \noindent
 {\bf Remark 1.} The restriction to only one $B$-particle at time 0 is
 for convenience only. The theorem remains valid if we start with any
 finite number of $B$-particles at (nonrandom) positions.
 \medskip \noindent
 {\bf Remark 2.} 
We already remarked that the theorem does not give local
 survival if $\la$ is sufficiently small. Neither does it tell us
 anything about the location of the $B$-particles as a function of $t$
 on the event that the $B$-particles survive forever.

 By a special argument one can show that \equ(1.2ab) holds for $d=1$
 and $\la < \la_c$ on the event that the $B$-particles survive forever. 
 \medskip \noindent
 {\bf Remark 3.} The proof that there is survival for small $\la >0$
 works even in the case in which an $A$-particle turns into a
 $B$-particle whenever it coincides with a $B$-particle, that is, 
 if we do not require that the $A$ or $B$-particle jumps before
 reinfection can occur after recuperation of a $B$-particle.
 \medskip \noindent
 {\bf Remark 4.}  A similar result for another variant of the model 
is obtained in \cite{AMP}. This article considers the so-called frog
 model in which only the $B$-particles move and
  the $A$-particles stand still. In \cite {AMP} time is taken discrete.
 It is assumed that
 each $B$-particle is removed from the system at its
 first recuperation. One could interpret this by means of the introduction
 of a third type of particles, namely imune ones which do not interact
 with any particles. When a particle recuperates from the infection it 
 becomes immune. This results also in some conclusions which differ
 from the ones in the present paper. In particular,
 \cite{AMP} shows that in their case there never is
 survival in dimension 1, if recuperation is allowed (i.e., $\la_c =
 0$, so that {\it there is no nontrivial phase
 transition in dimension 1}, in contrast to our model). 

 The fact that the $A$-particles can move in our model makes the
 analysis here much harder than in \cite {AMP}. This also forces us to
 stick to Poisson initial conditions, while \cite {AMP} can handle much
 more general initial conditions, as well as more general graphs as
 $\Bbb Z^d$.  

We note that our proof of survival 
in Section 3 still goes through if the $A$ and
$B$-particles perform the same random walk and $B$-particles are
immune after recuperation. In this case one also has extinction for
large $\la$ by Theorem 1 and monotonicity arguments as in Lemma 4 below. 
The system in which $B$-particles become immune lies stochastically
below the system we are investigating here (in the sense of Lemma
4). Thus Theorem 1 remains valid if $B$-particles are
immune after recuperation.
 \medskip \noindent
 {\bf Remark 5.} The following version of the frog model can still be
 analyzed to some extent. Take time continuous, and assume that the
 $A$-particles cannot move. Assume further
 that a $B$-particle turns any $A$-particle with which it coincides 
 instantaneously into a $B$-particle. $B$-particles turn back to
 $A$-particles at a constant rate $\la > 0$, but these recuperated
 particles stay in the system and act as any original $A$-particle. 
 For the initial state take the $N_A(x,0-)$ as i.i.d., mean $\mu_A$
 Poisson variables, and add one $B$-particle at $\bold 0$. We have not
 constructed such a process, but we take it
 for granted that this process can be properly defined so as to justify
 the argument below.

 \cite {AMP} proved survival for the process in discrete time,
in which the $A$-particles
 stand still, $\la$ is small, and 
in which particles which recuperate are removed from the
 system, and $d \ge 2$.
 We expect that this also holds for the process just described. 
 It is perhaps surprising, though, that the 
rules of the preceding paragraph imply that for large 
$\mu_A$ the process {\it always survives}.
 More precisely, we show that
 if $\mu_A >$ some $\mu_{A,d}$ (which depends on the dimension
 $d$ only), then the $B$-particles survive for
 all values of $\la$, so that {\it there is no phase transition}.

 The key observation for proving this lack of a phase transition 
 is that if there are several particles
 present at some space-time point $(x,t)$, then they are all of type $A$
 or all of type $B$. In the latter case, if one of the $B$-particles
 tries to recuperate, it is immediately reinfected by the other
 $B$-particles at the same location, and so, as long as 
there are at least two particles on one site, none of the particles at
 that site can change from type $B$ to $A$. This shows that
 $B$-particles can turn back to $A$-particles only at sites with no
 other particle. Since the $A$-particles stand still, it follows that,
at any fixed site, at most one $B$-particle can recuperate and stay of
 type $A$ forever after.

 We shall also use that for $\mu_A > $ some $\mu_{A,d}$ it holds
 $$
 \sum \Sb C \text{ connected}\\ \bold 0 \in C \endSb
 P\big\{\sum_{x \in C}N_A(x,0-) < \frac 12 \mu_A|C|\big\} < \infty.
 \teq(1.3aa)
 $$
 This follows from standard large deviation estimates for the Poisson
 distribution, since $\sum_{x \in C}N_A(x,0-)$
 has a Poisson distribution with mean $\mu_A|C|$, and from the fact
 that the number of connected sets $C$ with $\bold 0 \in C$ grows only
 exponentially in $|C|$. It follows from
 \equ(1.3aa) and the Borel-Cantelli lemma that for $\mu_A > \mu_{A,d}$,
 almost surely there
 exists some random $k_0$ such that for any connected set $C \subset
 \Bbb Z^d$ which contains $\bold 0$ and with $|C| \ge k_0$,
 $$
 \sum_{x \in C}N_A(x,0-)\ge \frac 12 \mu_A|C|.
 \teq(1.3ab)
 $$

 Assume now that there exist
 $k_0$ {\it distinct} particles $\rho_1, \dots, \rho_{k_0}$,
 and space-time points $(x_i,t_i)$, such that $\rho_i$ is
 at $x_i$ at time $t_i$ as a $B$-particle. (Some of the $x_i$ or $t_i$
 with different $i$ may have the same value, 
but the $\rho_i$ have to be distinct.) 
 Assume further that 
 $$
 A_0 :=\{x_1, \dots, x_{k_0}\} \text{ is connected}.
 \teq(1.3ac)
 $$
 Assume also
 that the infection dies out at some time $t_\infty < \infty$.
 Let $C_0$ be the collection
 of sites visited by one of the $\rho_i$ before the infection dies out,
 or more precisely
 $$
 C_0 := \{x:\text{ for some $1\le  i \le k_0,\;\rho_i$ visits $x$
 during $[t_i,t_\infty]$}\}.
 \teq(1.3ad)
 $$
 $C_0$ is again a connected set, because each particle $\rho_i$ moves
 by a simple random walk through a connected set.
 Next, let $D_0$ be the 
 collection of sites at which the $k_0$ particles $\rho_i$ are at time
 $t_\infty$ (and hence also at $t > t_\infty$, because each $\rho_i$
 must have type $A$ from the time of extinction of the infection
 on). Then, by the one but last paragraph, 
 $$
 |C_0| \ge |D_0| = k_0.
 $$
 Now, let $\ze$ be some particle at some $x \in
 C_0$ at time 0 (if such a particle exists). 
Then $x$ is visited by some $\rho_i$ at some time $s_i
 \in [t_i, t_\infty]$. Pick such an $i$ and let $s_i$ be the smallest time in
 $[t_i,t_\infty]$ at which $\rho_i$ is at $x$. We claim that $\rho_i$
 must have type $B$ at time $s_i$. Indeed, if $s_i = t_i$, this is
 true by our assumption on $\rho_i$ at $(x_i,t_i)$. If $s_i > t_i$,
 then $\rho_i$ must jump to $x$ at time $s_i$. But only $B$-particles
 do jump, so that our claim also holds in this case. Now, 
either $\ze$ has type $B$ at some time during $[0,s_i)$,
 or is of type $A$ and sits still at $x$ during all of $[0,s_i]$ 
and then it is turned into
 type $B$ by $\rho_i$ at $s_i$. In either case, the infection
 cannot die out before $\ze$ too recuperates for a last time. 
 But, by \equ(1.3ab) the number of particles in $C_0$ at time 0 is at least
 $$
 \sum_{x \in C_0}N_A(x,0-)\ge \frac 12 \mu_A|C_0| \ge \frac 12
 \mu_Ak_0 \ge 2 k_0,
 \teq(1.3abb)
 $$
 provided we take $\mu_A \ge \mu_{A,d}\ge 4$. Thus we now have found at
 least $2k_0$ particles  which must recuperate during $[0, t_\infty]$.
 We can repeat the
 argument with he collection $\rho_1, \dots, \rho_{k_0}$ replaced by
 the particles in $C_0$, and $A_0$ replaced by $A_1 := C_0$. $k_0$ is
 then replaced by some $k_1 \ge 2k_0$. By repeating this argument
 infinitely often we see that it is impossible for the infection
 to die out in finite time, if there is a $k_0$ such that \equ(1.3ab)
 holds for all connected $C \supset \{\bold 0\}$ with $|C| \ge k_0$, and
 particles $\rho_i, \; 1 \le i \le k_0$, as above. Here we have taken
 it for granted (without proof) that in any reasonable version of the process
 only finitely many $B$-particles can be formed in finite time.
 We apply the preceding remarks with $A_0 = \{\bold 0\}$
 and a large non-random $k_0$. This  shows that
 $$
 \align
 &P\{\text{infection does not die out}\}\\
 &\ge P\{\text{\equ(1.3ab)
 holds for all connected $C \supset \{\bold 0\}$ with $|C| \ge k_0$ and
 $N_A(\bold 0,0-) \ge k_0$}\}\\
 &\ge P\{\text{\equ(1.3ab)
 holds for all connected $C \supset \{\bold 0\}$ with $|C| \ge k_0$}\}
 P\{N_A(\bold 0,0-) \ge k_0\}\\ 
 &\text{(by the Harris-FKG inequality) }
 > 0.
 \endalign
 $$
 This argument is independent of the value of $\la$ and therefore 
 proves that there is no phase transition.

 \bigskip
 In the next section we begin with a monotonicity property which
 immediately implies that there exists a critical $\la_c$ with the
 properties stated in Theorem 1, except that $\la_c = 0$ or $\infty$ is
 still not excluded. In Section 3 we then show that $\la_c >0$
 and in Section 5 we show
 that $\la_c < \infty$. These two sections which show that 
 there are nontrivial regions of survival and extinction, respectively,
 form the core of this paper. Section 4 is a kind of interlude in which
 we prove that the maximal number of jumps during $[0,t]$ in 
 a certain class of paths is at most $O(t)$. This estimate is crucial
 for the proof of extinction in Section 5.

 Our methods are a combination of the
 multi-scale analysis of \cite{KSa}, \cite {KSb} 
and percolation arguments. To show that the infection survives for
 small $\la$ we introduce (in Section 3) a certain
 directed (dependent) percolation process  with the property that 
 if percolation occurs in this process, then the infection survives. 
 We then show that percolation occurs for sufficiently small $\la$
 by showing that there is only a very small probability that the origin is
 separated from $\infty$ by a distant separating set.
To show that the infection dies out when $\la$ is large we use a block
 argument (in Section 5). We show that with high probability, 
along ``almost all'' paths in space-time there have to be blocks which prevent
 the transmission of the infection. The paragraph
 following the statement of Proposition 24 in Section 5 gives some
 more details of this strategy.

A reader interested in the details of the 
proofs will have to refer to \cite {KSa}-\cite {KSc} 
a number of times.

 Throughout this paper we make the following convention about
 constants. $K_i$ will denote a strictly positive, finite constant,
 whose precise value is unimportant for our purposes. The value of the
 same $K_i$ may be different in different formulas. 
 We use $C_i$ for constants whose value remains fixed throughout the
 paper. They will again have values in $(0, \infty)$.
 If necessary, we
 indicate on what quantities a constant depends at the time when
 it is first introduced.
 Throughout $\|x\|$ denotes the $\l^\infty$ norm of the vector $x =
 (x(1), \dots, x(d)) \in \Bbb R^d$, i.e., 
 $$
 \|x\| = \max_{1 \le i \le d} |x(i)|.
 \teq(2.0)
 $$ 
 $$
 \cC(m) = \{x:\|x\|\le m\} = [-m,m]^d.
 \teq(2.0a)
 $$
 $\bold 0$ will denote the origin (in $\Bbb Z^d$ or $\Bbb R^d$);
 $|C|$ usually denotes the cardinality of the set $C$.

 \bigskip
 \noindent
 {\bf Acknowledgements.} We thank Ronald Meester for bringing 
 several of the questions studied here to our attention.

 Much of the research for this paper was 
 carried out by visits of one or both authors to Eurandom in Eindhoven
 and the Newton Institute for Mathematical Sciences in Cambridge.
 HK thanks Eurandom for appointing him as Eurandom Professor in the fall
 of 2002. 
 He also thanks Eurandom and the Newton Institute
 for their support and for their hospitality during his visits.
 Further support for HK
  came from the NSF under Grant DMS 9970943 and from Eurandom.

 VS thanks Cornell University and the
 Newton Institute for their support and hospitality during his visits
 to these institutions.
 His research was further supported by FAPERJ Grant E-26/151.905/2001, 
 Pronex (CNPq-Faperj).

 \subhead
 2. Two monotonicity properties
 \endsubhead
 \numsec=2
 \numfor=1
 We repeat that {\it we assume that all particles perform copies of 
 the same random walk}.
 In this section we show that increasing the recuperation rate decreases
 the number of infected particles. In addition we repeat a monotonicity
 property from \cite {KSb} for the system without recuperation.

 First some recapitulation of the  
 notation used in \cite {KSb}, \cite {KSc} for
 the construction of a suitable Markov process.
 $\Si_0$ is a subset of
 $\Si$ (defined in \equ(1.1)) which serves as the state space for a
 strong Markov process $\{Y_t\}_{t \ge 0}$ constructed as a suitable
 version of our infection process without recuperation. 
 For our purposes here we do not have
 to know the exact definition of $\Si_0$, but we merely have to know
 that the initial conditions, as described by the Poisson variables
 $N_A(x,0-)$, lie almost surely in $\Si_0$ (by Proposition 4 of \cite {KSb}),
 and that then the Markov chain takes values in $\Si_0$ for all times,
 almost surely. Moreover, we have from Section 2 in \cite{KSb}(see
 (2.18) there), that almost surely  
 $$
 \sup_{s\le t} \text{(number of $B$-particles at time $s$ in the
 process $\{Y_t\}$)} < \infty.
 \teq(2.0bb)
 $$ 
 $\Si_0$ will also be the state space for the infection
 process with recuperation. We write $\{Y_t(\la)\}$ for the process
 with recuperation rate $\la$, even when $\la =0$. The process $\{Y_t(0)\}$, 
 does not allow recuperation, but it is not the same as the process
 $\{Y_t\}$ of \cite {KSa}, \cite {KSb}. In the former process an
 $A$-particle turns into a $B$-particle only when one of these two
 particles jumps to the position of the other. In particular this
 process can have $A$ and $B$-particles at the same site. 
In the process $\{Y_t\}$ this is not possible, because an $A$-particle
 turns instantaneously to a $B$-particle when it coincides with a
 $B$-particle. The difference between these two processes, even
 though it is small, forces some extra work on us.

 To motivate our construction for $\{Y_t(\la)\}$ consider a particle $\rho$
 which is of type $B$ at time $s$ in the process $\{Y_t(\la)\}$, and
 which has changed type only finitely often in this process. Such a
 particle should have an analogue of a genealogical path as introduced in
 Proposition 4 in \cite {KSb} in $\{Y_t\}$.  Specifically, there should be
 space-time points $(x_i,s_i)$ with 
 $1\le i \le \l$ for some $\l$,
 and $0 < s_1 < \dots < s_\l < s$, and particles $\rho_i$ for $0 \le i \le
 \l+1$ with $\rho_{\l+1} = \rho$, such that at time $s_i$, 
 $\rho_i$ jumps to the position of $\rho_{i+1}$ or vice versa. Moreover,
 (with $s_{\l+1} = s$) $\rho_0$ should have type $B$ at time 0, and
 $\rho_i$ should have type $B$ and not recuperate
 during $[s_i, s_{i+1}]$ in $\{Y_t(\la)\}$.
 This last requirement was of course not
 present in \cite{KSb}, but nevertheless the backwards construction 
 of the genealogical path from \cite{KSb} works 
 with only trivial modifications. To be more specific, 
 start with $\rho$ of
 type $B$ at time $s$ and find the time $t_1 := \min\{u: \rho$ has 
 type $B$ in $\{Y_t(\la)\}$ during $[u,s]\}=\min\{u: \rho$ 
 does not recuperate in $\{Y_t(\la)\}$ during $[u,s]\}$.
 Then, either $t_1=0$ or $t_1 >0$. If $t_1 = 0$ then $\rho$
was of type $B$ at time 0 and did not recuperate during $[0,s]$ and 
we are done. If $t_1 > 0$, then there must have been some other particle 
 $\rho^{(1)}$ of type $B$ in $\{Y_t(\la)\}$, 
 and this $\rho^{(1)}$ must have jumped to the position of $\rho$,
 or vice versa, at time $t_1$. We then define $t_2 = \min\{u:
 \rho^{(1)} \text{ has 
 type $B$ in $\{Y_t(\la)\}$  during }[t_2,t_1]\}$, etc., 
 until we arrive, for some $\l$ 
 at time $t_\l$ and a 
 particle $\rho^{(\l+1)}$ which had type $B$ in $\{Y_t(\la)\}$ 
 during $[0, t_\l]$. The genealogical path for $\rho$ in $\{Y_t(\la)\}$ 
 is then obtained by using the $t_i$ and $\rho^{(i)}$ in
 reverse order for the $s_i$ and $\rho_i$. 
 Note that if $\rho$ is of type $B$ at time $s$ and 
 has a genealogical path of times $0 < s_1 < \dots < s_\l < s_{\l+1} =
 s$ and corresponding particles $\rho_i$, 
 in the process
 $\{Y_t(\la)\}$, then $\rho$ can also be regarded as a $B$-particle at
 time $s$ in the process $\{Y_t\}$.
 Indeed, one easily shows by induction on $i$ that each of the
 particles $\rho_i$ must have type $B$ at time $s_i$ in $\{Y_t\}$. 
 (Note that we are
 not saying that $\rho_i$ changes type from $A$ to $B$ at time $s_i$ in
 $\{Y_t\}$; the
 argument here does not rule out that $\rho_i$ is already of type $B$
 just before $s_i$, but this does not matter.)

 With the motivation provided by the preceding paragraph we construct
 $\{Y_t(\la)\}$ on the product of the probability space for $\{Y_t\}$
 with the probability space for all the recuperation processes
 $\{r(i,\rho)\}$. For a generic point $\si = \big(\si'(k),
 \si''(k)\big)$  in the state space $\Si$ (see \equ(1.1)) define
 $\overline \si$ to be the point obtained from $\si$ by taking
 $\si''(k) = B$ for all $k$ for which there is an $\l$ with $\si'(\l) =
 \si'(k)$ and $\si''(\l) = B$. This means that $\overline \si$ is
 obtained from $\si$ by changing to $B$ the type of all particles at   
a position which already has at least one $B$-particle. 
 We now describe the process $\{Y_t(\la)\}$ starting from a
 $\si$ for which $\overline \si \in \Si_0$. In \cite {KSb}, \cite
 {KSc} we defined 
the process $\{Y_t\}$ starting from $\overline \si$. This begins with
 assigning to each particle $\rho$ a random walk path $\pi_A(\cdot,
 \rho)$ and then giving $\rho$ the position $\pi(t,\rho) = \pi(0,\rho) +
 \pi_A(t,\rho)$ at time $t$, where $\pi(0,\rho)$ is just the initial
 position of $\rho$. We now assign to $\rho$ the same positions
 $\{\pi(t,\rho)\}_{t \ge 0}$ in $\{Y_t(\la)\}$.  
To complete the description we merely have to 
decide what type to assign to a
 particle $\rho$ as a function of time in $\{Y_t(\la)\}$. If $\rho$ has
 type $A$ at time $s$ in $\{Y_t\}$ starting from $\overline \si$, 
then we also assign it type $A$ at
 time $s$ in $\{Y_t(\la)\}$. In particular, since almost surely 
 only finitely many
 particles meet a particle of type $B$ during $[0,s]$ in $\{Y_t\}$
 (by \equ(2.0bb) and the fact that any particle which meets a
 $B$-particle 
before time $s$ has type $B$ at time $s$ in $\{Y_t\}$),
 this rule also assigns type $A$ during $[0,s]$  to all but finitely
 many particles in $\{Y_t(\la)\}$. Let $\rho^{(1)}, \dots,\rho^{(m)}$
 be the finitely many particles of type $B$ at time $s$ in $\{Y_t\}$.
 The particles which
 have type $A$ at time $s$ in $\{Y_t\}$ have no influence at all on the
 types of the $\rho^{(j)},\; 1 \le j \le m$, during $[0,s]$.   
 We can therefore construct the types of the finitely many $\rho^{(j)}$ in
 $\{Y_t(\la)\}$ by changing types appropriately at the only finitely
 many times during $[0,s]$ when one of these particles jumps to the
 position of another one, or when a recuperation event 
 $r(i,\rho^{(j)})$ occurs for some $j \le m$. It is not hard to check
that if $0 \le s_1 < s_2$, then the restriction of the process so
 constructed on $[0,s_2]$ to $[0,s_1]$ agrees with the process
 constructed on $[0,s_1]$. Indeed the only difference between 
the two constructions on $[0,s_1]$ could come from
the particles which have type $A$ at $s_1$, but type $B$ at
 $s_2$. However, these particles have not interacted with any particle
 during $[0,s_1]$. We shall not discuss the
 construction of the process $\{Y_t(\la)\}$ further, and in particular shall
 not verify that the above construction actually gives us a good
 version of $\{Y_t(\la)\}$.

 The preceding construction provides also a coupling of the processes
 $\{Y_t\}$ and $\{Y_t(\la)\}$. This coupling shows that $\{Y_t\}$ has
 more $B$-particles than the $\{Y_t(\la)\}$ process starting from
 $\si$, in the sense of the following lemma.
\proclaim{Lemma 2} Let $\{Y_t(\la)\}$ and $\{Y_t\}$ start at $\si$ and
 $\overline \si$, respectively, with $\overline \si \in \Si_0$. In
 particular, each
 particle is at the same position at time 0 in both processes and 
each particle which has type $B$ in $\{Y_t(\la)\}$ at time 0 also 
has type $B$ in $\{Y_t\}$ at time 0. Then the coupling 
 described above is such that
any particle present at a  space-time point
 $(x,s)$ in one of the processes $\{Y_t(\la)\}$ and $\{Y_t\}$  
is also present in the
 other. Moreover, if a particle at $(x,s)$ has type $B$ in 
$\{Y_t(\la)\}$,
 then it also has type $B$ in $\{Y_t\}$.
 \endproclaim

The lemma is immediate from the construction. The next lemma is very
 similar. It proves a monotonicity in the recuperation rate.
 \proclaim{Lemma 3} Let $0 \le \la_1 \le \la_2$ and let 
 $\{r_1(i,\rho)\}$ and $\{r_2(i,\rho)\}$ be Poisson processes 
 with the rates $\la_1$ and
 $\la_2$, respectively. Assume that these are coupled such that for
 each $\rho$ 
 $$
 \{r_1(i,\rho)\}_{i\ge 1}\subset \{r_2(i,\rho)\}_{i \ge 1}.
 \teq(2.1)
 $$
 Let $\{Y_t(\la_j)\}$ be the infection process corresponding to the 
recuperation rate $\la_j,\; j=1,2$, and assume that $\{Y_t(\la_1)\}$  
 and $\{Y_t(\la_2)\}$ are constructed from the same initial state
 $\si$ and the same set of random walk paths $\pi(\cdot,\rho)$, but potential 
 recuperation times $r_1(i,\rho)$ and $r_2(i,\rho)$, respectively.
 Assume that $\overline \si \in \Si_0$. Then the processes $\{Y_t(\la_1)\}$  
 and $\{Y_t(\la_2)\}$ are coupled in such a way
 that any particle present at a  space-time point
 $(x,s)$ in one of the $\{Y_t(\la_j)\}$ is also present in the
 other. Moreover, a.s. it holds for all $s$ that if a particle at 
$(x,s)$ has type $B$ in $\{Y_t(\la_2)\}$,
 then it also has type $B$ in $\{Y_t(\la_1)\}$.
 \endproclaim
 \demo{Proof} Clearly any particle $\rho$ present in
 one of the $\{Y_t(\la_j)\}$ at $(x,s)$ is also present 
 at $(x,s)$ in the
 other process since the position of any initial particle $\rho$ at
 time $s$ is $\pi(s,\rho)$ in both processes. 
\comment
Now let $\si$ be the (random) initial state of both $\{Y_t(\la_1)\}$
and $\{Y_t(\la_2)\}$ and let $\overline \si$ be as in Lemma 2. 
$\overline \si$ is formed by adding finitely many $B$-particles to 
the i.i.d. Poisson numbers $N_A(x,0-)$ of $A$-particles, and changing
all $A$-particles at the sites at which $B$-particles are added 
to $B$-particles. Then 
$\overline \si$ and $\si$ a.s. lie in $\Si_0$, by Proposition 5
and Lemma 15 of \cite {KSb}.
\endcomment

We can now couple the process $\{Y_t(\la_j)\}$ with a process
$\{Y_t\}$ which starts in $\overline \si$, as in Lemma 2. Then,
by Lemma 2,
the number of particles in $\{Y_t(\la_j)\}$ and in $\{Y_t\}$ at any
space-time point is the same, and the number of $B$-particles in
$\{Y_t(\la_j)\}$ is no more than in $\{Y_t\}$ at any space-time point.
This implies that a.s.,  for $j=1$ and for $j=2$, 
$$
 \sup_{s\le t} \big(\text{number of particles at 
 $(x,s)$ in  $\{Y_t(\la_j)\}\big) < \infty$ for all }x \in \Bbb Z^d, t \ge 0,
 \teq(2.08e)
 $$
and that there are only finitely many $B$-particles in 
$Y_t(\la_j)$ at any time $t$ (by virtue of Lemma 2 of \cite {KSb}).
\comment
As argued just before (2.33) in \cite {KSb} this implies that a.s.,
for $j=1,2$ 
$$
\min\{\wh \tau, \tau_\infty\} = \infty in \{Y_t(\la_j)\}.
\teq(2.08f)
$$
\endcomment
In particular, a.s. for all $s$, any $B$-particle at time
$s$ in $\{Y_t(\la_j)\}$ has an analogue of a genealogical path as above.

 Assume now that a particle $\rho$ has type $B$ at time $s$ in
 $\{Y_t(\la_2)\}$. Let its genealogical path in $\{Y_t(\la_2)\}$
be determined by the
 space-time points $(x_i,s_i)$ and by the particles $\rho_i$. 
That means that there  
 are space-time points $(x_i,s_i)$ with $1\le i \le \l$ for some $\l$,
 and $0 < s_1 < \dots < s_\l < s$ and particles $\rho_i$ for $0 \le i \le
 \l+1$ with $\rho_{\l+1} = \rho$, such that at time $s_i$, 
 $\rho_i$ jumps to the position of $\rho_{i+1}$ or vice versa. Moreover,
 $\rho_0$ has type $B$ at time 0, and $\rho_i$ does not recuperate
 during $[s_i, s_{i+1}]$ in $\{Y_t(\la_2)\}$ (with $s_{\l+1} = s$). 
 Note that, because $\rho_i$ stays of type $B$
 in $\{Y_t(\la_2)\}$ during $[s_i,s_{i+1}]$, 
 $r_2(j,\rho_i) \notin [s_i,s_{i+1}]$ for all $j$. But then $\rho_i$
 does not recuperate
 during $[s_i,s_{i+1}]$ in $\{Y_t(\la_1)\}$ either, by virtue of \equ(2.1). 
 It then follows by
 induction on $i$ that also in $\{Y_t(\la_1)\}$, 
 each $\rho_i$ is of type $B$ at time $s_i$
 and stays of type $B$ through time $s_{i+1}$. In particular
 $\rho=\rho_{\l+1}$ must have type $B$ at time $s$ in $\{Y_t(\la_1)\}$.
 \hfill $\blacksquare$
 \enddemo

\comment
 By taking $\la^{(1)} = 0$ we recover the following statement which
was a feature of our construction of $\{Y_t(\la)\}$:
 \proclaim{Corollary 3} The system without recuperation and the system
 with recuperation rate $\la > 0$ can be coupled in such a way that any
 particle at $(x,s)$ which has type $B$ in the system
 with recuperation is also a particle at $(x,s)$ 
 of type $B$ in the system without recuperation. 
 \endproclaim
\endcomment

 A consequence of Lemma 3 is that if the infection dies out for
 some value $\la^{(1)}$ of the recuperation rate, then it dies out for all
 larger recuperation rates. As already stated this shows that $\la_c$
 exists, but it may still have the value  0 or $\infty$.

 We will also need another monotonicity property for $\{Y_t(\la)\}$.
 Basically this says that if we increase the number of $B$-particles in
 the initial state, then this will increase the number of $B$-particles
 at any later time. The analogue of this result for $\{Y_t\}$ is 
in lemma 14 of \cite {KSb}.
\proclaim{Lemma 4} Let $\la \ge 0$ and let
$\si^{(2)}$ be such that $\overline \si^{(2)} \in
\Si_0$. Assume further that $\si^{(1)}$ lies below $\si^{(2)}$ in the following
sense:
$$
\text{for any site $x \in \Bbb Z^d$, all particles present in
$\si^{(1)}$ at $x$ are also present in $\si^{(2)}$ at $x$},
\teq(7.1a)
$$
and
$$
\text{any particle which has type $B$ in $\si^{(1)}$ 
also has type $B$ in $\si^{(2)}$}.
\teq(7.2a)
$$
Let $\pi_A(\cdot,\rho)$
be the random walk paths associated to the 
 various particles. Assume that the Markov processes
 $\{Y_t^{(1)}(\la)\}$ and $\{Y_t^{(2)}(\la)\}$ are constructed 
{\rom (}as explained in Section 2{\rom)} by means of
 the same set of paths $\pi_A(\cdot,\rho)$
and the same recuperation processes $\{r(1,\rho)\}$ for any $\rho$
present in $\si^{(1)}$. Assume further that $\{Y_t^{(i)}(\la)\}$
starts in $\si^{(i)},\; i=1,2$. Then, almost surely, 
$\{Y_t^{(1)}(\la)\}$ and $\{Y^{(2)}_t(\la)\}$ satisfy 
 \equ(7.1a) and \equ(7.2a) for all $t$,
with $\si^{(i)}$ replaced by $Y^{(i)}_t(\la),\;
 i=1,2$. Moreover, almost surely
  $$
 \sup_{s\le t} \big(\text{number of particles at 
 $(x,s)$ in  $\{Y^{(i)}_t(0)\}\big) < \infty$ for all }x \in \Bbb Z^d, t \ge 0,
 \teq(2.08g)
 $$
 for $i=1$ and for $i=2$. 
\endproclaim

\demo{Proof} It is clear that \equ(7.1a) holds with $\si^{(i)}$
replaced by $Y_t^{(i)}(\la)$, that is,
$$
\align
&\text{for any site $x \in \Bbb Z^d$, and $t \ge 0$, all particles present in
$Y_t^{(1)}$ at $x$}\\
&\phantom{MMMMMMMMM}\text{are also present in $Y_t^{(2)}$ at $x$}.
\teq(7.1az)
\endalign
$$

By Lemma 14 in \cite {KSb} $\overline \si^{(2)} \in \Si_0$ implies
that also $\overline \si^{(1)} \in \Si_0$. In the same way 
as in the second paragraph
of the proof of Lemma 3 one now shows that a.s. \equ(2.08g) holds and that
a.s. there are only finitely many $B$-particles in 
$Y^{(i)}_t(\la)$ at any time $t$. Also 
a.s. for all $s$ any $B$-particle at time
$s$ in $\{Y_t^{(i)}(\la)\}$ has an analogue of a genealogical path.

To prove \equ(7.2a) with $\si^{(i)}$ replaced by $Y_s^{(i)}(\la)$,
assume that $\rho$ has type $B$ at time $s$ in the
first process, i.e., in $\{Y_t^{(1)}(\la)\}$. Then it has a genealogical
path determined by space-time points 
$(x_j,s_j)_{1\le j \le \l}$ for some $\l$,
 and $0 < s_1 < \dots < s_\l < s$ and particles $\rho_j$ for $0 \le j \le
 \l+1$ with $\rho_{\l+1} = \rho$ and $s_{\l+1} = s$, such that at time $s_j$, 
 $\rho_j$ jumps to the position of $\rho_{j+1}$ or vice versa. Moreover,
all these $\rho_j$ and $\rho$ are present in $\si^{(1)}$ (and hence
are particles in $\{Y_t^{(1)}\}$), 
$\rho_0$ has type $B$ at time 0, 
and $\rho_j$ has type $B$ and does not recuperate
during $[s_i, s_{i+1}]$ in $\{Y_t^{(1)}(\la)\}$ (with $s_{\l+1} = s$). 
One then proves by induction on $j$ that each $\rho_j,\; 0 \le j \le
\l+1$, is also present and has type $B$ during $[s_j,s_{j+1}]$
in $\{Y^{(2)}_t\}$. In particular, $\rho = \rho_{\l+1}$ is present and
of type $B$ at time $s$ in $\{Y_t^{(2)}\}$. Thus, \equ(7.2a) holds.
\hfill $\blacksquare$
\enddemo
\comment
\proclaim{Lemma 4} Let $\si^{(2)} \in
 \Si_0$. Assume further that $\si^{(1)}$ lies below $\si^{(2)}$ 
in the following
 sense:
 $$
 \text{for any site $z \in \Bbb Z^d$, all particles present in
 $\si^{(1)}$ at $z$ are also present in $\si^{(2)}$ at $z$},
 \teq(7.1a)
 $$
 and
 $$
 \align
 &\text{at any site $z$ at which the particles in $\si^{(2)}$ have type $A$},\\
 &\text{the particles also have type $A$ in $\si^{(1)}$}.
 \teq(7.2a)
 \endalign
 $$
 Let $\pi(\cdot,\rho)$ be the random walk paths associated to the 
 various particles and assume that the Markov processes
 $\{Y_t^{(1)}\}$ and $\{Y_t^{(2)}\}$ are constructed by means of
 the same set of paths $\pi(\cdot,\rho)$ and
 starting with state $\si^{(1)}$ and $\si^{(2)}$, respectively 
 {\rom (}as defined in Section 2 of \cite {KSb}; see in particular 
 (2.6) and (2.7) of \cite
 {KSb} with $\pi_A(\cdot,\rho) = \pi_B(\cdot,\rho) =
 \pi(\cdot,\rho)${\rom )}. Then $\{Y_t^{(1)}\}$ and
 $\{Y^{(2)}_t\}$ satisfy almost surely for all $t$
 \equ(7.1a) and \equ(7.2a) with $\si^{(i)}$ replaced by $Y^{(i)}_t,\;
 i=1,2$. In particular, $\si^{(1)} \in \Si_0$ and almost surely
 $$
 \sup_{s\le t} \big(\text{number of particles at 
 $(z,s)$ in  $\{Y^{(i)}_t\}\big) < \infty$ for all }z \in \Bbb Z^d, t \ge 0,
 \teq(2.08e)
 $$
 for $i=1$ {\rom (}as well as for $i=2${\rom )}.
 \endproclaim
\endcomment

 \subhead
 3. Survival for small $\la$
 \endsubhead
 \numsec=3
 \numfor=1
 In this section we show that $0 < \la_c \le \infty$. 
 To introduce the directed percolation process which we promised in
 the introduction, we must describe certain
 blocks in $\Bbb Z^{d+1}$. $C_0$  will be the same large integer as in
 \cite {KSb} (see (4.18), (4.19) there). Without loss of generality we
 take $C_0$ even.
Also $\ga_0 \in (0,\infty)$ will be as in \cite {KSb}.
 Many constants $K_i$ and $p_i$ will appear in the
 proof. These will all depend only on $d, D, C_0,\ga_0, \mu_A$. 
All $K_i$ and $p_i$ are finite
 and strictly positive. These properties of the $K_i,p_i$ 
 will not be mentioned further. Throughout this section we think of $p$
 as fixed, and often suppress it in the notation; 
 we shall see at the end of the proof of Lemma 12 that any large
 enough value of $p$ will work for our purposes.
For the time being we only need to know that $p$ is an integer $\ge 1$.
We also fix 
 $$
 q = 2d+1
 $$ 
 and define
 $$
 \De_r = C_0^{6r}.
 $$
 For $\bold i=(i(1), \dots, i(d)) \in \Bbb Z^d$ and $k \in \Bbb Z$ 
 we take 
 $$
 \wh \cB_p(\bold i,k) = \prod_{s=1}^d [i(s)\De_p,
 (i(s)+1)\De_p) \times [kp^q\De_p,(k+1)p^q\De_p).
 \teq(3.1)
 $$
 This definition is similar to that of the blocks $\cB_r(\bold i,k)$
 used in \cite {KSa}-\cite {KSc}, but there are obvious differences  in
 the handling of the last coordinate in these definitions. We further
 define the {\it bottom} of the block $\wh \cB_p(\bold i,k)$ as
 $$
 \cZ_p(\bold i,k) = \prod_{s=1}^d[(i(s)-4d-1)\De_p, (i(s)+4d+2)\De_p) \times
 \{kp^q\De_p\}.
 \teq(3.2)
 $$
 The directed graph $\cD$ will be the graph with vertex set $\Bbb Z^d
 \times  \{-1,0,1,2,\dots\}$,  
 and with a directed edge from $(\bold i,k)$ to $(\bold
 j,\l)$ if and only if $\|\bold i-\bold j\| \le 1$ and  $\l = k+1$. 
(Recall that
 the first condition means $|i(s)-j(s)| \le 1$ for $1 \le s \le d$.) 
 We also need the graph
 $\cL$. It has vertex set $\Bbb Z^{d+1}$ and an edge between $v$
 and $w$ if and only if $\|v-w\| \le 1,
 v \ne w$. We shall call the edges of $\cD$ and $\cL$,
 $\cD$-edges and $\cL$-edges, respectively.
 We shall call $(\bold i,k)$ a {\it parent} of $(\bold j,k+1)$ if 
 there is a $\cD$-edge from $(\bold i,k)$ to $(\bold j,k+1)$.

 For any set $A$ in the vertex set of $\cL$ (i.e., $ A \subset \Bbb  Z^{d+1}$)
 we define the following pieces of its boundary:
 $$
 \align
 \partial_{ext}A = \{v \in \Bbb Z^{d+1}:\;& v \text{ is adjacent on $\cL$ 
 to some $w \in A, v \notin A$,
 and there}\\
 &\text{exists a path on $\Bbb
 Z^{d+1}$ from $v$ to $\infty$ which avoids }A\};
 \endalign
 $$
 $$
 \align
 \partial_{ext}^+ A:= \{v \in \partial_{ext}: &\text{ there is some $w \in A$
 such that the edge}\\
 &\text{ from $w$ to $v$ is a $\cD$-edge}\};
 \endalign
 $$
 $$
 \partial_{ext}^* A:= \{v \in \partial_{ext}:v+e_{d+1} \in A\}.
 $$
 Note that $\partial_{ext}^+$ and $\partial_{ext}^*$ are not
 disjoint in general. If $A,S \subset \Bbb Z^{d+1}$, then we say that
 $S$ {\it separates }$A$ {\it from} $\infty$ {\it on} $\Bbb Z^{d+1}$
 if $S \cap A = \emptyset$ and 
every path on $\Bbb Z^{d+1}$ from $A$ to $\infty$ contains a point
 of $S$.

 The next lemma is of a 
 topological nature only. 
 \proclaim{Lemma 5} Let $A\subset \Bbb Z^d \times \{0, 1, 2 \dots\}$ 
 be a finite, non-empty, $\cL$-connected set. Then
 $$
 \partial_{ext} A \text{ is $\Bbb Z^{d+1}$-connected and 
 separates $A$ from $\infty$ on $\Bbb Z^{d+1}$}
 \teq(3.3)
 $$
 and
 $$
 |\partial_{ext} A| \le 6|\partial_{ext}^+ A|.
 \teq(3.4)
 $$
 \endproclaim
 \demo{Proof} Relation \equ(3.3) is just a special case of
 lemma 2.23 in \cite {Kb} (with $d$ replaced by $d+1$).
\cite {Kb} does not state the fact that $\partial_{ext} A$ separates
 $A$ from $\infty$ in the generality of the present lemma.
However, the proof on the top of p. 144 of \cite {Kb} shows easily that
the separation property in \equ(3.3) holds.

 To prove \equ(3.4), assume that $v \in \partial_{ext}A$. 
 Then $v$ is adjacent on $\cL$ to some $w \in A$
 and there exists some path $\pi$ from $v$ to $\infty$ on $\Bbb
 Z^{d+1}$ which
 is disjoint from $A$. We distinguish three main cases according to the
 value of $v(d+1)-w(d+1)$ ($v(d+1)$ is the last coordinate of
 $v$); the last two
 cases are split into two subcases.
 \newline
 Case a): $ v(d+1) = w(d+1) + 1$. In this case the edge from $v$ to $w$ 
 is a $\cD$-edge, so that
 $v \in \partial_{ext}^+A$. Thus the number of vertices $v \in
 \partial_{ext} A$ which are in case a) is at most $|\partial_{ext}^+|$. 
 \newline  
 Case b): $v(d+1) = w(d+1)$. Subcase bi): $v + e_{d+1} \notin A$. Here we
 abuse notation somewhat. $e_j$ denotes the $j$-th coordinate vector and
 the $s$-th component of $v+e_{d+1}$  equals $v(s)$ if $s \le d$
 and equals $v(d+1)+1$ if $s = d+1$. In this subcase, the path on $\Bbb
 Z^{d+1}$ consisting of the edge from $v+e_{d+1}$ to $v$ followed by
 $\pi$ is a path on $\Bbb Z^{d+1}$ 
 from $v+e_{d+1}$ to $\infty$ which is disjoint from $A$. 
 Moreover, $(v + e_{d+1})(d+1) = w(d+1)+1$ and the edge from $w$ to $v+e_{d+1}$
 is a $\cD$-edge.
 Thus $v+e_{d+1} \in \partial_{ext}^+ A$ and again, the number of
 vertices $v \in \partial_{ext} A$ which are in case bi) 
 is at most $|\partial_{ext}^+|$.
 Subcase bii) applies when $v + e_{d+1} \in A$. 
 Then the edge from $v+e_{d+1}$ to $v$ goes from a point of $A$  
 to a point of $\partial_{ext}A$, but the last coordinate decreases
 by one along this edge. Thus, $v \in \partial_{ext}^* A$ in this case.
 Thus the number of vertices $v \in
 \partial_{ext} A$ which are in case bii) is at most
 $|\partial_{ext}^*|$. To complete the handling of this subcase we
 prove that in general
 $$
 |\partial_{ext}^* A| \le |\partial_{ext}^+ A|
 \teq(3.5)
 $$ 
 for any finite $A \subset \Bbb Z^d \times \{0,1, \dots\}$. To see \equ(3.5)
 consider any line parallel to the last coordinate axis of the form
 $\{v_0+ ne_{d+1}: n \in \Bbb Z\}$. The points of this line are in
 the unbounded component of $\Bbb Z^{d+1} \setminus A$ for large $n$
 both in the positive and negative direction.
 Therefore, as one lets $n$ run from $-\infty$ to $+\infty$,
 there are as many transitions from the unbounded component in $\Bbb
 Z^{d+1}$ of $\Bbb Z^{d+1} \setminus A$ 
 to $A$ as there are transitions from $A$ to
 the unbounded component of $\Bbb Z^{d+1} \setminus A$. The former
 transitions go from a vertex $v$ outside $A$ to a vertex in $A$ by
 adding $e_{d+1}$, and therefore occur for $v \in \partial_{ext}^*$. 
 The latter transitions are along a $\cD$-edge from a vertex of $A$ to
 a vertex $v$ outside $A$ and therefore occur when $v+e_{d+1} 
 \in \partial_{ext}^+$. The numbers of the two types of transitions are
 equal, and this holds for any choice of $v_0$. \equ(3.5) follows.
 \newline
 Case c): $v(d+1) = w(d+1)-1$. Again this has the subcases ci) with
 $v+e_{d+1} \notin A$ and cii) with $v+e_{d+1} \in A$.  
 In case ci) one easily checks (by the argument for case bi) 
 that $\wt v :=v+e_{d+1} \in \partial_{ext} A$, and that $\wt v$ is in
 case b). Thus, by the results for case b) 
  the number of
 vertices $v \in \partial_{ext} A$ which are in case ci) 
 is at most $2|\partial_{ext}^+|$.
 Finally, if v is in subcase cii), then replace $w$ by $\wt w =
 v+e_{d+1}$. In
 this situation, $v$ is adjacent on $\cL$ to $\wt w \in A$  and
 therefore $v$ lies in $\partial_{ext}^*$. \equ(3.5) therefore shows
 that also
 the number of vertices $v \in
 \partial_{ext} A$ which are in case cii) is at most $|\partial_{ext}^+|$.
 The inequality \equ(3.4) follows by adding the contributions of the
 various cases.
 \hfill $\blacksquare$
 \enddemo

 We can now set up our percolation problem on the graph $\cD$. We 
 define $m(\bold i) = m_p(\bold i) \in (\Bbb Z+\frac 12)^d \De_p$ as the
 point with components 
 $$
  m(\bold i)(s) = (i(s)+1/2)\De_p,\;1 \le s \le d.
 $$
 $m(\bold i)$ is in some sense the midpoint of $\prod_{s=1}^d
 [(i(s)\De_p, (i(s)+1)\De_p)$, which constitutes the spatial part of
 $\wh \cB_p(\bold i, k)$. $m(\bold i)(s)$ is an integer because we
 took $C_0$ even. For purposes of the proof of survival of the
 infection, it turns out to be convenient to change the initial conditions
 of the $B$-particles slightly. {\it For the rest of this section 
 we will assume that we do
 not add a $B$-particle at the origin at time 0, but instead add
 a $B$-particle at $m(\bold 0)$}. 
 Thus we take the state at time 0 to satisfy
 $$
 N_A(x,0) = N_A(x,0-)  \text{ if } x \ne m(\bold 0),
 \teq(3.11)
 $$
 $$
 N_A(x,0) + N_B(x,0) = N_A(x,0-)+1 \text{ if } x = m(\bold 0),
 \teq(3.11a)
 $$
 $$
 N_B(x,0) =  0  \text{ if } x \ne m(\bold 0),
 \teq(3.12)
 $$
 and
 $$
 1 \le N_B(x,0) \le N_A\big(m(\bold 0),0-\big) +1 \text{ if } x = m(\bold 0).
 \teq(3.12b)
 $$
 Clearly \equ(1.1a) holds with the original initial condition if and only
 if it holds in this modified system.
 Thus it suffices for showing $\la_c > 0$ that 
 $$
 \align
 P\{\text{there are $B$-particles} &\text{  at all times in the system}\\
 &\text{which starts with \equ(3.11)-\equ(3.12b)}\} > 0.
 \teq(3.12a)
 \endalign
 $$  
 It will be necessary in the proofs of Lemma 6 and 7 to consider initial
 conditions in which a $B$-particle is added at time 0 at a finite
 number of sites $m(\bold c_1),\dots,m(\bold c_r)$. In this situation
 $m(\bold 0)$ in \equ(3.11)-\equ(3.12b) has to be replaced by
 $m(\bold c_1), \dots, m(\bold c_r)$.
 Till the end of Lemma 7 we shall allow this, but will
 indicate the location of the initial particles in the notation only
 where it is crucial.

 We further define
 $$
 t(k) = t_p(k) = kp^q \De_p,
 \teq(3.12abc)
 $$
 and 
 $$
 Z_p(\bold i) = \prod_{s=1}^d[(i(s)-4d-1)\De_p, (i(s)+4d+2)\De_p) 
 \subset \Bbb Z^d,
 \teq(3.12cde)
 $$
so that $\cZ_p(\bold i,k) = Z_p(\bold i) \times \{t(k)\}$.
We also define
 $x(\bold i,k) \in \Bbb Z^d$ as the nearest (in the $\l^\infty$ sense 
on $\Bbb Z^d$) 
 site to $m(\bold i)$ which contains a
 $B$-particle at time $t(k)=kp^q\De_p$ in our infection process
 $\{Y_t(\la)\}$. If there are several possible choices for $x(\bold
 i,k)$, then we use some deterministic rule to break the tie. If there
 are no $B$-particles in $\{Y_t(\la)\}$ at time $t_k$, then we leave
 $x(\bold i,k)$ undefined.
If a $B$-particle is added at $m(\bold c)$ at
 time 0, then we take $x(\bold c,0) = m(\bold c)$.
 We call
 the vertex $(\bold i,k)$  of $\cD$ {\it active} 
(or more explicitly $\la$-{\it active}) if there is a  
 site $x \in m(\bold i) + \cC(\frac 18\De_p)$ which is occupied by at
 least one
 $B$-particle at time $t(k)$ in our infection process with
 recuperation $\{Y_t(\la)\}$
 (see \equ(2.0a) for $\cC$). By convention, 
 if a $B$-particle is added at $m(\bold c)$ at
 time 0, the vertex $(\bold c,0)$ is active.

 We now want to define when certain $\cD$-edges are open. 
 To this end we first define
 the $\cZ_p(\bold i,k)$-{\it process started at }
 $\big(x,t(k)\big)$ for any $x \in Z_p(\bold i)$.
 This process is defined only
 from time $t(k)$ on and it will use only particles which are in 
 $Z_p(\bold i)$ at time $t(k)$.
  Also, we only define this process if $x$ is occupied
 by some particle at time $t(k)$. 
 To define this process we first reset
 the types of the particles in $Z_p(\bold i)$ at time $t(k)$. All
 particles in $Z_p(\bold i)\setminus \{x\}$ 
are given type $A$. One particle at $x$ is given
 type $B$. Denote this particle by $\rho(x,t(k))$. 
All other particles at $x$ (if any) are given type
 $A$. If there are $B$-particles at $\big(x,t(k)\big)$, 
then $\rho(x,t(k))$ is chosen from these $B$-particles, but
 apart from this restriction $\rho(x,t(k))$ 
can be selected from the particles at $(x,t(k))$ in any way which does 
 not depend on the future paths of the particles in $Z_p(\bold i)$ 
at time $t(k)$. 
 The $\cZ_p(\bold i,k)$-process started at $\big(x,t(k)\big)$ is then 
 the evolution of the particles which are in $Z_p(\bold i)$ at 
time $\{t(k)\}$ with
 the reset types according to the rules for $\{Y_t(0)\}$, that is,
 there is {\it no recuperation}, but we still insist that an
 $A$-particle turns into a $B$-particle only if it jumps onto a
 $B$-particle or a $B$-particle jumps onto it.
 Note that in this process all particles outside
 $Z_p(\bold i)$ at time $t(k)$ are ignored.

 We now say that the $\cD$-edge from $(\bold i,k)$ to $(\bold j,k+1)$ 
 is {\it open}
 if the following three events (3.13)-(3.15) occur.
 $$
 (\bold i,k) \text{ is active}.
 \teq(3.7)
 $$
 $$
 \align
 A(\bold i,k,\bold j):=
 \{&\text{the $\cZ_p(\bold i,k)$-process started at $\big(x(\bold i,k),
 t(k)\big)$ has at least}\\
 &\text{one $B$-particle in $m(\bold j) +\cC(\frac 18
 \De_p)$ at time $t(k+1)$}\}
 \teq(3.9)
 \endalign
 $$
(see Figure 1).
 If $A(\bold i, k, \bold j)$ occurs, then there exists in $\{Y_t(0)\}$
 a genealogical
 path from some $B$-particle at $\big( x(\bold i,k), t(k)\big)$ to some
 particle in $m(\bold j) + \cC(\frac 18 \De_p)$. Among all such paths
 choose the first one in some deterministic ordering of such paths. Let
 this be determined by the times $s_i,\; 1 \le i \le \l$, 
 and particles $\rho_0, \dots,\rho_\l$, in the sense that $t(k) < s_1<
 \dots < s_\l < t(k+1)$, $\rho_0$ is some $B$-particle at  
 $\big( x(\bold i,k), t(k)\big)$ and $\rho_\l$ is located in
 $m(\bold j) + \cC(\frac 18 \De_p)$ at time $t(k+1)$; moreover, at time
 $s_i, 1 \le i \le \l$, one of $\rho_i$ and $\rho_{i-1}$ jumps to the
 position of the other. All the particles $\rho_i,\; 0 \le i \le \l$,
 are in $Z_p(\bold i)$ at time $t(k)$. Note that by our definition of
 the $\cZ_p(\bold i,k)$-process, the path here is chosen without
 reference to the recuperation events. 
 The last required event for the edge from $(\bold i,k)$ to $(\bold
 j,k+1)$ to be open is 
 $$
 \align
 B(\bold i,k,\bold j,\la) := \;& \{\text{with $s_i$ and $\rho_i$ as in the
 preceding lines, the particle $\rho_i$}\\
 &\text{has no recuperation event in
 $\{Y_t(\la)\}$ during $[s_i,s_{i+1}]$,}\\
 &\text{that is, $r(h,\rho_i) \notin
 [s_i,s_{i+1}]$ for all $h$ and $0 \le i \le \l$}\\
 &\text{(with $s_{\l+1} = t(k+1))$}\}.
 \teq(3.8)
 \endalign  
 $$

\epsfverbosetrue
\epsfxsize=300pt
\epsfysize=150pt
\centerline {\epsfbox{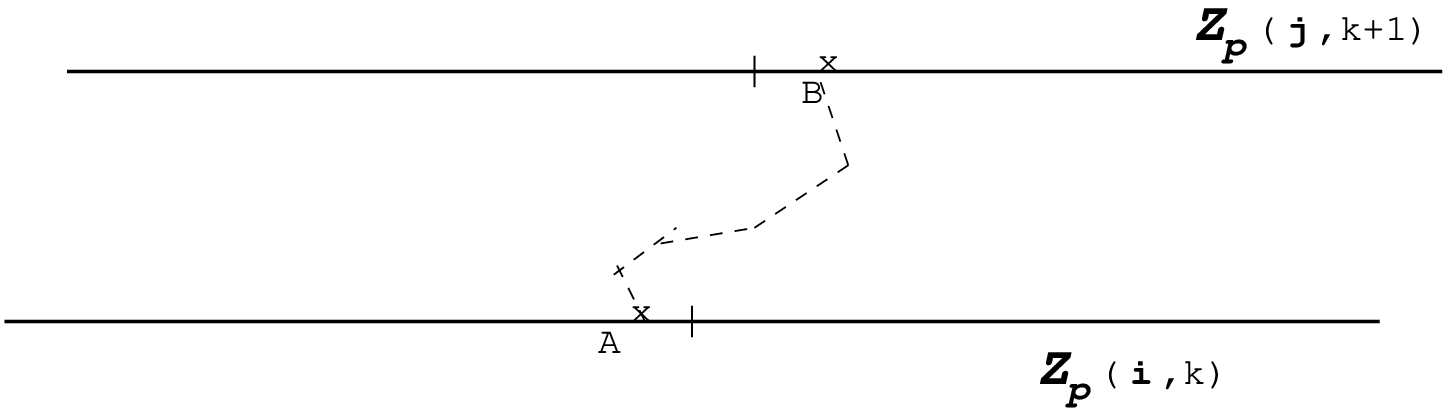}}

\botcaption
{Figure 1} 
Relative location of the sets $\cZ_p(\bold i,k), \cZ_p(\bold j,k+1)$
for $d=1$, where $(\bold i,k)$ is a parent of $(\bold j,k+1)$. 
The points marked by a small
vertical bar are $\big(m(\bold i),t(k)\big)$ and $\big(m(\bold j),
t(k+1)\big)$, respectively.
The point $A$, marked by an x is $\big(x(\bold i,k), t(k)\big)$. 
$B$, also marked by an x, is the endpoint in $\big(m(\bold j)+ \cC(\frac
18 \De_p), t(k+1)\big)$ 
of the dashed curve. This curve represents 
the genealogical path along which the infection is transmitted in the
$\cZ_p(\bold i,k)$-process started at $A$, so that $A(\bold i,k,\bold
j)$ occurs.
\endcaption
\bigskip

 Note that this definition applies only if there exists a $\cD$-edge from
 $(\bold i,k)$ to $(\bold j, k+1)$, that is if $\bold j = \bold i$ or
 $\bold j= \bold i \pm e_s$ for some $s\in \{1, \dots, d\}$. If any of 
 the three conditions \equ(3.7)-\equ(3.8) fails, then the $\cD$-edge from
 $(\bold i,k)$ to $(\bold j, k+1)$ is called {\it closed}.

 By definition of ``active'', the
 infection survives if there are with positive probability infinitely
 many active sites. The next lemma is a tool for finding active sites. 

 \proclaim{Lemma 6} Start the infection process off by adding
 a $B$-particle at $m(\bold c_1),\dots, m(\bold c_r)$ at time 0.
 If some vertex $(\bold i,k)$ is active and the $\cD$-edge from 
 $(\bold i,k)$
 to $(\bold j,k+1)$ is open, then $(\bold j,k+1)$ is also active. 
 \endproclaim
 \demo{Proof} 
 By assumption \equ(3.7), 
 $$
 x(\bold i,k) \in m(\bold i)+\cC(\frac 18\De_p) 
 \subset Z_p(\bold i)
 $$ 
 (this is even true for $(\bold i,k) = (\bold c_j,0)$, since we
 interpret $x(\bold c_j,0)$ as $m(\bold c_j)$). Now apply Lemma 4 with
 $\si^{(2)}$  the true state of the process $\{Y_t(\la)\}$ at time $t(k)$
 and $\si^{(1)}$ the state obtained by resetting the
 type of the particles in $\cZ_p(\bold i,k)$ and ignoring the particles
 outside $\cZ_p(\bold i,k)$ to form the
 $\cZ_p(\bold i,k)$-process started at $\big(x(\bold i, k),t(k)\big)$.
 This procedure only involves removing particles and changing the type
 of some particles from $B$ to $A$. Indeed, according to our
 construction, the one particle at $x(\bold i, k)$
 which is given type $B$ at $\big(x(\bold i, k),t(k)\big)$ 
 is chosen from the particles which already had type $B$ at time $t(k)$
 in $\{Y_t(\la)\}$.
  Therefore, 
 \equ(7.1a) and \equ(7.2a) hold for these choices of $\si^{(1)}$ and
 $\si^{(2)}$. If we now let the particles in the full system (i.e.,
 the collection of all the particles, including the ones outside
 $\cZ_p(\bold i,k)$) develop
 till time $t(k+1)$ according to the rules for $\{Y_t(0)\}$, 
then Lemma 4 tells us that
 at time $t(k+1)$  the $\cZ_p(\bold
 i,k)$-process started at $\big(x(\bold i,k),t(k)\big)$ will still lie
 below the full process.
 Moreover, by assumption \equ(3.9), the $\cZ_p(\bold
 i,k)$-process started at $\big(x(\bold i,k),t(k)\big)$ has at least
 one $B$-particle in $m(\bold j) + \cC(\frac 18\De_p)$ at time
 $t(k+1)$.

 The last few lines tell us that there will be some
 $B$-particles in $m(\bold j)+\cC(\frac 18 \De_p)$
at time $t(k+1)$ in  the full process, if it develops
according to the rules for $\{Y_t(0)\}$ (i.e., without recuperation) 
during $[t(k),t(k+1)]$. In fact we can choose
 times $s_i$ and particles $\rho_i$ as in the lines between \equ(3.9)
 and \equ(3.8). Then $\rho_i$ has type $B$ during $[s_i,s_{i+1}], 0 \le
 i \le \l$, and in particular $\rho_\l$ will have type $B$ at $t(k+1)$,
 if we suppress recuperation during $[t(k),t(k+1)]$.
 Finally we note that by induction on $i$, 
 the occurrence of $B(\bold i,k, \bold j,\la)$ implies that
 $\rho_i$ still has type $B$ during $[s_i, s_{i+1}]$
 even if recuperation is allowed. Indeed, if $\rho_i$ is of type $B$
 at time $s_i$ in the process $\{Y_t(\la)\}$, then it will keep type
 $B$ during $[s_i,s_{i+1}]$, even when recuperation is suppressed, 
since it has no recuperation during this
 interval anyway if $B(\bold i,k,\bold j,\la)$ occurs. Then $\rho_i$ and  
 $\rho_{i+1}$ coincide at time $s_{i+1}$, and therefore $\rho_{i+1}$
 will be of type $B$ at $s_{i+1}$ even in $\{Y_t(\la)\}$.
 In particular $\rho_{\l+1}$ has type $B$ at time $t(k+1)$.

 We conclude from the last two paragraphs that even in the process
 $\{Y_t(\la)\}$ there is a $B$ particle in 
 $m(\bold j) + \cC(\frac 18\De_p)$ at time
 $t(k+1)$.
 In particular, $x\big (\bold j,t(k+1)\big)$, the 
 position of the nearest $B$-particle to 
 $m(\bold j)$ at time $t(k+1)$ in $\{Y_t(\la)\}$, must lie in 
 $m(\bold j) + \cC(\frac 18\De_p)$.
 \hfill $\blacksquare$
 \enddemo

 Again assume that we add a $B$-particle at $m(\bold c_1), \dots,
 m(\bold c_r)$ at time 0. Now define the open cluster of the set
 $\{(\bold c_1,0),\dots, (\bold c_r,0)\}$  on $\cD$ as
 $$
 \align
 \frak C = \frak C(\bold c_1, \dots, \bold c_r) =
 \;\{&v \in \cD: \exists \text{ path $v_0, v_1, \dots, v_n=v$
 from some $v_0 = (\bold c_j,0)$}\\
 &\text{to $v$ on $\cD$ for some $n$ such that each edge
 $\{v_i,v_{i+1}\}$,}\\
 & 0 \le i \le n-1, \text{ is open}\}.
 \teq(3.13)
 \endalign
 $$
 We always include each $(\bold c_j,0)$ in $\frak C(\bold c_1, \dots,
 \bold c_r)$. From Lemma 6 it follows
 that each vertex $(\bold i,k)$ in $\frak C$ must be active.
 Thus, if $\frak C$ is infinite with positive probability, then the
 infection survives. If $\frak C$
 is finite (and nonempty, since it contains $(\bold c_j,0)$), 
then $\partial_{ext}\frak C$  
 is $\Bbb Z^{d+1}$-connected and separates $\frak C$ from $\infty$ on
 $\Bbb Z^{d+1}$, by Lemma 5. 
 Moreover, each $\cD$-edge from some vertex $(\bold
 i,k) \in \frak C$ to some vertex $(\bold j, k+1) \in \partial_{ext}^+
 \frak C$ 
 must be closed. This is true by definition, because if $(\bold
 i,k) \in \frak C$ and the edge from $(\bold i,k)$ to $(\bold j, k+1)$
 were open, then also $(\bold j,k+1)$ would belong to $\frak C$.
 These observations indicate that it will be useful to have an upper
 bound for the probability that the edge from $(\bold i, k)$ to $(\bold
 j, k+1)$ is closed. To derive such a bound we generalize the
 definitions of the events $A$ and $B$ in \equ(3.9) and \equ(3.8).
 For $x \in m(\bold i)+\cC(\frac 18\De_p)$ we define
 $$
 \align
 \wt A\big(x,t(k), \bold j\big) := \big\{&x \text{ is not occupied 
 at time $t(k)$}\big\} \cup \big\{\text{$x$ is occupied at $t(k)$ }\\
 &\text{and the $\cZ_p(\bold i,k)$-process started at }
 \big(x,t(k)\big)
 \text{ has at least}\\
 &\text{one $B$-particle in }m(\bold j)+\cC(\frac
 18\De_p) \text{ at time }t(k+1)\big\}.
 \teq(3.13x)
 \endalign
 $$
 If $\wt A\big(x,t(k),\bold j)$ occurs and $x$ is occupied at time
 $t(k)$, then we define $\wt B(x,t(k), \bold j, \la)$ as
 follows. First we reset the  types of the
 particles in $\cZ_p(\bold i,k)$ to form the $\cZ_p(\bold i,k)$-process
 started at $\big(x,t(k)\big)$. We then 
 pick a genealogical
 path in this process from the unique $B$-particle (after resetting) at 
 $\big(x,t(k)\big)$ to a vertex in $m(\bold j)+\cC(\frac 18\De_p)$ at
 time $t(k+1)$. 
 As in the lines between \equ(3.9) and \equ(3.8) let this path be
 determined by times $s_1, \dots, s_\l$ and particles $\rho_0, \dots,
 \rho_\l$, which have to be in $Z_p(\bold i)$ at time $t(k)$.
 We then also define (with $s_{\l+1} = t(k+1)$)
 $$
 \align
 \wt B\big(x,t(k),\bold j,\la) := \big\{&x \text{ is not occupied 
 at time $t(k)$}\big\} \cup \big\{\big(x,t(k)\big) 
 \text{ is occupied},\\
 &\text{and for $0 \le i \le \l$ the particle $\rho_i$
 has no recuperation}\\
 &\text{event in
 $\{Y_t(\la)\}$ during $[s_i,s_{i+1}]$}\big\}.
 \teq(3.13y)
 \endalign
 $$
By these definitions, on the event $\{(i,k)$ is active\},
$A(\bold i,k,\bold j) =
 \wt A\big(x(\bold i,k),t(k),\bold j\big)$ and $B\big(\bold i,k,\bold j,
 \la)= \wt B\big(x(\bold i,k),t(k), \bold j,\la)$.
 Finally,  if there exists a $\cD$-edge from $(\bold i,k)$
 to $(\bold j,k+1)$ and if $(\bold i, k)$ is active, we set
 $$
 C(\bold i,k,\bold j,\la) := 
 \bigcup_{x \in m(\bold i) + \cC(\frac 18\De_p)} \Big[\wt A\big(x,t(k), 
 \bold j\big) \cap \wt B(x,t(k),\bold j,\la) \Big]^c.
 \teq(3.13z)
 $$
 If there is no $\cD$-edge from $(\bold i,k)$
 to $(\bold j,k+1)$, or $(\bold i,k)$ is not active, then we set
 $$
 C(\bold i,k,\bold j,\la) := \emptyset.
 $$
 Finally,  let $\frak C_0$ be some finite $\cD$-connected subset of
 vertices of $\cD$. The vertices of $\cD$ can also be thought of as vertices 
 of $\Bbb Z^{d+1}$, so that we can also think of $\frak C_0$ as a
 subset of $\Bbb Z^{d+1}$. We shall continue to denote generic vertices
 of $\Bbb Z^{d+1}$ as $(\bold i,k)$ with $\bold i \in \Bbb Z^d$ and $k
 \in \Bbb Z$, because the last coordinate always plays the 
special role of time.
 We call a subset $S$ of $\Bbb Z^{d+1}$ a 
 $\frak C_0${\it -barrier} if $S$ is $\Bbb Z^{d+1}$-connected,  
 separates $\frak C_0$ from $\infty$ on $\Bbb Z^{d+1}$ and 
 satisfies the following condition:
 $$
 \align
 &\text{$S$ contains at least $|S|/6$ vertices $(\bold j,k+1)$ (with $k
 \in \Bbb Z$}\\
 &\text{arbitray) which have a parent $(\bold i,k)$ such that the edge from}\\
 &\text{$(\bold i,k)$ to $(\bold j,k+1)$ is a $\cD$-edge for which
 $C(\bold i,k,\bold j,\la)$ occurs}.
 \teq(3.14x)
 \endalign
 $$
 For $\frak C_0$ as above we define the quantity
 $$
 \Up(p,\la,\frak C_0)= \sum_{n \ge 1} P\{\text{there 
 exists a $\frak C_0$-barrier $S$ with $|S| = n$}\},
 \teq(3.15)
 $$
 where the probability is
 calculated in the $\{Y_t(\la)\}$-process which starts with a
 $B$-particle added to each $m(\bold c_i,0)$ in $\frak C_0$. (Recall
 that $p$ is the parameter which appears in the size of the blocks $\wh \cB$.)

 \proclaim{Lemma 7} Assume that there is a $\cD$-edge from $(\bold i,k)$
 to $(\bold j,k+1)$ and that $(\bold i, k)$ is active. Then 
 $$
 \{\text{the edge from $(\bold i,k)$ to $(\bold j, k+1)$ is closed}\}
 \subset 
 C(\bold i,k,\bold j,\la).
 \teq(3.13ab)
 $$

 If there exists some choice for 
 $p_0 \in \{1,2,\dots\}, \la_0 \in (0,\infty)$ and for the $\Bbb
 Z^d$-connected set $\{\bold c_1,
 \dots, \bold c_r\} \subset \Bbb Z^d$ such that 
 with $\frak C_0 = \{(\bold c_1,0), \dots, (\bold c_r,0)\}$ it holds
 $$
 \Up(p_0,\la_0, \frak C_0) < 1, 
 \teq(3.16)
 $$
 then 
 $$
 \align
 P\{&\text{infection survives in the process $\{Y_t(\la_0)\}$ which
 starts}\\
 &\qqquad \text{with a $B$-particle added at $m(\bold 0)$}\} > 0,
 \teq(3.17)
 \endalign
 $$
 and consequently $\la_c \ge \la_0$.
 \endproclaim

 \demo{Proof} Let there be a $\cD$-edge from $(\bold i,k)$
 to $(\bold j,k+1)$ and let $(\bold i, k)$ be active.
 By definition of
 an open edge, the edge from $(\bold i,k)$
 to $(\bold j,k+1)$ can be closed only if 
 $A(\bold i,k,\bold j) \cap B(\bold i,k,\bold j,\la) =
 \wt A\big(x(\bold i,k),t(k), \bold j\big) \cap \wt B\big(x(\bold i,k),t(k),
 \bold j, \la \big)$ fails. In addition $x(\bold i,k) \in m(\bold i) +
 \cC(\frac 18 \De_p)$ if \equ(3.7) holds. Thus the inclusion
 \equ(3.13ab) holds.

To prove \equ(3.17), assume that \equ(3.16) holds and 
 fix $p_0, \la_0$ and $\bold c_1, \dots, \bold c_r$
 such that $\frak C_0 = \{(\bold c_1,0),\dots, (\bold c_r,0)\}$ is a
 finite $\cD$-connected set for which \equ(3.16) holds.
\comment
 $$
 \Up(p_0,\la_0, \frak C_0) < 1.
 \teq(3.17abc)
 $$
 This can always be done, because if \equ(3.16) holds, then we can
 replace $\frak C_0$ by a large box $[-N,N]^d
 \times\{-1, 0, \dots, N]$ which contains $\frak C_0$, and is so large
 that \equ(3.17abc) holds (since the sets $S$ which appear in the 
 sum $\Up(p_0,\la_0,\frak C_0)$ must then lie outside $[-N,N]^d
 \times\{-1, 0, \dots, N]$).  
\endcomment
 Now consider the process $\{Y_t(\la_0)\}$ which has
 recuperation rate $\la_0$, but start it by adding at time 0 a $B$-particle to
 each site $m(\bold c_j),\; 1\le j \le r$. 
 Then all $(\bold c_j,0)$ are active. Let $\frak C$ be given by \equ(3.13)
 and view it as a subset of $\Bbb Z^{d+1}$.
 This cluster is the open cluster of $\frak C_0$. It contains $\frak
 C_0$ and is $\cL$-connected (see the lines following \equ(3.2) for $\cL$).
 Moreover, by Lemma 6, each vertex in $\frak C$ must be
 $\la_0$-active. 
Write $S$ for $\partial_{ext} \frak C$
 If $\frak C$ is finite, then $S$ is $\Bbb
 Z^{d+1}$-connected and separates $\frak C_0$ from infinity, by Lemma 5.
 Moreover, as observed after \equ(3.13),
 each $\cD$-edge  
 from some vertex $(\bold i,k) \in \frak C$ to some vertex $(\bold j,
 k+1) \in \partial_{ext}^+ \frak C$ must be closed.
 By virtue of \equ(3.4), $S$ then must contain at least $|S|/6$ vertices
 $(\bold j,k+1)$ which have a parent $(\bold i,k) \in \frak
 C$ with a closed $\cD$-edge from $(\bold i,k)$ to $(\bold j,k+1)$,
  and
 therefore such that $C(\bold i,k,\bold j,\la_0)$ occurs (by
 \equ(3.13ab) and the fact that $(\bold i,k)$ is active).
 Thus $S$ must have
 property \equ(3.14x). This implies that 
 $$
 \align
 &P\{\frak C \text{ is finite in the process $\{Y_t(\la_0)\}$
 which starts with}\\
 &\phantom{MM}\text{ a $B$-particle added to $m(\bold c_j)$ for each $1
 \le j \le r$}\}\\
 & \le \Up(p_0,\la_0, \frak C_0) < 1,
 \endalign
 $$
 or equivalently,
 $$
 \align
 &P\{\text{the infection survives in $\{Y_t(\la_0)\}$ if one 
 adds at time 0 a }\\
 &\phantom{MMMMMMM}\text{$B$-particle at each $m(\bold c_j), 
1 \le j \le r$}\}\\
 &\ge 
 P\{\frak C \text{ is infinite in $\{Y_t(\la_0)\}$} \text{ if one 
 adds at time 0 a }\\
 &\phantom{MMMMMMM}\text{$B$-particle at each $m(\bold c_j), 
1 \le j \le r$}\} > 0.
 \teq(3.19)
 \endalign
 $$

 It remains to show that the probability of survival of the infection
 remains strictly positive if we add a
 $B$-particle at time 0 only at $m(\bold 0)$. 
 In fact, this statement is still ambiguous, because, so far, we
 have only mentioned the locations were we add a $B$-particle at
 time 0, but we haven't specified how many particles we turn into 
 $B$-particles at these locations at time 0. 
To discuss this we remind the reader that $\overline \si$ was defined
 before Lemma 2 as the state obtained from $\si$ by changing all
 $A$-particles which coincide with a $B$-particle in the state $\si$ 
to $B$-particles. Now if we start in a random state $\si$ obtained by
 choosing $N_A(x,0-)$ $A$-particles at $x$ for i.i.d. Poisson
 variables $N_A(x,0-)$, and adding finitely many $B$-particles to the
 system,  then $\overline \si$ lies a.s. in $\Si_0$ (by Proposition 4
 in \cite {KSb}). Lemma 4 therefore shows 
 that the more particles we turn into $B$-particles at time 0, the more likely
 survival is. Therefore, the strongest conclusion to prove is that
 $$
 \align
 P\{&\text{the infection survives  if one 
 adds at time 0 a $B$-particle}\\
 &\text{at $m(\bold 0)$ only, and turns no $A$-particles to type $B$}\} > 0.
 \teq(3.19aa)
 \endalign
 $$
 And the weakest statement to start from is 
 $$
 \align
 P\{&\text{the infection survives  if one 
 adds at time 0 a }\\
 &\text{$B$-particle at each $m(\bold c_j), 1 \le j \le r$, and turns}\\
 &\text{all particles at
 these sites to $B$-particles}\} > 0.
 \teq(3.19bb)
 \endalign
 $$
 We shall prove \equ(3.19aa) from \equ(3.19bb). 
 Our argument for this is inspired by the proof at the
 bottom of p. 79 in \cite {D}. For simplicity we take the $\bold c_j,
 1 \le j \le r$, distinct. Only minor modifications are needed if some
 pairs of the $\bold c_j$ can be equal.
 Assume that we add a $B$-particle $\ze_j$ at $m(\bold c_j)$ for $1\le
 j \le r$ and
 that we turn all $A$-particles at these sites to $B$-particles at time 0. 
 Consider a sample point in the process $\{Y_t(\la_0)\}$
 with these initial conditions in which the 
 infection survives. If $v=(x,t)$ is a space-time point occupied by a
 $B$-particle in this process, then it has a genealogical path starting
 at some $B$-particle at one of the $\big(m(\bold c_j),0\big)$. 
 Thus, there exist times
 $0 < s_1 < \dots < s_\l < t$ and particles $\rho_0, \rho_1, \dots,
 \rho_\l$ such that $\rho_0$ is a $B$-particle at some $(\bold c_j,0),
 \; \rho_\l$ is at $v$ at time $t$, $\rho_i$ has type $B$ during
 $[s_i,s_{i+1}]$ and has no recuperation during this interval. 
Consider instead the process starting with a particle $\ze_j$ added at each  
$(m(\bold c_j),0)$, but now with only the particle 
$\rho_0$ of type $B$ and all other particles of type $A$. As
 before, induction on $i$ shows that $\rho_i$ is of type $B$ during
 $[s_i, s_{i+1}]$. 
\comment
 In other words, $\rho_\l$ has type $B$ at time $t$ even if
at time 0 one gives type $B$ only to the added particle $\ze_j$ and
 type $A$ to all other added particles, and if one turns no 
$A$-particles at any of the
 $m(\bold c_k)$ to $B$-particles. 
 In particular, $\rho_\l$ is even of type $B$ at time $t$ in the
 process which starts with only $\rho_j$ added as
a $B$-particle at $\big(m(\bold
 c_j),0\big)$, but with $\rho_k$ added as an
 $A$-particle at the point $\big(m(\bold
 c_k),0\big)$ for $k \ne j$, and and in which no types are changed at
 time 0.
\endcomment
In particular, there will be a $B$-particle at $(x,t)$ in this process 
with modified starting types.
 It follows from this same argument that  for
 any choice of the $a_k \ge 0, b_k \ge 0, \; 1 \le k \le r$,
 and with $\de(k,j) = 1$ or 0 according as
 $k=j$ or $k\ne j$, 
$$
 \align
 &P\{\text{infection survives}\big|N_A\big(m(\bold c_k),0\big)= a_k, 
 N_B\big(m(\bold c_k),0\big)= b_k,\,1 \le k \le r\}\\ 
 &\le \sum_{j=1}^r b_j P\{\text{infection survives}\big| N_A\big(m(\bold
 c_k),0\big)= a_k + b_k-\de(k,j),\\
&\phantom{MMMMMMMMMMMM} N_B\big(m(\bold c_j),0\big)= \de(k,j), 1 \le
 k \le r\},
 \teq(3.20)
 \endalign$$
 If \equ(3.19bb) holds, then 
 the left hand side of \equ(3.20) must be strictly positive for some
 choice of $a_k \ge 0, b_k \ge 1$. Consequently for some $j$
 and choice of the $a_k,b_k$, 
 $$
 \align
 &P\{\text{infection survives if one adds a $B$-particle at $\big(m(\bold
 c_{j}),0\big)$ only,}\\
 &\phantom{MMMMMMMMMMM}\text{but changes no $A$-particles to $B$-particles}\}\\
 &\ge P\{\text{infection survives }\big| N_A\big(m(\bold
 c_k),0\big)= a_k
 +b_k-\de(k,j), N_B\big(m(\bold c_k,0\big)= \de(k,j)\}\\
 &\phantom{MMMMMMMMMMM}\times  P\{N_A\big(m(\bold c_k),0-\big)
 = a_k+b_k-\de(k,j), 1 \le k \le r\}\\
 &> 0.
 \endalign
 $$
 By translation invariance we may assume $c_j = \bold 0$, so that
 \equ(3.19aa) follows.
 \hfill $\blacksquare$
 \enddemo

 To conclude the proof of $\la_c > 0$ we shall now establish that
 \equ(3.16) holds for suitable $p_0, \la_0, \frak C_0$. 
 Its proof has much in common with that of Proposition 3 in \cite
 {KSb}. First some
 more notation and definitions. 
 For purposes of comparison it is useful to couple our system with the system 
 in which there are no
 $B$-particles and in which all original
 $A$-particles move forever without interaction. In this system, which
 we shall denote by $\cP^*$, an $A$-particle $\rho$ which starts at $z$ will
 have position $z + \pi_A(t,\rho)$ for all $t$. 
 We write $N^*(x,s)$
 for the number of particles at the space-time point $(x,s)$ in the
 system $\cP^*$. $N^*(x,0)$ is taken equal to $N_A(x,0-)$, 
the initial number of
 $A$-particles at $x$. No initial $B$-particles are introduced in
 $\cP^*$ and all particles are of type $A$ at all times in $\cP^*$.
 Note that $N^*(x,s)$ is independent of the recuperation rate $\la$,
 because it depends only on the paths, and not the types, of the particles.
 It is easy to see that
 $$
 \align
 N^*(x,s) = &N_A(x,s)+N_B(x,s)- [\text{number of $B$-particles
 introduced}\\
 &\text{(in $\{Y_t(\la)\}$) at time 0 which are at $x$ at time $s$}].
 \teq(3.21cde)
 \endalign
 $$
 In particular
 $$
 N^*(x,s) \le N_A(x,s)+N_B(x,s). 
 \teq(3.21def)
 $$
 Next we define
 $$
 \cQ_p(x) = \prod_{s=1}^d [x(s),x(s)+C_0^p),
 \teq(3.20ab)
 $$
 and
 $$
 U_p(x,v)= \sum_{y \in \cQ_p(x)} N^*(y,v) = \sum \Sb y:x(s) \le y <
 x(s)+C_0^p\\ 1 \le s \le d \endSb  N^*(y,v).
 \teq(3.20bc)
 $$
 We call the bottom $\cZ_p(\bold i,k)= Z_p(\bold i) \times \{t_p(k)\}$ 
{\it good} if 
 $$
 U_p\big(x,t_p(k)\big) \ge \ga_0 \mu_AC_0^{dp}
 \text{ for all $x$ for
 which $\cQ_p(x)\subset Z_p(\bold i)$},
 \teq(3.20cd)
 $$
 where  $\ga_0$ is the constant introduced in (4.10), (4.16) and (4.17)
 of \cite{KSb}. 
 We also need the following technical estimate of some random walk
 probabilities. 
 \proclaim{Lemma 8} There exists a $p_1=p_1(d,D) <
 \infty$ such that if $\cZ_p(\bold i,k)$ is good and $\De_p \le u \le
 p^q\De_p, x \in \prod_{s=1}^d [(i(s)-4d)\De_p,
 (i(s)+4d+1)\De_p)$, then
 $$
 \sum_{y \in Z_p(\bold i)} [N_A\big(y,t(k)\big)+ N_B\big(y,t(k)\big)] 
P\{y+S_u = x\} \ge \frac 34 \ga_0\mu_A \text{ for } p \ge p_1.
 \teq(3.25a)
 $$
 \endproclaim
 \demo{Proof} This is nearly a copy of the proof of Lemma 5 in \cite {KSb}. 
 We introduce the blocks
 $$
 \cM(\pmb \l) := \prod_{s=1}^d [\l(s)C_0^p, (\l(s)+1)C_0^p).
 $$
 In our previous notation $\cM(\pmb\l) = \cQ_p(z)$ with $z(s) =
 \l(s)C_0^p$.
 These blocks have edge length only $C_0^p$, and 
 the set $Z_p(\bold i)$
 is a disjoint union of $(8d+3)^dC_0^{5dp}$ of these smaller blocks.
 Let 
 $$
 \La =\La(\bold i,p) = \{\pmb\l \in \Bbb Z^d:
 \cM(\pmb\l)\subset Z_p(\bold i)\}.
 $$ 
 Also, for each $\pmb\l \in \La$, let $y_{\pmb\l} \in \cM(\pmb\l)$ 
 be such that
 $$
 P\big\{y_{\pmb\l} +S_u =x\big\}
 = \min_{y \in \cM(\pmb\l)} 
 P\big \{y +S_u \in x\big\}.
 $$ 
 Then the left hand side of \equ(3.25a) equals
 $$
 \align
 &\sum_{\pmb \l\in \La} \sum_{y \in \cM(\pmb \l)}
 [N_A\big(y,t(k)\big)+ N_B\big(y,t(k)\big)]P\{y+S_u =x\}\\
 &\ge \sum_{\pmb\l \in \La} \sum_{y \in \cM(\pmb\l)}
 N^*\big(y,t(k)\big)P\{y_{\pmb\l}+S_u =x\}.
 \teq(3.25c)
 \endalign
 $$
 Since $\cZ_p(\bold i, k)$ is assumed to be good, we have
 $$
 \sum_{y \in \cM(\pmb\l)}N^*\big(y,t(k)\big)
  = U_p\big(\pmb\l C_0^p, t(k)\big)\ge \ga_0 \mu_A C_0^{dp}
 = \sum_{y \in \cM(\pmb\l)}\ga_0\mu_A.
 $$
 We can therefore continue \equ(3.25c) to obtain that the left hand side
 of \equ(3.25a) is at least
 $$
 \align
 &\sum_{\pmb\l \in \La} \sum_{y \in \cM(\pmb\l)}
 \ga_0\mu_A P\big\{y_{\pmb\l}+S_u =x \}\\ 
 &\ge  \sum_{\pmb\l \in \La} \sum_{y \in \cM(\pmb\l)}
 \ga_0\mu_A P\big\{y+S_u =x\}\\
 &\phantom{MM}-\sum_{\pmb\l \in \La} \sum_{y \in \cM(\pmb\l)}
 \ga_0\mu_A \big|P\big\{y_{\pmb\l}+S_u =x\}
 -P\big\{y+S_u = x\}\big|.
 \teq(3.21)
 \endalign
 $$

 Now, since $x \in \prod_{s=1}^d [(i(s)-4d)\De_p,
 (i(s)+4d+1)\De_p)$,
 the first multiple sum in the right hand side 
 of \equ(3.21) is at least
 $$
 \align
 &\sum_{w \in [-\De_p, \De_p)^d}
 \ga_0\mu_A P\{S_u = w\}\\
 &= \ga_0\mu_A\big[1- P\{S_u \notin [-\De_p, \De_p)^d\}\big]\\
 &\ge \ga_0\mu_A\big[ 1 - K_2\exp[-K_3^{-1}p^{-q}\De_p] \big]
 \teq(3.25e)
 \endalign
 $$
 for some constants $K_2(d,D), K_3(d,D)$. In the last 
 inequality we used simple large deviation estimates
 for $S_u$ (see for instance (2.40) in \cite{KSa}) and the fact that $u \le 
 p^q\De_p$. 

 The second multiple sum in the right hand side of
 \equ(3.21) is at most
 $$
 \ga_0 \mu_A\sum_{v \in \Bbb Z^d}\;\sup_{w:\|w-v\| \le dC_0^p} 
\big|P\{S_u = v\} - P\{S_u = w\}\big|.
 \teq(3.25f)
 $$
 This sum has already been estimated in the proofs of
 Lemmas 6 and 12 of \cite {KSa} (see in particular (5.26) there). 
 This sum is at most $K_4\ga_0 \mu_AC_0^p[\log u]^d u^{-1/2}$ 
 for some constant $K_4(d,D)$.

 For some $p_1(d,D)$ and all $p \ge p_1$ we finally have 
 from \equ(3.25e) and \equ(3.25f) that the left hand side of \equ(3.25a)
 is at least
 $$
 \ga_0\mu_A\Big[ 1 - K_2\exp[-K_3^{-1}p^{-q}\De_p]
 -K_4C_0^p[\log \De_p]^d \De_p^{-1/2} \big] 
 \Big] \ge \frac 34  \ga_0\mu_A.
 \tag "$\blacksquare$"
 $$
 \enddemo

 We define the  $\si$-fields
 $$
 \align
 \cH(\bold i,k) =\;& \cH_p(\bold i,k) =\si\text{-field 
 generated by the positions and
 types of}\\
 &\text{all particles at time 0, by
 all paths $\pi(\cdot, \rho)$ during $[0, t_p(k)]$}\\
 &\text{and by the paths  for all times of
 all particles outside $Z_p(\bold i)$}\\
 &\text{at time $t_p(k)$, 
and by all recuperation times $r(i,\rho)$ during }[0,t].
 \endalign
 $$
 We note that all $N_A(x,0), N_B(x,0)$ and the types of all particles
 at time $t(k)$ are $\cH(\bold i,k)$-measurable. Also the event that a
 given $x$ is occupied at time $t(k)$ belongs to $\cH(\bold i,k)$.
 The next lemma contains the crucial estimate for establishing \equ(3.16).
It proves that in the reset processes $\cZ_p(\bold i,k)$ 
with sufficiently large $p$, the infection will with high probability 
be transmitted, in a certain sense, 
 along a given $\cD$-edge. Recall that recuperation is ignored in 
the $\cZ_p(\bold i,k)$-process, so this is very similar to showing
 that the infection spreads with a certain minimal speed if 
 recuperation is not possible, as done in \cite {KSb}. For the
 infection to reach a certain cube $\cC$ (of size $\De_p/8$) we 
 define (as in \cite {KSb}) a random path along which a
 ``distinguished'' $B$-particle has a ``drift towards $\cC$.''
We use a corresponding martingale to show that with high probability 
the distinguished  $B$-particle has to follow the drift and will reach $\cC$.

 \proclaim{Lemma 9} Assume that there is a $\cD$-edge from $(\bold i,k)$
 to $(\bold j,k+1)$. 
 There exists a constant $p_2$ $($independent of $\bold i,k,\bold j,\la)$
  such that on the event 
 $$
 \{\cZ_p(\bold i,k) \text{ is good}\},
 \teq(3.23aa)
 $$
it holds
 $$
 \sum_{x \in m(\bold i) + \cC(\frac 18 \De_p)}
 P\{\wt A\big(x,t(k),\bold j\big) \text{ fails}|\cH(\bold
 i,k)\}
 \le \De_p^{-1} \text{ for } p \ge p_2.
 \teq(3.25g)
 $$
 \endproclaim
 \demo{Proof}
 Note that $\{\cZ_p(\bold i,k)\text{ is
 good\}} \in \cH(\bold i,k)$, because this event
 is defined in terms of the initial conditions, and
 paths during $[0,t(k)]$ only. 
  We divide the proof into 6 steps. A number of constants $p_i$ and
 $K_i$ will
 appear in this proof. These all depend only on $d,D,C_0$ and $\ga_0\mu_A$. 
We shall not make any further mention of this.
 \medskip 
 Step (i). We begin with a lower bound for the number of particles 
 in certain intervals in the $\cZ_p(\bold i,k)$-process started
 at $\big(x,t(k)\big)$ for some $x \in m(\bold i) + \cC(\frac
 18\De_p)$. To this end we define
 $$
 \align
 \wt K_p(z,v) = \;&\wt K_p\big(z,v;x,t(k)\big)
 = \big[\text{total number of particles in 
$\prod_{s=1}^d [z(s),z(s) + K_5p)$}\\
 &\text{at time
 $v$ in the $\cZ_p(\bold i,k)$-process started
 at $\big(x,t(k)\big)$}\big],
 \endalign
 $$
 for a constant $K_5$ to be chosen soon. We are interested in space
 time points $(z,v)$ satisfying 
 $$
 \align
 z \in \prod_{s=1}^d [(i(s)-4d)\De_p, &(i(s)+4d+1)\De_p - K_5p) 
 \text{ and }\\ 
 &v \in [t_p(k)+\De_p, t_p(k+1)],\quad v \in \Bbb Z.
 \teq(3.22)
 \endalign
 $$
 In this step we shall prove that
 we can choose $K_5$ and $p_3$ 
 such that on the event \equ(3.23aa) and for $p \ge p_3$
 $$
 P\big\{\wt K_p(z,v) < \frac 12 \ga_0  \mu_A \lfloor K_5p\rfloor ^d
 \text{ for some $(z,v)$ satisfying \equ(3.22)}
 \big|\cH(\bold i,k)\big\}
 \le \De_p^{-d-1}.
 \teq(3.23)
 $$

 Note that we are only interested in numbers of particles in \equ(3.23),
 irrespective of their types. 
 \comment
 We shall prove that, under \equ(3.22), 
 the probability in \equ(3.23) is at most
 $$ 
 \teq(3.25)
 $$
 \endcomment
 To see \equ(3.23), fix some $(z,v)$ in the set \equ(3.22).
 Now note that 
 if the $N_A\big(y,t_p(k)\big) +N_B\big(y,t_p(k)\big)$
 for $y \in Z_p(\bold i)$ are given, then the conditional distribution of 
 $\wt K_p(z,v)$
 for any fixed $(z,v)$, given $\cH(\bold i,k)$, 
 equals the distribution of $\sum_{y,n}X(y,n)$, where
 the $X(y,n)$ are
 independent binomial variables with
 $$
 \align
 P\{X(y,n) = 1\} &= 1-P\{X(y,n)=0\}\\
 &= P\{y+ S_{v-t(k)} \in \prod_{s=1}^d [z(s), z(s) + K_5p)\}\\
 &= \sum_{w \in \prod_{s=1}^d [z(s), z(s) + K_5p)} P\{S_{v-t(k)}= w-y\},
 \endalign
 $$
 and $\sum_{y,n}$ runs over $y \in Z_p(\bold i)$ 
 and, for given $y$, over $1 \le n \le N_A\big(y,t(k)\big)+
 N_B\big(y,t(k)\big)$. By Lemma 8, on \equ(3.23aa),
 $$
 \sum_{y,n}P\{X(y,n) = 1 \} \ge \frac 34 \ga_0\mu_A \Big|\prod_{s=1}^d
 [z(s), z(s) + K_5p)\Big| \ge
 \frac 34 \ga_0\mu_A \lfloor K_5p\rfloor ^d,
 $$
 provided $p \ge p_1$.
 Standard large deviation arguments now give that 
 $$
 P\big\{\wt K_p(z,v) < \frac 12 \ga_0  \mu_A \lfloor K_5p\rfloor^d
 \big|\cH(\bold i,k)\big\}
 \teq(3.23bb)
 $$
 is at most
 $$
 \exp\big[-\frac 18 \th_0 \ga_0 \mu_A \lfloor K_5p \rfloor ^d\big],
 $$ 
 for any $\th_0 > 0$ which satisfies $(3/4)\big(1-\exp[-\th_0]\big)\ge
 (5/8)\th_0$
 (compare (4.28) and the lines preceding it in \cite {KSb}).
 Thus we can choose $K_5 = K_5(d,D,C_0, \ga_0\mu_A)$ and $p_3 \ge p_1$ 
 such that 
 \equ(3.23bb) is at most $(8d+1)^{-d}p^{-q}\De_p^{-2d-3}$ for $p \ge p_3$.
 Since \equ(3.22) allows no more than $(8d+1)^dp^q\De_p^{d+1}$
 possible choices for $(z,v)$, \equ(3.23) then follows.
 \comment
 $$
 \align
 &P\{\cZ_p(\bold i,k)\text{ is good, but }
 \wt K_p(z,v) < \frac 12 \ga_0  \mu_AC_0^{dp} \text{ for some $(z,v)$
 satisfying \equ(3.22)}\}\\
 &\le K_9K_2p^q[\De_p]^{d+1}\exp\big[-\frac 18 \th_0 \ga_0 \mu_A \lfloor
 K_4p \rfloor^d\big] \\
 &\le K_9K_2 p^q \exp\big[p(d+1)\log C_0 
 -\frac 18 \th_0 \ga_0 \mu_A \lfloor K_4p \rfloor^d\big].
 \teq(3.27)
 \endalign
 $$
 \endcomment
 \medskip
 Step (ii). In this step we largely imitate Lemma 9 of \cite
 {KSb}. We define a path $\la(\cdot, \bold j)$ and use it to construct
 a martingale which shows that $\la(\cdot, \bold j)$ has a drift towards 
 $m(\bold j)$. $\la(s,\bold j)$ will be chosen for $t(k) \le s \le t(k+1)$
 according to the rules (i)-(v) below. In general, these rules do not determine
 $\la$ uniquely. Throughout this proof
$x \in m(\bold i) +\cC(\frac 18 \De_p)$ is
 a fixed vertex, occupied at time $t(k)$, and we are working in the
 $\cZ_p(\bold i,k)$-process started at $\big(x,t(k)\big)$. In
 particular, we do not allow recuperation and only consider particles
 which are in $Z_p(\bold i)$ at time $t(k)$, and the types of the
 particles refer to the types they have in the
$\cZ_p(\bold i,k)$-process started at $\big(x,t(k)\big)$. 
Here are rules (i)-(v):
 $$
 \align
 &(i) \phantom{iii} \la\big(t(k), \bold j) = x;\\
 &(ii)\phantom{ii} \text{for all $s \in [t(k), t(k+1)]$ there 
 is a distinguished $B$-particle,
 $\wh \rho (s)$ say, at} \\
 &\phantom{MMM}\text{$\la(s, \bold j)$;\;$\wh \rho(s)$ is a particle in the 
 $\cZ_p(\bold i,k)$-process
 started at $\big(x,t(k)\big)$};\\
 &\phantom{MMM}\text{at time $t(k),\; \wh \rho\big(t(k)\big)$ is the unique 
 $B$-particle at $x$ in the $\cZ_p(\bold i,k)$-process}\\
 &\phantom{MMM}\text{started at $\big(x,t(k)\big)$};\\
 &(iii) \; s \mapsto \la(s,\bold j) \text{ can jump only at times when $\wh
 \rho(s-)$ jumps away from $\la(s, \bold j)$}\\
 &\phantom{MMM}\text{and $\la(\cdot, \bold j)$
 is constant between such jumps};\\
 &(iv)\phantom{ii} \text{if $\wh \rho(s-)$ jumps from 
 $\la(s-,\bold j)= w$ to $w'$ at some time $s$},\\
 &\phantom{MMM}\text{and if this was the only $B$-particle 
 at $w$ at time $s-$,}\\
 &\phantom{MMM}\text{then $\la(\cdot, \bold j)$ also jumps to $w'$ at 
 time $s$ (so that
 $\la(s, \bold j) = w'$),}\\
 &\phantom{MMM}\text{and $\wh \rho(s) = \wh \rho(s-)$, the
 particle which jumped at time $s$};\\
 &(v)\phantom{iii} \text{if $\wh \rho(s-)$ jumps from $\la(s-, \bold j)= w$
 to $w'$ at some time $s$ such that there}\\
 &\phantom{MMM}\text{is at least one other
 $B$-particle $\rho'$ at $w$ at time $s$, then $\la(\cdot, \bold j)$ 
jumps to}\\
 &\phantom{MMM}\text{$w'$ at time $s$ if and only if 
 $\|w'-m(\bold j)\|_2 < \|w-m(\bold
 j)\|_2$, and in this}\\
 &\phantom{MMM}\text{case again $\wh \rho(s) = \wh
 \rho(s-)$; if, however,
 $\|w'-m(\bold j)\|_2 \ge \|w-m(\bold j)\|_2$,}\\
 &\phantom{MMM}\text{then $\la(\cdot, \bold
 j)$ does not jump at time $s$ and we take $\wh \rho(s) = \rho'$}.
 \endalign
 $$
 Note that rules (iv) and (v) depend on whether there is another $B$-particle 
 than $\wh \rho$ at $\la(s-,\bold j)$. In \cite {KSb} any particle at
 the same space-time point as $\wh \rho$ automatically had type $B$,
 but this is not the case in the present setup, because a jump is
 required before an $A$-particle can turn into a $B$-particle. This
 fact will necessitate a few extra remarks in the next step.

 As in \cite {KSb} (4.42), (4.43) we now define
 $$
 \align
 I_1(u)& = I[N_B\big(\la(u,\bold j),u\big) = 1]\\
& = I[\wh \rho(u)\text{ is the
 only $B$-particle present at }(\la(u,\bold j),u)];\\
 I_{\ge 2}(u) &= I[N_B\big(\la(u,\bold j),u\big) \ge 2];
 \endalign
 $$
 and with $e_{d+i} = -e_i$ for $1 \le i \le d$,
 $$
 \align
 \Ga_1(u)& = \frac 1{2d} \sum_{i=1}^{2d} \big[\|\la(u,\bold
 j)+e_i-m(\bold j)\|_2
  - \|\la(u,\bold j)-m(\bold j)\|_2\big];\\
 \Ga_{\ge 2}(u) &= \frac 1{2d} \sum {}^* 
 \big[\|\la(u,\bold j)+e_i-m(\bold j)\|_2
  - \|\la(u,\bold j)-m(\bold j)\|_2\big],
 \teq(3.50abc)
 \endalign
 $$
 where $\sum {}^*$ is the sum over those $i \in \{1, \dots,2d\}$ for which
 $$
 \|\la(u,\bold j)+e_i-m(\bold j)\|_2 - \|\la(u,\bold j)-m(\bold j)\|_2 < 0.
 $$
 Finally we define, for $t \ge t(k)$,
 $$
 M(t) = M(t,\bold j):=\|\la(t,\bold j)-m(\bold j)\|_2
 -D\int_{t(k)}^t \big[I_1(u)\Ga_1(u) +
 I_{\ge 2}(u) \Ga_{\ge 2}(u)\big] du
 \teq(3.64abc)
 $$
(recall that $D$ is the jump rate af all particles).
 The result of this step is that $M(t)$
 is a right continuous martingale under the measure 
 which governs the $\cZ_p(\bold i,k)$-process started at
 $\big(x,t(k)\big)$, conditioned on $\cH(\bold i,k)$, or conditioned on
 the 
 $$
 N_A\big(y,t(k)\big)+ N_B\big(y,t(k)\big) \text{ for }y \in Z_p(\bold i).
 \teq(3.30aa)
 $$
 The proof is essentially the same as that of Lemma 9 in \cite {KSb},
 so we leave this to the reader. We merely remark that the
 $\cZ_p(\bold i,k)$-process started at $\big(x,t(k)\big)$
 depends only on the variables in \equ(3.30aa),
 the paths on $[t(k), \infty)$ of the particles in $Z_p(\bold i)$
 at time $t(k)$, and lastly, on the choice of which particle at
 $\big(x,t(k)\big)$ is given type $B$. However, changing this choice from
one particle to another amounts to interchanging the roles of two
 particles at $\big(x,t(k)\big)$.
Since all particles move and recuperate in the same way,
 such an interchange does not influence the distribution of $M(t)$ for
 later $t$.
Thus, conditioning on $\cH(\bold i,k)$ or on the
 variables in \equ(3.30aa) gives the same distribution for
 $\{M(s)\}_{t(k) \le s \le t(k+1)}$.
 \medskip
 Step (iii). In this step we start on  lower bounds for
 $$
 Z(t)= Z(t,k) := \int_{t(k)}^{t(k)+t} I_{\ge 2}(u)du
 \teq(3.30abc)
 $$
 and for
 $$
 \int_{t(k)}^{t(k)+t} \big[I_1(u)\Ga_1(u) +
 I_{\ge 2}(u) \Ga_{\ge 2}(u)\big] du.
 $$
 These lower bounds are essentially the same as in \cite
 {KSb}.
 Following \cite {KSb} we define 
 for an integer $L \ge 2$
 $$
 \be(L,d)= \cases 1 &\text{ if } d =1\\ [\log L]^{-1} 
  &\text{ if } d = 2\\ L^{2-d}  &\text{ if } d \ge 3,
 \endcases
 \teq(3.80cba)
 $$
 $$
 \align
 \cE_n = &\cE_n(\bold j,k) =\{\text{there is some particle } 
 \rho' \ne \wh \rho\big(t(k)+3L^2(n-1)\big) \\
&\text{of the $\cZ_p(\bold i,k)$-process started at $(x,t(k))$ in}\\
 &\la\big(t(k)+3L^2(n-1),\bold j\big)+ [-L,L]^d \text{ at time }
 t(k)+3L^2(n-1)\}, 
 \teq(3.80abc)
 \endalign
 $$
 and 
 $$
 \align
 J_n = I\Big[&\text{at some time }u \in \big(t(k) + 3L^2(n-1), 
t(k)+ L^2(3n-1)\big],
  \wh \rho(u) \text{ coincides}\\
 &\text{with another $B$-particle
 in the $\cZ_p(\bold i,k)$-process started at $\big(x,t(k)\big)$} \Big].
 \endalign
 $$
 The only differences of any consequence 
 with the definitions just before Lemma 11 of \cite {KSb} is the
 insistence in the definition of $J_n$ that $\wh \rho$ coincide 
 with another  $B$-particle, and that this be a particle in the
  $\cZ_p(\bold i,k)$-process started at $\big(x,t(k)\big)$.
 A $B$-particle is necessary because it is also possible for an $A$ and
 $B$-particle to be at the same space-time point. 
 The particle is required to belong to the
  $\cZ_p(\bold i,k)$-process started at $\big(x,t(k)\big)$, because these are
  the only particles under consideration at the moment.
 Despite these differences we have just
 as in Lemma 11 of \cite {KSb} that for all $L \ge 1$
 $$
 E\{J_n|\cF_{3L^2(n-1)}\} \ge K_7\be^2(L,d) \text{ on the event } \cE_n,
 \teq(3.80a)
 $$
 where this time, for $s \ge t(k)$,
 $$
 \align
 \cF_s = \si &\text{-field generated by  $\cH(\bold i,k)$ and
 the paths of}\\
 &\phantom{MMMMMMM} \text{all particles in $\cZ_p(\bold i,k)$ on $[t(k),s]$},
 \endalign
 $$
 and the expectation is with respect to the $\cZ_p(\bold i, k)$-process
 started at $\big(x,t(k)\big)$; $K_7$ is some constant which depends on
 $d,D$ only. To prove \equ(3.80a) we first observe that there is no loss
 of generality to assume that $\wh \rho$ and $\rho'$ are at different 
 locations at time $3L^2(n-1)+1$ on the event $\cE_n$. This is so 
 because there is a probability of at least $e^{-D}(1-e^{-D})$
 that $\wh \rho$ does not jump,
 and  $\rho'$ has one jump in $[3L^2(n-1),3L^2(n-1)+1)$.
 If $\wh \rho$ and $\rho'$ are at distinct sites,
 then $\rho'$ will have to have type $B$ no later than the first time
 when these particles get together. With this obervation the proof 
 of \equ(3.80a) is the same as the proof of Lemma 11 in 
 \cite {KSb}.
\comment
 that $J_n=1$ on the intersection of the events (4.68) and
 (4.69). (The argument after that in Lemma 12 of \cite{KSb} is only
 used for the case $D_A \ne D_B$, which is not considered
 here. Therefore, one may even replace $\ga^2$ by $\ga$ in \equ(3.80a),
 but this is of no consequence for us.)
\endcomment

 This time we take 
 $$
 L= \lceil K_5p\rceil \text{ and } \overline t = p^q\De_p = t(k+1) - t(k).
 \teq(3.80cde)
 $$
 These values for $L$ and $t$ will remain fixed for the rest of this proof.
 Without loss of generality we take $p_3$ so
 large that  $L \ge 1$ for $p \ge p_3$. Still following \cite {KSb}
 we define
 $$
 V( \overline t,L) =V( \overline t,L, \bold j,k):= 
 \sum_{1+ \De_p/(3L^2) \le n \le [t(k+1)-t(k)]/(3L^2)} I[\cE_n(\bold j,k)].
 $$
 We further define the event
 $$
 D(x,t(k), \bold j) :=\big\{\|\la(u, \bold j) -m(\bold j)\|_2 \ge \frac
 1{16}\De_p \text{ for all } u \in [t(k),t(k+1)]\big\}.
 \teq(3.78)
 $$
 This event was not used in \cite {KSb}. Nevertheless,
 with $Z$ given by \equ(3.30abc), we can essentially copy the 
proof of Lemma 12 in \cite {KSb}. We use that
$$
E\{J_n|\cF_{3L^2(n-1)}\} \ge K_7\be(L,d) I[\cE_n(\bold j,k)],
$$
and consequently
$$
E\big\{\min\{1, \int_{3L^2(n-1)}^{3L^2n} I_{\ge 2}(u)du\}
\big|\cF_{3L^2(n-1)} \big\} \ge e^{-2D}K_7\be(L,d) I[\cE_n(\bold j,k)],
$$
as in \cite {KSb}. This yields
\comment
we have, with 
 $$
 \ep = \ep(d,D) = \min\Big[\frac {K_7}{15 e^{2D}},1\Big ]
 $$
 in Lemma 13 in \cite {KSb},
 \endcomment
for $p \ge p_3,\; 0 < \ep \le 1$ and $K_3$ a universal constant
 $$
 \align
 &P\{Z( \overline t) \le \ep\be(L,d)L^{-2} \overline t 
\text{ and $D(x,t(k), \bold j)$
 occurs}|\cH(\bold i,k)\} \\
 &\le 
 P\big\{V( \overline t,L) \le \frac {2\ep}{K_7}e^{2D}L^{-2} \overline t
 \text{ and $D(x,t(k),\bold j)$ occurs}\big|\cH(\bold i,k)\}\\
 &\phantom{MMMMMMMMMMMMMM} +
 2\exp\big[-\frac{K_3}3 \ep^2\be^2(L,d)L^{-2} \overline t\big].
 \teq(3.79)
 \endalign
 $$

 Next we recall for the reader the bounds (4.67)
 and (4.68) of \cite {KSb}. These bounds say that for $\la, x \in \Bbb Z^d$
 $$
 \frac 1{2d} \sum_{i=1}^{2d} \big[\|\la+e_i-x\|_2
  - \|\la-x\|_2\big] \le \frac{K_{12}}{\|\la -x\|_2+1},
 \teq(3.93a)
 $$
 and, with $\sum {}^*$ as in \equ(3.50abc),
 $$
 \frac 1{2d} \sum {}^* \big[\|\la+e_i-x\|_2
  - \|\la-x\|_2\big] \le
  -K_{13}+\frac{K_{12}}{\|\la -x\|_2+1}
 \teq(3.93b)
 $$
 for some constants $K_{12}, K_{13}$ which depend on $d$ only.
 Moreover the left hand sides of \equ(3.93a) and \equ(3.93b) are at most
 1 in absolute value.

 It is immediate from \equ(3.93a), \equ(3.93b) and \equ(3.50abc)
 that on $D(x,t(k),\bold
 j)$ it holds
 $$
 \align
 &\int_{t(k)}^{t(k)+ \overline t} \big[I_1(u)\Ga_1(u) +
 I_{\ge 2}(u) \Ga_{\ge 2}(u)\big] du\\
 &\le [t(k+1)-t(k)]\frac {16 K_{12}}{\De_p}
 -K_{13}Z\big(t(k+1)-t(k)\big)\\ &=
 16K_{12}p^q -K_{13}Z\big(t(k+1)-t(k)\big).
 \teq(3.93cd)
 \endalign$$

 \medskip
 Step (iv). Here we use the martingale $M(\cdot)$ to estimate
 $P\{D(x,t(k),\bold j)|\cH(\bold i,k)\}$ 
 in terms of the distribution of $V(t,L)$.
 To this end we note first that
 $$
 \align
 M\big(t(k)\big) &= \|\la\big(t(k),\bold j\big) - m(\bold j)\|_2 = \|x-m(\bold
 i) +m(\bold i) - m(\bold j)\|_2 \\
 &\le \sqrt d \frac 18 \De_p+ \|\bold i
 - \bold j\|_2 \De_p \le \sqrt d \frac 98 \De_p,
 \teq(3.48ab)
 \endalign
 $$ 
where we used the fact that there is a $\cD$-edge from $(\bold i,k)$
 to $(\bold j,k+1)$.
 On the other hand, on the event $D(x,t(k),\bold j)$ we have from
 \equ(3.93cd) that
 $$
 \align
 M\big(t(k+1)\big) &\ge -D \int_{t(k)}^{t(k)+ \overline t} 
\big[I_1(u)\Ga_1(u) +
 I_{\ge 2}(u) \Ga_{\ge 2}(u)\big] du\\
 &\ge -16DK_{12}p^q +K_{13}DZ\big(t(k+1)-t(k)\big).
 \teq(3.48bc)
 \endalign
 $$
 Further we have the martingale inequality (4.55) of \cite {KSb}: for
 some constants $K_{14}-K_{16}$ which depend on $D$ only, $a \ge 2+2D$,  
 $0 \le b \le 1$ and $T \ge t(k)$
 $$
 \align
 P\Big\{\sup_{t(k)\le s \le T} \big|M(s) &-M\big(t(k)\big)\big| 
 \ge a+b\big(T-t(k)\big)\Big|\cH(\bold i,k)\Big\} \\
 &\le K_{14}\exp\big[-K_{15}\big(T-t(k)\big)\big] + 2\exp[-K_{16}ab].
 \teq(3.47)
 \endalign
 $$
 For $a= \De_p/2, b=\De_p/[2\big(t(k+1) - t(k)\big)]$ and $T =
 t(k+1)$ this implies
 $$
 P\big\{\big|M\big(t(k+1)\big) -M\big(t(k)\big)\big| \ge
 \De_p\big|\cH(\bold i,k)\}
 \le 4\exp\big[ -\frac{K_{16}\De_p}{4p^q}\big],
 \teq(3.47aa)
 $$
 provided $p \ge p_4$ for some constant $p_4\ge p_3$.
 Combined with \equ(3.48ab) and \equ(3.48bc) this gives
 $$
 \align
 &P\{D(x,t(k),\bold j)|\cH(\bold i,k)\} \\
 &\le P\{M\big(t(k+1)\big)-M\big(t(k)\big) \ge \De_p|\cH(\bold i,k)\}\\
 & \phantom{MMM}+ P\big\{-16DK_{12}p^q +K_{13}DZ\big(t(k+1)-t(k)\big) \le
 M\big(t(k)\big)+\De_p \\
 &\phantom{MMMMMMMMMMMMMMMMMMM}\text{ and $D(x,t(k),\bold j)$ occurs}
 |\cH(\bold i,k)\big\}\\
 &\le  4\exp\big[ -\frac{K_{16}\De_p}{4p^q}\big]\\
 & \phantom{MMM}+P\big\{Z\big(t(k+1)-t(k)\big) \le \frac 1{K_{13}D}
 \big[\big(\sqrt d \frac
 98 +1\big)\De_p +16DK_{12}p^q\big]\\
 &\phantom{MMMMMMMMMMMMMMMMMMM}\text{ and $D(x,t(k),\bold j)$ occurs}
 \big|\cH(\bold i,k)\big\}\\
 &\le  4\exp\big[ -\frac{K_{16}\De_p}{4p^q}\big]
 +P\big\{Z\big(t(k+1)-t(k)\big) \le \frac {3\sqrt
 d}{K_{13}D}\De_p \\
 &\phantom{MMMMMMMMMMMMMMMMMMM}\text{ and
 $D(x,t(k),\bold j)$ 
 occurs}\big|\cH(\bold i,k)\big\},
 \teq(3.47ab)
 \endalign
 $$
 for $p \ge $ some $p_5$.
 Finally we note that $\be(L,d) \ge L^{1-d}$
 (see \equ(3.80cba)), so that for 
$$
\ep = \ep(d,D) = \min\Big[\frac{K_7}{15e^{2D}}, 1\Big]
$$ 
and $p \ge$ a suitable
 constant $p_6$,
 $$
 \ep \frac{\be(L,d)}{L^2} p^q\De_p \ge \frac{3\sqrt d}{K_{13}D}
 \De_p \text{ (recall that $q = 2d+1$)}.
 $$
 Thus, by \equ(3.79), we can continue \equ(3.47ab) to
 obtain
 $$
 \align
 &P\{D(x,t(k),\bold j)|\cH(\bold i,k)\} 
 \le 4\exp\big[ -\frac{K_{16}\De_p}{4p^q}\big]\\
 &\phantom{MMM}+P\big\{V( \overline t,L) \le \frac 2{15}L^{-2}
 \overline t \text{ and $D(x,t(k),\bold j)$ occurs}
 \big|\cH(\bold i,k)\}\\
 &\phantom{MMMMMMMMMMMMM} + 2\exp\big[-\frac{K_3}3
 \ep^2\be^2(L,d)L^{-2}
 \overline t\big].
 \teq(3.47bc)
 \endalign
 $$

 \medskip
 Step (v). In this step we shall estimate the probability in the 
 right hand side of \equ(3.47bc). This will be done by 
 using the following direct consequence of the definitions of $\wt K$
 and of $\cE_n$:
 $$
 \wt K_p\big(\la(t(k) + 3L^2(n-1),\bold j),t(k) + 3L^2(n-1)\big) 
  \ge 2
 $$
 implies that $\cE_n$ occurs (recall that $L = \lceil K_5p \rceil$). 
 This will allow us to use the bound
 \equ(3.23) on the probability that $\wt K$ is `small'. We turn to 
 the details. Take $p_6$ so large that for all $p \ge p_6$ 
 $$
 \frac 12 \ga_0\mu_A (K_5p)^d \ge 2. 
 $$ 
Assume now that the events
 $$
\big \{\wt K_p(z,v) \ge \frac 12 \ga_0  \mu_A \lfloor K_5p\rfloor^d
 \text{ for all $(z,v)$ satisfying \equ(3.22)}\big\}
 \teq(3.31)
 $$ 
 and 
 $$
\big \{\la(u, \bold j) \in \prod_{s=1}^d [(i(s)-4d)\De_p, 
(i(s)+4d+1)\De_p - K_5p)  \text{ for } t(k) \le u < t(k+1)\big\}
 \teq(3.32)
 $$
 occur. 
 Then for
 $$
 1+ \De_p/(3L^2) \le n \le [t(k+1)-t(k)]/(3L^2),
 \teq(3.32bc)
 $$
 it holds
 $$
 \wt K_p\big(\la(t(k) + 3L^2(n-1),\bold j),t(k) + 3L^2(n-1)\big)  
 \ge \frac 12 \ga_0  \mu_A \lfloor K_5p \rfloor^d \ge 2.
 \teq(3.33)
 $$
 As observed, this implies that $\cE_n$ also occurs for the $n$ in
 \equ(3.32bc) and then also 
 $$
 V( \overline t,L) \ge \frac{t(k+1)-t(k)-\De_p}{3L^2} -2 >
 \frac 2{15L^2}p^q\De_p
 \teq(3.32cd)
  $$ 
 for $p \ge $ some constant $p_7$.

 However, we know from \equ(3.23) that on the event \equ(3.23aa),
 \equ(3.31) indeed holds outside a set of conditional 
 probability $\De_p^{-d-1}$.
 To estimate the probability that \equ(3.32) fails for a relevant value
 of $u$, we introduce the random time
 \comment
 $$
 \si := \inf\{t \in [t(k),t(k+1)]:
 \la(t, \bold j) \notin 
 \prod_{s=1}^d [(i(s)-4d)\De_p, (i(s)+4d+1)\De_p - K_5p)\},
 \teq(3.35)
 $$
 $$
 \tau = \inf\{t \in [t(k),t(k+1)]:\la(t, \bold j) \in m(\bold j)
 +\cC(\frac 18 \De_p)\}.  
 \teq(3.36)
 $$
 \endcomment 
 $$
 \tau := \inf\{w \in [t(k),t(k+1)]:
 \la(w, \bold j) \notin 
 \prod_{s=1}^d [(i(s)-4d)\De_p, (i(s)+4d+1)\De_p - K_5p)\}.
 \teq(3.35)
 $$
 This definition for $\tau$ holds 
 if the set in the right hand
 side of \equ(3.35)
 is not empty; otherwise we set $\tau$ equal to
 $t(k+1)$. 
 We shall prove in the remainder of this step that
 $$
 \align
 &P\{\la(u, \bold j) \notin 
 \prod_{s=1}^d [(i(s)-4d)\De_p, (i(s)+4d+1)\De_p - K_5p)\\
 &\phantom{MMMM}
 \text{ for some } u \in \big[t(k),t(k+1)\big)
 \text{ and $D(x,t(k),\bold j)$ occurs}
 |\cH(\bold i,k)\}\\
 &= P\{\tau < t(k+1) \text{ and $D(x,t(k),\bold j)$ occurs}
 \big|\cH(\bold i,k)\}\\
 &\le 4\exp\big[- \frac{K_{16}\De_p}{4p^q}\big]
 \teq(3.34)
 \endalign
 $$
 for $p \ge$ some constant $p_8$. As we observed
 above, this will imply
 $$
 \align
 &P\{V( \overline t,L) \le \frac 2{15L^2}p^q\De_p
 \text{ and $D(x,t(k),\bold j)$ occurs}|\cH(\bold i,k)\}\\
 &\le P\{\text{\equ(3.31) fails}|\cH(\bold i,k)\}
 +P\{\equ(3.32)\text{ fails but $D(x,t(k),\bold j)$ occurs}|\cH(\bold i,k)\}\\
 &\le \De_p^{-d-1}+ 4\exp\big[- \frac{K_{16}\De_p}{4p^q}\big],
 \teq(3.34bb)
 \endalign
 $$
 for $p \ge p_8$ and on the event \equ(3.23aa).

 To prove \equ(3.34) we note that on the event $D(x, t(k),\bold j)$
 it holds
 $$
 \int_{t(k)}^\tau \big[I_1(u)\Ga_1(u) +
 I_{\ge 2}(u) \Ga_{\ge 2}(u)\big] du
 \le 16K_{12}p^q
 $$ 
 (see \equ(3.93cd)). Consequently, on $\{\tau < t(k+1)\} \cap
 D(x,t(k),\bold j)$ we also have
 $$
 \align
 M(\tau) &\ge \|\la(\tau,\bold j)-m(\bold j)\|_2 - 16DK_{12}p^q\\
 &\ge \|\la(\tau,\bold j)-m(\bold i)\|_2 - \|m(\bold i) - m(\bold j)\|_2 -
  16DK_{12}p^q\\
 &\ge 4d \De_p - K_5p - \sqrt d \De_p -16DK_{12}p^q 
 \ge \De_p +\sqrt d \frac 98 \De_p,
 \endalign
 $$
 provided $p \ge$ some constant $p_9$.
 In particular, on $\{\tau < t(k+1)\} \cap
 D(x,t(k),\bold j)$ 
 $$
 \sup_{t(k) \le s \le t(k+1)}|M(s)-M\big(t(k)\big)| \ge M(\tau) -
 M\big(t(k)\big) \ge \De_p \text{ (see \equ(3.48ab))}.
 $$
 We already proved in \equ(3.47) and \equ(3.47aa) that 
 $$
 P\{\sup_{t(k) \le s \le t(k+1)}|M(s)-M\big(t(k)\big)|\ge
 \De_p|\cH(\bold i,k)\} \le
 4\exp\big[-\frac{K_{16}\De_p}{4p^q}\big].
 $$
 This implies the promised inequality \equ(3.34).

 \medskip
 Step (vi). We finally prove the lemma in this step, by combining
 the preceding steps.
 It follows from the definition of $\wt A(x,t(k),\bold j)$ that
 $[\wt A(x,t(k), \bold j)]^c$ can occur only if 
 $\big(x,t(k)\big)$ is occupied,   
 but 
 $$
 \align
 &\text{the $\cZ_p(\bold i,k)$-process
 started at $\big(x,t(k)\big)$ does not have}\\
 &\phantom{MMMMMMM}\text{a $B$-particle in $m(\bold
 j) + \cC(\frac 18\De_p)$ at time }t(k+1). 
 \teq(3.50)
 \endalign
 $$
 In turn, this last event can occur only if
 $D(x,t(k),\bold j)$ occurs, or if 
 $$
 \|\la(u,\bold j)-m(\bold j)\|_2 \le \frac 1{16}\De_p \text{ for some } u
 \in [t(k),t(k+1)],   
 \teq(3.51)
 $$
 as well as \equ(3.50) occur. However, the probability of the
 intersection of \equ(3.51)
 and \equ(3.50) is small. Indeed, if 
 $$
 \si :=\inf\{u \ge t(k):\|\la(u,\bold j)-m(\bold j)\|_2 \le \frac 1{16}\De_p\},
 $$
 then on the intersection of \equ(3.51) and \equ(3.50), 
 the particle $\wh \rho (\si)$ is within distance
 $\De_p/(16)$ of $m(\bold j)$ at time $\si$, but outside $m(\bold j) +
 \cC(\frac 18 \De_p)$ at time $t(k+1)$. (Recall that we are working
 with the $\cZ_p(\bold i,k)$-process started at
 $\big(x,t(k)\big)$. This has no recuperation, so the particle $\wh
 \rho(\si)$ will still be a $B$-particle in this process at time $t(k+1)$.)
 In other words, the particle $\wh \rho(\si)$, 
which is the distinguished one at time $\si$, travels a distance at
 least $\De_p/8 - \De_p/16 = \De_p/16$ during $[\si, t(k+1)]$. Consequently,
 $$
 \align
 &P\{\text{\equ(3.51) and \equ(3.50) occur}|\cH(\bold i,k)\}\\
 &\le P\{\sup_{s \le t(k+1)-t(k)} \|S_{t(k+1)}-S_s\|_2 \ge \frac 1{16}\De_p
 |\cH(\bold i,k)\}\\
 &\le 8d \exp\big[ - K_{17}\frac{\De_p}{p^q}\big] \text{ (see (2.42)
 in \cite {KSa})}
 \endalign
 $$
 for some constant $K_{17} = K_{17}(d,D)$. It follows that
on the event \equ(3.23aa)
 $$
 \align
 &P\{\equ(3.50) \text{ occurs}|\cH(\bold i,k)\}\\
 & \le 
 P\{D(x,t(k),\bold j)|\cH(\bold i,k)\} +
  8d \exp\big[ - K_{17}\frac{\De_p}{p^q}\big]\\
 &\le 4\exp\big[ -\frac{K_{16}\De_p}{4p^q}\big]\\
 &\phantom{MM}+P\big\{V(\overline t,L) \le \frac 2{15L^2}\overline t
 \text{ and $D(x,t(k),\bold j)$ occurs}
 \big|\cH(\bold i,k)\}\\
 &\phantom{MM} + 2\exp\big[-\frac{K_3}3 \be^2(L,d)L^{-2}\overline t\big]
 + 8d \exp\big[ - K_{17}\frac{\De_p}{p^q}\big]
 \text{ (by \equ(3.47bc))}\\
 &\le 4\exp\big[ -\frac{K_{16}\De_p}{4p^q}\big]
  + 2\exp\big[-\frac{K_3}3 \be^2(L,d)L^{-2}\overline t\big]\\
 &\phantom{MM}+ 8d \exp\big[ - K_{17}\frac{\De_p}{p^q}\big]
 +\De_p^{-d-1}+ 4\exp\big[- \frac{K_{16}\De_p}{4p^q}\big]\\
  \endalign
 $$
 (see \equ(3.34bb)) for $p \ge$ a suitable constant $p_2$. 
 Summation of this estimate over the at most $[1+\De_p/4]^d$
 possible $x \in m(\bold i) +\cC(\frac 18\De_p)$ now 
 proves \equ(3.25g).
 \hfill $\blacksquare$
 \enddemo

 For the remainder of this section we shall only consider $p \ge p_2$.
 \proclaim{Lemma 10} For $p \ge p_2$ we can choose $\la_0 = \la_0(p) > 0$
 such that on the event \equ(3.23aa),
 $$
 \sum_{x \in m(\bold i) +\cC(\frac 18\De_p)}
 P\{[\wt A(x,t(k),\bold j) \cap \wt B(x,t(k),\bold j,\la_0)]^c
\text{ \rom {in} } \{Y_t(\la_0)\} |\cH(\bold i,k)\} \le 2\De_p^{-1},
 \teq(3.60)
 $$ 
 and 
 $$
 P\big\{C(\bold i,k,\bold j,\la_0)\text{ \rom{in} } \{Y_t(\la_0)\}
|\cH(\bold i,k)\big\} \le 2\De_p^{-1}.
 \teq(3.61)
 $$
 \endproclaim
 \demo{Proof} The left hand side of \equ(3.60) equals
 $$
 \align
 \sum_{x \in m(\bold i) +\cC(\frac 18\De_p)}
 &P\{[\wt A(x,t(k),\bold j)]^c|\cH(\bold i,k)\}\\
 &+\sum_{x \in m(\bold i) +\cC(\frac 18\De_p)}
 P\{\wt A(x,t(k),\bold j)\cap [\wt B(x,t(k),\bold j,\la_0)]^c|\cH(\bold i,k)\}.
 \endalign
 $$
 In view of Lemma 9 it therefore suffices for \equ(3.60) to prove that
 on the event $\{\cZ_p(\bold i,k)$ is good\}
 $$
 \sum_{x \in m(\bold i) +\cC(\frac 18\De_p)}
 P\{\wt A(x,t(k),\bold j)\cap [\wt B(x,t(k),\bold j,\la_0)]^c
\text{ in } \{Y_t(\la_0)\}|\cH(\bold i,k)\}
 \le \De_p^{-1}.
 \teq(3.62)
 $$
 We claim that each summand in \equ(3.62) is at most 
 $$
 1-\exp\big[-\la_0[t(k+1)-t(k)]\big]\le \la_0[t(k+1)-t(k)].
 $$
 Indeed, by the definitions \equ(3.13x), \equ(3.13y) of $\wt A$ and $\wt B$, 
 once we know that $\wt A(x,t(k),\bold j)$ occurs, the event 
 $[\wt B(x,t(k),\bold j,\la_0)]^c$ can occur only if some particle
 $\rho_i$ has a recuperation event in $\{Y_t(\la_0)\}$ during
 $[s_i,s_{i+1}]$, for $0 \le i \le \l$. 
 Here $\rho_i$ are certain particles and the $s_i$ are increasing and
 such that $s_{\l+1} - s_0 = t(k+1)-t(k)$. These $\rho_i$ and $s_i$ 
 are determined by the $\cZ_p(\bold i,k)$-process started at
 $\big(x,t(k)\big)$, and therefore are independent of the recuperation
 events during $[t(k),\infty)$. 
 This proves our claim. \equ(3.62) now follows for some small
 $\la_0(p)$, since there are at most $[1+\De_p/4]^d$ terms in the sum in
 \equ(3.62).

 The preceding paragraph proves \equ(3.60). \equ(3.61) is now an
 immediate consequence of \equ(3.13z) and the fact that $C(\bold i, k,
 \bold j,\la) = \emptyset$ if $(\bold i,k)$ is not active.
 \hfill $\blacksquare$
 \enddemo
 Lemma 10 will help us to bound the probability that there are many sites 
 $(\bold i,k)$ in an open cluster $\frak C$ with a good bottom and with
 a closed edge from $(\bold i,k)$ to a site in $\partial_{ext} C$. In
 order to obtain \equ(3.16) from such a bound we first have to show
 that there is only a small probability that a $\frak C_0$-barrier $S$  
 has of order $|S|$ sites with a parent (on $\cD$) with a bad bottom. 
 (A bottom $\cZ_p(\bold i,k)$ is called {\it bad} if it is not good.) 
 This will be the goal of
 Lemmas 11 and 12. For a $\Bbb Z^{d+1}$-connected set $S$ and
 integer $p \ge 2$ we define
 $$
 S_p^* = \bigcup \Sb (\bold {i'},k): \|\bold {i'}- \bold j\| \le 4d-1\\
 \text{ for some }(\bold j,k+1) \in S
 \endSb  \wh \cB_p(\bold i',k).
 \teq(3.63)
 $$
 Recall that we used the vertex $(\bold i,k)$ as a kind of renormalized
 site to replace the block $\wh B_p(\bold i,k)$. Forming $S_p^*$ from
 $S$ is a construction going in the other direction. From the
 collection $S$ of renormalized sites we reconstruct the blocks which
 are near the ones represented by sites in $S$. We define
 further for any positive integer $\nu$ and $r \ge p$ the blocks
 $$
 \cL_{r,\nu}(\bold m,u) := \prod_{s=1}^d [\nu m(s)\De_r,
 \nu (m(s)+1)\De_r) \times [\nu u\De_r, \nu(u+1)\De_r).
 $$
 Note that the blocks $\cL_{r,\nu}(\bold m,u)$ form a partition  of
 $\Bbb Z^{d+1}$ into disjoint cubes.
 \proclaim{Lemma 11} Let $r \ge p \ge 2$ and $q = 2d+1$ (as before).
 There exists a constant $K_{18}$, depending on $d$ only, such that for 
 each $\Bbb Z^{d+1}$-connected set $S$ and
 each integer $\nu \ge 1$,
there exists a $\Bbb Z^{d+1}$-connected
 set $\La_{p,r}(S,\nu)  \subset \Bbb Z^{d+1}$ such that
 $$
 |\La_{p,r}(S,\nu)| \le K_{18}\Big[\frac {|S|\De_p}{\nu \De_r} +1\Big]p^q
 \teq(3.64)
 $$
 and such that 
 $$
  \bigcup_{(\bold m,u) \in \La_{p,r}(S,\nu)} 
 \cL_{r,\nu}(\bold m,u) \supset S^*_p.
 \teq(3.65)
 $$
 \endproclaim
 \demo{Proof} This lemma is essentially the same as lemma 1 in \cite {CGGK}.
 For the convenience of the reader we repeat the main steps of the
 proof. Let $|S|=n$. Since $S$ is connected it has a spanning tree 
 with $n-1$ edges, and then there exists a path $(v_0,v_1, \dots, v_a)$
 on $\Bbb Z^{d+1}$ of length $a \le 2n-2$ whose vertices are exactly
 the vertices of $S$ (some vertices are repeated; the path is not self-avoiding,
 in general). 
 For $0 \le u \le a$ let $v_u = (\bold i_u,k_u)$.
 Now set 
 $$
 \mu = \nu \frac{\De_r}{\De_p} = \nu \De_{r-p},
 $$
 and consider the vertices $v_{j\mu}$ for $0 \le j \le a/\mu$ (note
 that $\mu$ is an integer, by our choice of $C_0$ and $\De_r$).
 For $0 \le j \le a/\mu$
 let $(\bold m_j,u_j)$ be the unique vertex in 
 $\Bbb Z^{d+1}$ such that 
 $$
 (\bold i_{j\mu}\De_p, k_{j\mu}p^q\De_p) \in \cL_{r,\nu}(\bold m_j,u_j).
 \teq(3.66)
 $$
 We now take 
 $$
 \align
 \La_{p,r}(S,\nu) =  \bigcup \big\{&(\bold m, u): 
 \text{ there exists a $0 \le j \le a/\mu$ 
 such that}\\
 &\text{$\|\bold m - \bold m_j\| \le (4d+2)$ and $|u - u_j| 
 \le 3p^q$}\big\}.
 \endalign
 $$
 Then, since $(\bold m_j,u_j)$ takes
 at most $(a/\mu +1) \le (2n\De_p/\nu\De_r +1)$ values, it holds
 $$
 |\La_{p,r}(S, \nu)| \le \Big[\frac {2n\De_p}{\nu\De_r} + 1\Big]
(8d+5)^d (6p^q+1),
 $$
 and \equ(3.64) holds.

 Next we verify \equ(3.65).
 Assume that $y$ is a vertex in $S^*_p$. Then there is some $v_u
 = (\bold i_u,k_u)$ and some $(\bold i',k')$ with $|i'-i_u| 
 \le 4d-1,\; k' = k_u-1$
 such that $y \in \wh \cB_p(\bold i',k')$. In particular, $|y(s)-
 i_u(s)\De_p| \le 4d\De_p, 1 \le s \le d$, and $|y(d+1)- k_u p^q\De_p|
 \le p^q\De_p$. 
 Also, there exists some $j$ such that
 $j\mu\le u < (j+1)\mu$, and $0 \le j \le a/\mu$, and 
consequently $\|(\bold i_u,k_u) - (\bold
 i_{j\mu},k_{j\mu})\| \le \mu$. Finally, by virtue of \equ(3.66),
 $$
 \align
 \|(\bold i_{j\mu}\De_p,k_{j\mu}p^q\De_p)&- 
 (\mu \bold m_j\De_p, \mu u_j\De_p)\|\\
 &=\|(\bold i_{j\mu}\De_p,k_{j\mu}p^q\De_p)- 
 (\nu\bold m_j\De_r, \nu u_j\De_r)\| \le \De_r
 \endalign$$
 and
 $$
 \align
 &|y(s) - \nu m_j(s)\De_r| \\
 &\qquad \le |y(s) - i_u(s)\De_p| + |i_u(s)\De_p -
 i_{j\mu}(s) \De_p| + |i_{j\mu}(s)\De_p - \nu m_j(s)\De_r|\\
 &\qquad \le 4d\De_p+(\nu+1)\De_r \le (4d+2)\nu\De_r, \; 1 \le s \le d,
 \endalign
 $$ 
 and similarly
 $$
 |y(d+1) - \nu u_j\De_r| \le 3p^q \nu \De_r.
 $$
 The relation \equ(3.65) now follows easily.

 Finally, the connectedness of $\La_{p,r}(S,\nu)$ follows from the fact
 that $\La_{p,r}(S,\nu)$ is the union of the rectangular boxes 
 $$
 \prod_{s=1}^d [m_j(s)-4d-2,m_j(s)+4d+2] 
\times [u_j-3p^q,u_j+3p^q], \; 1\le j
 \le a/\mu.
 $$
 These boxes are clearly $\Bbb Z^{d+1}$-connected and the boxes
 corresponding to the two successive values $j$ and $j+1$
 intersect, because they have the point $(\bold m_{j+1}, u_{j+1})$ in
 common. Clearly, this point lies in the $(j+1)$-th box.
It also lies in the $j$-th box, because 
$m_j(s) = \lfloor i_{j\mu}(s)/\mu \rfloor$ and
$m_{j+1}(s) = \lfloor i_{(j+1)\mu}(s)/\mu \rfloor,\, 1 \le s \le d$,
 (by \equ(3.66)) and
$|i_{(j+1)\mu}(s) - i_{j\mu}(s)| \le
 \|v_{(j+1)\mu}-v_{j\mu}\| \le \mu$.
Similarly, $u_j= \lfloor k_{j\mu}p^q/\mu \rfloor$
and $|k_{(j+1)\mu} - k_{j\mu}| \le \mu$.
 \hfill $\blacksquare$
 \enddemo

In the next lemma $\frak C_0$ always will be such that
$$
\align
&\frak C_0 \subset \cD,\;\bold 0 \in \frak C_0, \text{ and}\\
&\quad\text{$\frak C_0$ is
$\Bbb Z^{d+1}$-connected when viewed as a 
subset of $\Bbb Z^{d+1}$}.
\teq(3.68zz)
\endalign
$$
\proclaim{Lemma 12} We can choose $\la_0 > 0, p_0$ and $\frak C_0$ such that
\equ(3.68zz) and \equ(3.16) are satisfied.
 \endproclaim
 \demo{Proof} The hard part of the work was done in \cite {KSa} and
 \cite {KSb}. It is too long to repeat and we shall be content with
 reducing the lemma to some results in those references.
 \medskip
 Step (i). Basically we are going to show 
that there is only a small  probability that there
 exists a $\frak C_0$-barrier with
 `many' vertices which are $\cD$-adjacent to a vertex with a bad
 bottom (see \equ(3.20cd) for the definition of a good bottom).
In this step we reduce the bounding of the number of
 barriers with a large number of vertices that have a parent 
 $(\bold i,k)$ in $\cD$ with a bad bottom to estimates in \cite {KSb}.

 Let $\frak C_0$ have the properties in \equ(3.68zz).
It is easy to see that then any  
 $\frak C_0$-barrier $S$ must contain some vertices $(k_+, 0, \dots,0)$
 and $(k_-,0,\dots,0)$ on the positive and negative first coordinate
 axes, respectively. Moreover, if $|S| = n$, then
 $$
 \text{diameter}(S) = \max_{x,y \in S}  \|x-y\| \le n.
 $$
 It follows that $1 \le k_{\pm} \le n$ and that $S \subset
 [-n,n]^{d+1}$. Since $S$ must be $\Bbb Z^{d+1}$ connected and must
 contain $(k_+,0,\dots,0)$ for some $1 \le k_+ \le n$, there are
 at most $n[K_{19}]^n$ possibilities for $S$, (with $K_{19}$ depending
 on $d$ only; see for instance \cite
 {Ka}, formula (5.22)). 
 \comment
 Once $S$ is given there are at most $2^n$ possibilities
 for the subset of $(\bold j,k+1)$ here, and then since there has to be a 
 $\cD$-edge from $(\bold i,k)$ to $(\bold j,k+1)$, there are in total at most 
 $n[2\cdot3^dK_{19}]^n$ choices for the collection of $(\bold i,k)$. 
 Some of the corresponding blocks $\cB_p(\bold i,k)$ may have a good
 bottom and others a bad bottom.

 Again there are at most $n$ vertices
 $(\bold i,k)$ which appear here, and therefore at most $2^n$ choices
 for the subset of $(\bold i,k)$ for which $\cB_p(\bold i,k)$ has a bad bottom.
 Now fix
 a set $T$ of at most $n$ vertices which is one of the sets of
 $(\bold i,k)$ which can appear as the set of (indices of) blocks with
 a bad bottom in the preceding argument. From the above, there are at
 most $n[4\cdot3^dK_{19}]^n$ for $T$. 
 We want to show that uniformly in $T$ with 
 $$
 |T| \ge \frac 1{12\cdot 3^d}n
 \teq(3.67)
 $$
 $$
 P\{\cB_p(\bold i,k) \text{ has a bad bottom for all } (\bold i,k) \in T\} 
 \le 
 \teq(3.68)
 $$
 \endcomment
 We remind the reader that in \cite {KSb} we also defined bad
 $r$-blocks of the form 
 $$
 \cB_r(\bold m, \l) := \prod_{s=1}^d [m(s)\De_r, (m(s)+1)\De_r)
 \times [\l\De_r, (\l+1)\De_r),
 \teq(3.68)
 $$
and their {\it pedestals}
$$
\cV_r(\bold m, \l) := \prod_{s=1}^d[(m(s)-3)\De_r, (m(s)+4)\De_r)
\times \{(\l-1)\De_r\}.
\teq(3.68aa)
$$
Note that $\cB_r(\bold m,\l r^q) \subset \wh \cB_r(\bold m, \l)$.
The block in \equ(3.68) is called bad (in the sense of \cite {KSb})
 if (see \equ(3.20ab) and \equ(3.20bc) for $\cQ_p$ and $U_p$)
 $$ 
 \align
 &U_r(x,v) < \ga_r\mu_AC_0^{dr} \text{ for some $(x,v)$ with integer
 $v$ for which}\\
 &\cQ_r(x) \times \{v\}\subset \prod_{s=1}^d
 [(m(s)-3)\De_r,(m(s)+4)\De_r) \times [(\l-1)\De_r, (\l+1)\De_r).
 \endalign
 $$
Similarly, the pedestal in \equ(3.68aa) is called bad (in the sense of 
\cite {KSb}) if
$$
\align
&U_r(x,(\l-1)\De_r) < \ga_r\mu_AC_0^{dr} \text{ for some $x$ for which }\\
&\cQ_r(x) \subset \prod_{s=1}^d
 [(m(s)-3)\De_r,(m(s)+4)\De_r).
\endalign
$$
 Here $U_r$ is given by \equ(3.20bc)
 and the $\ga_r$ are increasing in $r$ and satisfy 
 $$
 0 < \ga_0 \le \ga_r \le \ga_\infty \le \frac 12,\; r \ge 0
 $$ 
 for some $\ga_0, \ga_\infty$. The precise form of the $\ga_r$ used in
 \cite {KSb} is not important at the moment.
 If $\wh \cB_p(\bold i,k)$ has a bad bottom as defined in
 \equ(3.20cd), then
 $$
 U_p\big(x,t_p(k)\big) < \ga_0 \mu_AC_0^{dp}
 \text{ for some $x$ for
 which $\cQ_p(x)\subset Z_p(\bold i)$}.
 $$
 In this case, $\cQ_p(x) \subset \prod_{s=1}^d [(i'(s)-3)\De_p,
 (i'(s)+4)\De_p)$ for some $i'$ with 
 $$
 i(s)-4d+2 \le i'(s) \le i(s)+4d-2, \; 1 \le s \le d.
 \teq(3.70)
 $$
 Therefore, $\cB_p(\bold {i'},kp^q)$ is bad in the 
sense of \cite {KSb} for
 some $\bold{i'}$ satisfying \equ(3.70).

 Now suppose $\wh \cB_p(\bold i,k)$ has a bad bottom 
 and a $\cD$-edge to $(\bold j,k+1) \in S$. Then there is an $\bold
 {i'}$ with $\|\bold {i'} - \bold i\| \le 4d-2$, and hence
 $\|\bold{i'} - \bold j\| \le 4d-1$, such that $\cB_p(\bold
 {i'},kp^q)$ is bad in the sense of \cite {KSb}. By the definition \equ(3.63)
 $\wh \cB_p(\bold {i'},k) \subset S_p^*$, and therefore, $\cB_p(\bold
 {i'},kp^q) \subset \wh \cB_p(\bold {i'},k) \subset S_p^*$.  Moreover, 
$\|\bold {i'} - \bold j\| \le
 4d-1$, so that at most $(8d-1)^d$ vertices $(\bold j,k+1) \in S$ can
 give rise to the same $(\bold {i'},k)$. This shows that 
 $$
 \align
 &\text{[the number of $(\bold j,k+1) \in
 S$ with a parent $(\bold i,k)$ for which}\\
 &\phantom{MMMMMMMMMMMMMMM}\text{$\wh \cB_p(\bold i,k)$ has a bad bottom]}\\
 &\le (8d-1)^d \times
 \text{ [the number of bad $p$-blocks $\cB_p(\bold {i'},kp^q)$ in 
the sense of}\\
&\phantom{MMMMMMMMMMMMMMM}\text{\cite {KSb} contained in $S_p^*]$}.
 \teq(3.71)
 \endalign
 $$
In the
 next step we apply \cite {KSb} to estimate the right hand side of \equ(3.71).

 \medskip
 Step (ii). In analogy with \cite {KSb} we now make the following
 definitions for a barrier $S$. In these definitions, an $r$-block 
is of the form \equ(3.68) and `good' or `bad' are meant in the sense
 of \cite {KSb}.
 $$
 \wh \phi_r(S^*_p) = \text{number of bad $r$-blocks which intersect $S^*_p$},
 \teq(3.72)
 $$
$$ 
\align
\wh \psi_r(S^*_p) = \;&\text{number of $r$-blocks which intersect $S^*$
 and }\\
&\text{which have a good pedestal, but contain a bad $(r-1)$-block},
\endalign
$$
 $$
 \wh \Phi_r(n) = \wh \Phi_r(n, \frak C_0) = \sup\{\wh \phi_r(S^*_p): S 
 \text { a $\frak C_0$-barrier 
 of cardinality $n$}\},
 \teq(3.73)
 $$
and
$$
\wh \Psi_r(n) = \wh \Psi_r(n, \frak C_0) = \sup\{\wh \psi_r(S^*_p): S 
 \text { a $\frak C_0$-barrier 
 of cardinality $n$}\}.
 \teq(3.73aa)
 $$
 In this step we shall prove that for
 every choice of $K$ and $\ep_0$, there exist constants $p_0, n_0$
 such that for all $p \ge p_0, n \ge n_0$,
 $$
 P\{\wh \Phi_r(n,\frak C_0) \ge  \ep_0n \text{ for some } 
 r \ge p\} \le \frac 2 {n^K}. 
 \teq(3.74bb)
 $$
 It is crucial that this estimate is uniform in $\frak C_0$ 
satisfying \equ(3.68zz). In fact,
 $p_0,n_0$ depend only on $d,\ga_0, \mu_A, \ep_0,K$, but not on $\frak C_0$.

 We saw in step (i) that all $\frak C_0$-barriers 
 $S$ of cardinality $n$ have to satisfy
 $$
 S \subset [-n,n]^{d+1}, 
 \teq(3.74)
 $$
 so that we may restrict the sup in \equ(3.73) to $S$ which satisfy the
 condition \equ(3.74). The quantities $\wh \phi_r$ and $\wh \Phi_r$ are
 analogues of the following quantities introduced in \cite {KSb}
 $$
 \align
 \phi_r(\wh \pi) := \;&\text{number of bad (in the sense of \cite {KSb}) 
 $r$-blocks which}\\
 &\text{intersect the space-time path $\wh \pi$}
 \endalign
 $$
 and
 $$ 
 \Phi_r(\l) = \sup_{\wh \pi \in \Xi(\l,t)} \phi_r(\wh \pi),
 \teq(3.75)
 $$
 with
 $$
 \align
 \Xi(\l,t) = \{&\wh \pi: \wh \pi\text{ is a space-time path over the
 time interval $[0,t]$ and}\\
 &\text{located in $\cC(t\log t)$, with exactly $\l$ jumps 
 during $[0,t]$}\}.
 \teq(3.75b)
 \endalign
 $$
 We showed in \cite {KSb}, Proposition 8,  that for any choice of $K$
 and $\ep_0 > 0$, there exist constants $r_0,t_1$, such that for all $t
 \ge t_1$
 $$
 P\{\Phi_r(\l) \ge \ep_0C_0^{-6r} (t+\l)\text{ for some } r \ge r_0, \l
 \ge 0\} \le \frac 2{t^K}.
 \teq(3.74aa)
 $$
 One can check that the lengthy proof of \equ(3.74aa) uses the
 restriction that $\wh \pi \in \Xi(\l,t)$ in the sup in \equ(3.75) only
 for the bound in (4.32) in \cite {KSb}. This bound says (after a small
 change to the present notation) that for integers
 $\nu \ge 1$ and $r \ge p$, the
 number of blocks $\cL_{r,\nu}(\bold m,u)$
 which intersect any given $\wh \pi \in \Xi(\l,t)$ is at most 
 $$
 \la(\l) := 3^d\Big (\frac{t+\l}{\nu \De_r}+2 \Big).
 \teq(3.76)
 $$
 In the present case we can replace this estimate by 
 \equ(3.64).
 This tells us that for $r \ge p$, 
 any set $S^*_p$ defined by \equ(3.63) for $S$ a 
$\frak C_0$-barrier of cardinality $n$,
 intersects at most 
 $$
 K_{18} \big[\frac{n\De_p}{\nu \De_r} +1\big]p^q
 \teq(3.77)
 $$ 
 blocks $\cL_{r,\nu}(\bold m,u)$. Apart from an insignificant change
 from the factor $3^d$ to $K_{18}$ this takes the place of the bound
 \equ(3.76), provided we replace  $(t+\l)$ by $n\De_p p^q$. 
We further have to replace $R(t)$ of (4.16) in \cite{KSb} by $\wh
 R(n)$, which we take to be the unique integer $R$ for which
$$
C_0^R \ge [K_4\log n]^{1/d} > C_0^{R-1}.
$$
If diameter $(S) = n$, then by \equ(3.74) and \equ(3.63)
$$
S_p^* \subset [(n-4d+1)\De_p, (n+4d-1)\De_p)^d \times [-p^q\De_p,
np^q\De_p).
$$
Simple estimates for the Poisson distribution (compare Lemmas 5 and 9
in \cite {KSa}) show then that we can
take $K_4=K_4(d, \mu_A,K)$ so large that
$$
P\{\wh \Phi_r(n) > 0 \text{ for any }r \ge \wh R(n)\lor \log p\} \le \frac
1{n^K}, \quad n\ge 1.
$$ 
This estimate takes the place of (4.17) in \cite {KSb}. We can then
follow the proof of Lemma 7 in \cite {KSb} with only trivial changes
to show that there exist constants $C_5, \ka_0, n_0$ which depend on
$d,\ga_0, \mu_A,K$ (but not on $p,r,n$ or $\frak C_0$), such that for
all $n \ge n_0, \ka \ge \ka_0, p \le r \le \wh R(n)-1$,
$$
\align
P\big\{\wh \Psi_{r+1} (n) \ge \frac{\ka
n}{\De_{r+1}}&p^q\De_p[\rho_{r+1}]^{1/(d+1)}\\
& \le 
\exp\Big[-nC_5\ka p^q\De_p\exp\big[-\frac{\ga_0\mu_A} {2(d+1)}
C_0^{(d-3/4)r}\big]\Big],  
\endalign
$$
where 
$$
\rho_{r+1}= 3^{d+1}C_0^{6(d+1)(r+1)}\exp[-\frac
12\ga_r\mu_AC_0^{(d-3/4)r}].
$$
With this estimate in hand one can copy the proof of Proposition 8
in \cite {KSb} with the simple replacement of $\ka_0(t+\l)/\De_{r+1}$
by $\ka_0 n p^q\De_p/\De_{r+1}$.
This yields that
$$
\wh \Phi(n) \le n6\ka_0
C_0^{6(d+1)}\De_p p^q\exp\Big[-\frac{\ga_0\mu_A}{2(d+1)}C_0^{(d-3/4)r}\Big]
< \ep_0n
$$
outside a set of probability $2n^{-K}$, for $n \ge n_0$ and $
r_0(d, \ga_0,\mu_A, \ep_0)\lor p \le r \le \wh R(n)-1$. 
By taking $p_0 \ge r_0$ and $r
\ge p \ge p_0$ one obtains \equ(3.74bb).
 \medskip
 Step (iii). Without loss of generality we take $p_0 \ge p_2$
 (which was determined in Lemma 10).
For $p \ge p_0, K = 2$ and $r = p$, \equ(3.71) and \equ(3.74bb),
 imply that for any $n_1 \ge n_0$ and any $\frak C_0$ satisfying \equ(3.68zz)
 $$
 \align
 \sum_{n \ge n_1} P\{&\text{there exists a $\frak C_0$-barrier $S$ with $|S|=n$
 and
 at least}\\ 
 &\text{$(8d-1)^d \ep_0n$ vertices which have
 a parent $(\bold i,k)$}\\
 &\text{such that $\wh \cB_p(\bold i,k)$ has a bad bottom}\}
 \le \sum_{n \ge n_1} \frac 2{n^2}.
 \endalign
 $$
 We shall take $\ep_0 = \ep_0(d)$ such that  
 $(8d-1)^d\ep_0 = 1/(12)$. We further take $n_1$ so large that
  $\sum_{n \ge
 n_1}2n^{-2} \le 1/3$. Finally, we fix
$$
\frak C_0 = \{(k,0,\dots,0):0 \le k \le n_1\}.
$$
It is clear that there do not exist any sets $S$ of fewer than $n_1$
elements which separate 
this segment of the first coordinate axis from $\infty$ in $\Bbb
Z^{d+1}$. 
Thus
$$
\align
 \sum_{n \ge 1} P\{&\text{there exists a $\frak C_0$-barrier $S$ with $|S|=n$
 and}\\
  &\text{at least $n/12$ vertices which have
 a parent }\\
 &\text{$(\bold i,k)$ such that $\wh \cB_p(\bold i,k)$ has a bad bottom}\}
 \le \frac 13.
\endalign
$$

Now, for $S$ to be a $\frak C_0$-barrier, it
 must contain a subset of at least $n/6$ vertices $(\bold j,k+1)$ 
 which have a parent $(\bold i,k) \in \frak C_0$ such that $C(\bold
 i,k,\bold j, \la_0)$ occurs (see \equ(3.14x)).
In view of our last estimate,
it suffices for \equ(3.16) to prove that
 $$
 \align
 &\sum_{n \ge n_1} P\{\text{there exists a $\Bbb Z^{d+1}$-connected 
set $S$ with $|S| = n$, which  
 separates $\frak C_0$ }\\
 &\qquad\text{from $\infty$ on $\Bbb Z^{d+1}$ and which
 contains at least $n/12$ vertices $(\bold j,k+1)$ with}\\
 &\qquad\text{a parent
 $(\bold i,k)$ such that $\wh \cB_p(\bold i,k)$ has a good bottom and
 $C(\bold i,k, \bold j,\la)$ occurs}\}\\
 & < \frac 23.
 \teq(3.80)
 \endalign
 $$
 In this step will shall prove \equ(3.80).

 Now suppose we are given any set of vertices $(\bold i_1, k_1),\dots, (\bold
 i_m,k_m)$ and further $\bold j_r, \; 1\le r \le m$, such that $(\bold
 i_r,k_r)$ is a parent of $(\bold j_r,k_r+1)$.
 Assume that 
 $$
 \|\bold i_r -\bold i_s\| \ge 8d+7  \text{ for all $r,s$ with $k_r = k_s$}.
 \teq(3.82)
 $$
 We claim that then for $p \ge p_0$ and $\la =\la_0(p)$ 
 $$
\align
& P\{\wh \cB_p(\bold i_r, k_r) \text{ has a good bottom, but 
 $C(\bold i_r,k_r, \bold j_r,\la_0)$ occurs}\\
&\phantom{MMMMMM}\text{in $\{Y_t(\la_0)\}$ 
  for all } 1 \le r \le m\} \le [2\De_p^{-1}]^m. 
 \teq(3.81)
 \endalign
$$
 Recall that $\la$ is the recuperation
 rate and $p$ is the parameter determining the block sizes used to
 define $\wt A$ and $\wt B$; $\la_0$ was
 determined in Lemma 10. \equ(3.81) is immediate from
 \equ(3.61). To see this, assume without loss of generality that $k_r \le
 k_s$ for all $1 \le r \le s\le m$. Then for $r <s$,
 $$
 C(\bold i_r,k_r, \bold j_r, \la) \in \cH(\bold i_s,k_s).
 \teq(3.83)
 $$
 Indeed, for $k_r< k_s$ this follows from the fact that 
 $C(\bold i_r,k_r,\bold j,\la)$
 depends on information during $[0, t_p(k_r+1)] \subset [0,t_p(k_s)]$
 only. For $k_r = k_s$ but $r < s$ \equ(3.83) follows from the fact that 
 $C(\bold i_r,k_r, \bold j_r, \la)$ depends only on information during
 $[0,t_p(k_r)]$ and on particles in  $Z_p(\bold i_r)$ at time
 $t_p(k_r) = t_p(k_s)$, and $Z_p(\bold i_r) \cap Z_p(\bold i_s) = \emptyset$ by
 virtue of \equ(3.82). Thus \equ(3.83) holds. We already remarked that
 also $\{\cZ_p(\bold i_r,k_r) \text{ is good}\} \in \cH(\bold i_r,k_r)$.
 Therefore, 
 $$
 \align
 P\big\{\cZ_p(\bold i_s,k_s) \text{ is good and }&\text{ 
$ C(\bold i_s, k_s,\bold j_s,\la_0)$ in }\{Y_t(\la_0)\}\big|
\cZ_p(\bold i_r,k_r) \text{ is}\\
 &\text{good and } C(\bold i_r, k_r,\bold j_r,\la_0)
 \text{ for }r < s\} \le 2[\De_p]^{-1},
 \endalign
 $$
 by \equ(3.61). \equ(3.81) follows.

 The rest of the proof is routine. We already
 saw in Step (i) that there are at most $n[K_{19}]^n$ possible $\frak
 C_0$-barriers $S$ of size $n$. For such a set $S$ to have the property
 in \equ(3.80), it must have a subset of at least $n/12$ vertices
 $(\bold j,k+1)$ with a parent $(\bold i,k)$ such that $C(\bold
 i,k,\bold j, \la)$ occurs. There are at most $2^n$ choices for the set of
 $(\bold j, k+1)$, since $S$ has only $2^n$ subsets. When these  
 $(\bold j, k+1)$ have been chosen, they have at most $3^dn$  
 parents $(\bold i,k)$, because any vertex has
 at most $3^d$ parents. A subset of these parents must have good bottoms. 
 Thus the total number of choices for the set of 
 $(\bold i,k)$ for which $\wh \cB_p(\bold i,k)$ must have a good bottom,
 while also $C(\bold i, k,\bold j, \la_0)$ occurs, is at most
 $$
 n[K_{19}2\cdot2^{3^d}]^n.
 $$
The number of parents $(\bold i,k)$
 needed so that each of $n/12$ vertices has at least one parent among  
 these $(\bold i,k)$ is at least $3^{-d}(n/12)$. 
 And for at least one choice of $\bold a = \big(a(1), \dots,
 a(d)\big)$, 
the residue class 
 $\{i_r(s)\equiv a(s) \pmod {8d+7}, 1 \le s \le d\}$ has at least 
 $(8d+7)^{-d}3^{-d}n/12$ members.
 By \equ(3.81), the probability that for all these $(\bold i,k)$
 $\wh \cB_p(\bold i,k)$ has a good bottom, while $C(\bold i,k,\bold j,\la_0)$
 occurs, is at most
 $$
 \big[2\De_p^{-1}\big]^{(8d+7)^{-d}3^{-d}(n/12)}.
 $$
 Thus, the $n$-th summand in \equ(3.80) is at most
 $$
 n[K_{19}2\cdot2^{3^d}]^n \big[2\De_p^{-1}\big]^{(24d +21)^{-d}(n/12)} \le 
 n\Big[K_{20}[\De_p]^{-1/(12 \cdot (24d+21)^d)}\Big]^n,
 $$
 for some constant $K_{20}(d)$. This shows that \equ(3.80) holds for large
 enough $p$ and completes the proof of the Lemma.
 \hfill $\blacksquare$
 \enddemo

 As pointed out in Lemma 7, \equ(3.16) implies that $\la_c >0$.

\subhead
4. The maximal number of jumps in a path
\endsubhead
\numsec=4
\numfor=1
We need a few definitions to state the purpose of this section. 
In this section we consider only the system of $A$-particles and there
is no interaction between any particles. Accordingly, recuperation
plays no role in this section. We start as usual with the 
$N_A(x,0-)$ as i.i.d. mean $\mu_A$ Poisson variables. Sometimes we
will add one $A$-particle at the origin at time 0. Thus 
$N_A(x,0) = N_A(x,0-)$ or $N_A(x,0) = N_A(x,0-) + \de(x,\bold 0)$. 
A {\it J-path} is a space-time path $\wh \pi:[0,\infty)
\to \Bbb Z^d \times \Bbb R_+$ such that at all times $t \ge 0, 
\wh \pi(t)$ is the
space-time position of some $A$-particle and such that each jump in
$\wh \pi$ coincides with a jump of some $A$-particle. Thus, such
a path at all times follows an $A$-particle. It may switch from
following one particle to following another particle $\rho$ at any
time when it is at the same space-time point as $\rho$.
(B. Toth suggested that one should think of the $A$-particles as
horses; the path always rides some horse, but may change from one
horse to another when the two horses are at the same place at the same
time.) The
`$J$' in the designation of these paths is to indicate the importance
of the jumps. In fact, we are interested in the following random
variables:
$$
j(t, \wh \pi) := \text{ number of jumps of $\wh \pi$ during
$[0,t]$},
\teq(4.9)
$$
and
$$
J(t,x) := \sup\{j(t, \wh \pi): \wh \pi \text{ is a $J$-path with }
\wh \pi(0) = (x,0)\}.
\teq(4.10)
$$
If $x$ is unoccupied at time 0 (i.e., $N_A(x,0) = 0$), then we take
$J(t,x) \equiv 0$. 
In this section we shall show that $J(t, x)$ is $O(t)$ a.s. We note
that this is obvious in the discrete time setting. The problem only
arises in continuous time and we have only found a quite elaborate
proof of this result.
\proclaim{Proposition 13} There exists a constant $C_{12} < \infty$
such that, for each fixed $x$,
$$
\limsup_{t \to \infty} \frac 1t J(t,x)   \le C_{12} \text{ a.s.}
\teq(4.1)
$$
\endproclaim
We present the proof in a string of lemmas. First we recall some
notation and results of \cite {KSa}. $C_0$ is a large integer chosen as
in \cite {KSa} (6.3)-(6.5) and, as before, $\De_r = C_0^{6r}$.
As in \equ(3.68) and \equ(3.68aa) we define the $r$-block
$$
\cB_r(\bold i,k) := \prod_{s=1}^d[i(s)\De_r, (i(s)+1)\De_r) \times
[k\De_r, (k+1)\De_r).
$$
We further define
$$
V_r(\bold i) := \prod_{s=1}^d [(i(s) -3)\De_r, (i(s) + 4)\De_r),
$$
and the pedestal of $\cB_r(\bold i,k)$ is then
$$
\cV_r(\bold i,k) = V_r(\bold i) \times \{(k-1)\De_r\}.
$$
$\cQ_r(x)$ and $U_r(x,v)$ are as defined in \equ(3.20ab) and \equ(3.20bc).
\comment
For $\bold i \in \Bbb Z^d, k \in
\Bbb Z$ we define the $r$-blocks
$$
\cB_r(\bold i, k) := \prod_{s=1}^d [i(s)\De_r, (i(s)+1)\De_r) \times
[k\De_r, (k+1)\De_r).
$$
The {\it pedestal} of $\cB_r(\bold i,k)$ is 
$$
\cV_r(\bold i,k) := \prod_{s=1}^d [((i(s) - 3)\De_r, (i(s) +4)\De_r) 
\times \{(k-1)\De_r\}.
$$
$$
\cV_r(\bold i, k) = V_r(\bold i) \times \{(k-1)\De_r\},
$$ 
where
For $x \in \Bbb Z^d$ we define further
$$
\cQ_r(x) = \prod_{s=1}^d[x(s), x(s)+C_0^r)
$$ 
$$
U_r(x,v)= \sum_{y \in \cQ_r(x)} N_A(y,v) = \sum \Sb y:x(s) \le y <
x(s)+C_0^r\\ 1 \le s \le d \endSb  N_A(y,v).
\teq(4.3)
$$
\endcomment
For want of a better term, we shall talk about good blocks and good
pedestals. However, the term `good' here does not have the same meaning as in
the good bottoms, good blocks and good pedestals used in Section 3. In
Section 3 a good object contained `many' particles, whereas here a
good object will be one containing `few' particles. 
Since the definitions of Section 3 will not
be used further in this paper we hope that this does not lead to confusion.
Also the constants $C_0,\ga_i$ will be as in (6.2)-(6.5) and (5.10) of \cite
{KSa} (rather than as in \cite {KSb}). The only property of them which
is important here is that the $\ga_r$ now are decreasing in $r$ and that 
$$
0 < \ga_\infty \le \ga _r \le \ga_0.
\teq(4.4)
$$
The $r$-block $\cB_r(\bold i,k)$ is called {\it good} if 
$$
\align
U_r(x,v) \le \ga_r\mu_AC_0^{dr} &\text{ for all $(x,v)$ for which}\\  
&\cQ_r(x) \subset V_r(\bold i) \text{ and } v \in [(k-1)\De_r,
(k+1)\De_r).
\endalign
$$
Similarly, the pedestal $\cV_r(\bold i,k)$ is called good, if
$$
U_r(x,(k-1)\De_r) \le \ga_r\mu_AC_0^{dr} \text{ for all $x$ for which  }
\cQ_r(x) \subset V_r(\bold i).
$$
A {\it bad} block or pedestal is one which is not good. 
Finally, 
$$
\align
\phi_r(\wh \pi):= &\text{number of bad $r$-blocks which intersect 
$\wh \pi|_{[0,t]}$,}\\
&\text{the restriction of $\wh \pi$ to $[0,t]$}.
\endalign
$$
and $\Phi_r(\l)$ and $\Xi(\l,t)$  are exactly as
in \equ(3.75)-\equ(3.75b). These definitions and notations all agree
with \cite {KSa}.
\comment
IS THIS USED ??????????????????
The following result is in Lemma 15 of \cite {KSa}.
There exist some $t_0$ and constants $K_{13}, K_{14}$ such that for
all $t \ge t_0$
$$
P\{\Phi_r(\l) \ge K_{14}\ka_0(t+\l)\exp[-K_{13}C_0^{r/4}] \text{ for
some }r \ge d, \l \ge 0\} \le \frac 2{t^2}.
\teq(4.5)
$$
$\ka_0$ is a strictly positive, finite 
constant introduced in Lemma 11 of \cite {KSa}, whose
value is of no importance here. 

Strictly speaking, the result of \cite {KSa}
is not directly applicable, because the size of the pedestal and
definition of a good block here differ from those in \cite {KSa} 
in the size of some blocks. However, the proofs in \cite {KSa} apply
essentially without change with the present choice of the pedestal and
a good block.  
\endcomment

For simplicity we shall think of $t$ as fixed and 
abbreviate  $\wh \pi|_{[0,t]}$ to $\wh \pi$ if it is clear that only the
restriction of $\wh \pi$ to $[0,t]$ plays a role. $(\bold i,k) \equiv
(\bold a,b)$ for $\bold a \in \{0, 1, \dots,11\}^d$ and $b = 0$ or $1$
will mean that $i(s) \equiv a(s) \pmod {12}$ and $k \equiv b \pmod 2$.

We shall bound $j(t, \wh \pi)$  by a number of sums of the form
$$
\sum_{r\ge 1} \sum_{(\bold i,k)}{}^{(\wh \pi,r)} M(r, \bold i,k) 
I(\wh \pi,r, \bold i,k),
\teq(4.11)
$$
where $\sum_{(\bold i,k)}^{(\wh \pi,r)}$ is a sum over all $(\bold i,k)$
for which $\cB_r(\bold i,k)$ is a good $r$-block which intersects $\wh
\pi |_{[0,t]}$; $I(\wh \pi, r, \bold i,k)$ is the indicator 
function of some
event, and several different choices for $M$ will be made below. 
Let $C_1$ be the constant in Theorem 1 in \cite {KSb} and let
$H_1, H_2$ be the events
$$
H_1(t) := \{\text{all J-paths starting at $(\bold 0,0)$ stay
 in $\cC(C_1t)$ during $[0,t]$}\}
$$
and, for $\wh \pi \in \Xi(\l,t)$,
$$
H_2(\wh \pi, r) := \big\{\sum_{(\bold i,k)}{}^{(\wh \pi,r)} 
I(\wh \pi,r, \bold i,k)
\le \ep_{r} (t+\l)\big\}
\teq(4.12)
$$
for some small $\ep_{r}$. As we shall see, we shall be able to get a
bound on $P\{[H_1 \cap H_2]^c\}$ in several cases. Finally,
it will be the case in our applications that
$$
\align
&\text{for fixed $r \ge 1$ and a 
nonrandom collection $\cS(\bold a, b)$ of $(\bold
i,k)$}\\ 
&\text{with $(\bold i,k) \equiv (\bold a,b)$,
there exists a collection of independent}\\
&\text{random variables
$\{\wt M(r,\bold i,k): (\bold i,k) \in \cS(\bold a,b)\}$ which satisfy}\\
&\wt M(r, \bold i,k) \ge |M(r, \bold i, k)|
\text{ and $E\exp[\th_r\wt M(r, \bold i,k)] \le \Ga_r$ for}\\ 
&\text{some constants $\th_r >0$ and $1 \le \Ga_r < \infty$}.
\teq(4.13)
\endalign
$$
Lemma 15 shows how to handle such a sum, but first we need some
information on the location of J-paths starting at the origin.

Even though $J(t,x)$ does not
involve $B$-particles, we shall make use of $B$-particles 
in the proof of Lemma 20. Also in the proof of that lemma, we shall
need to consider initial conditions which are not of the form of
i.i.d. Poisson variables $N_A(x,0-)$ plus some extra particles at time
0. We therefore do not make this assumption in the next lemma.
In particular, we only assume that the
$\{Y_t\}$-process (which has no recuperation) 
is formed by adding one $B$-particle at the origin
at time 0, and that this process is coupled with the $A$-system by
giving the same path to each $A$-particle present at time $0-$ 
in this process as in the
$A$-system. (This is exactly as in Section 2.) 
In addition, the initial $N_A(x,0-)$ have to be such that $Y_0
\in \Si_0$ a.s. $\Si_0$ is the state space introduced in \cite {KSb},
\cite {KSc}.
All particles still perform independent continuous time
simple random walks. 
 
\proclaim{Lemma 14} Under the conditions just described
$$
\align 
&\{\text{there is a J-path from $(\bold 0,0)$ through $(x,s)$}\}\\
&\subset\{\text{there is a $B$-particle at $(x,s)$ in the $\{Y_t\}$-process}\}.
\teq(4.16)
\endalign
$$
In particular, 
$$
\align
&P\{[H_1(s)]^c\} = 
P\{\text{some J-path starting at $(\bold 0,0)$ leaves $\cC(C_1s)$ during
$[0,s]$}\}\\
&\le 2P\{\text{some $J$-path which started 
at $(\bold 0,0)$} \text{ is outside $\cC(C_1s)$
at time $s$}\}\\
&\le 2E\{(\text{number of $B$-particles outside $\cC(C_1s)$ at time
$s$, in the $\{Y_t\}$-process})\}.
\teq(4.17)
\endalign
$$
\endproclaim
\demo{Proof} Clearly adding an $A$-particle to the $A$-system can only
increase the collection of J-paths, so that we may assume that we
start the $A$-system with $N_A(x,0) = N_A(x,0-) + \de(x,\bold 0)$. 
(We repeat that the $N_A(x,0-)$ do not have to be i.i.d. Poisson
variables in this lemma.) We
can then couple the $A$-system and the $\{Y_t\}$-process so that they
have the same particles and so that each particle follows the same
path in both processes. The only difference between the proceeses is
that in the $A$-system all particles have type $A$, while in
$\{Y_t\}$ there are particles of both types.

Now assume that there is a J-path $\wh \pi$ in the $A$-system from
$(\bold 0,0)$ to $(x,s)$. Then there exists some sequence of times
$s_0=0 < s_1 < \dots < s_\l < s$ and particles $\rho_i$ such that 
$\wh \pi$ agrees with the path of $\rho_i$ during $[s_i,
s_{i+1}], \; 0\le i \le \l$ (with $s_{\l+1} = s$). In addition
$\rho_{i+1}$ and $\rho_i$ are at the same position at time $s_i$,
while $\rho_0$ starts at $(\bold 0,0)$ and $\rho_\l$ is at $x$ at time
$s$. Now in the $\{Y_t\}$-process all particles at $\bold 0$ are
given type $B$ at time $s_0=0$. But then $\rho_0$ has type $B$ for all
$t \ge 0$. Then $\rho_1$ will have type $B$ at least from time $s_1$
on. One then sees by induction on $i$ that $\rho_{i+1}$ has type $B$
on $[s_i,\infty)$. In particular, $\rho_\l$ has type $B$ at the time $s >
s_\l$, at which time it is at $x$. This implies \equ(4.16).

Next, we have in the $A$-system 
$$
\align
&P\{\text{some J-path starting at $(\bold 0,0)$ leaves $\cC(C_1s)$ during
$[0,s]$}\}\\
&\le 2P\{\text{some J-path starting at $(\bold 0,0)$
is outside $\cC(C_1s)$ at time $s$}\}.
\teq(4.17z)
\endalign
$$
This follows from a reflection argument, as in  the proof of
Proposition 3 in \cite {KSb}. The last inequality in \equ(4.17) 
then follows from \equ(4.16).
\hfill $\blacksquare $
\enddemo

We now return to the usual initial conditions, that is we take the
$\{N_A(x,0-): x \in \Bbb Z^d\}$ as i.i.d., mean $\mu_A$ Poisson variables.
We also add an extra particle to the system at the origin at time
0. We note that Proposition 5 and Remark 2 after it in \cite {KSb}
show that in this case $Y_0 \in\Si_0$ a.s., so that we can apply Lemma
14 in this case. If the $N_A(x,0-)$ are i.i.d. Poisson variables, then
\equ(4.17), together with (1.3) in \cite {KSb}, shows that
for all large $t$
$$
P\{\text{some J-path starting at $(\bold 0,0)$ leaves $\cC(C_1t)$ during
$[0,t]$}\} \le 2e^{-t}.
\teq(4.17aa)
$$
This will allow us to restrict our further estimates to J-paths
$\wh \pi$ which stay in $\cC(C_1t)$ during $[0,t]$. If we take $t$ so
large that $C_1t \le t \log t$, then these paths also stay in $\cC(t
\log t)$ and therefore belong to $\cup_{\l \ge 0} \Xi(\l,t)$ (see
\equ(3.75b) for $\Xi$). This explains why the next few lemmas speak
about such paths only. In fact, it is useful to make a further
reduction. To this end, we define, as in (6.10) of \cite {KSa},
$R=R(t)$ as the integer for which
$$
C_0^R \ge [\log t]^{1/d} > C_0^{R-1}.
\teq(4.19hij)
$$
We then have just as in Lemma 9 of \cite {KSa} that
$$
P\{\text{for some $r \ge R$ there exists a bad $r$-block which intersects
$\cC(t \log t)$}\} \le \frac 1{t^2}.
\teq(4.22)
$$
Accordingly we define the event
$$
\Th(t) = \{\text{for all $r \ge R(t)$ no bad $r$-block 
intersects }\cC(t\log t)\}.
\teq(4.13aa)
$$
We then have 
$$
P\{[\Th(t)]^c\} \le \frac 1{t^2},
\teq(4.13bb)
$$
so that we can restrict most estimates to subevents of $\Th(t)$.
Since we are only concerned with the existence of certain $J$-paths it
is convenient to define
$$
\Xi(J,\l,t):= \text{ the collection of $J$-paths in }\Xi(\l,t).
\teq(4.13bcb)
$$

Many constants $K_i, R_i$ and $t_i$ appear in the remainder of this
section. It is crucial that these {\it do not depend on $t$ or $\l$},
even though we usually do not state this explicitly. 

\proclaim{Lemma 15} Assume that \equ(4.13) is satisfied. Then,
there exist constants $0< K_1-K_3
< \infty$, all depending on $d$ only, 
such that for $t \ge 2$ and for each $\l \ge 0$,
$$
\align
&P\{\text{$\Th(t)$ and for some } \wh \pi \in \Xi(J,\l,t),\;
\sum_{(\bold i,k)}{}^{(\wh \pi,r)} M(r, \bold i,k) I(\wh \pi, r, \bold
 i,k) \ge x\}\\
&\le P\{\text{for some }\wh \pi \in \Xi(J,\l,t),
\;\Th(t) \cap [H_2(\wh \pi,r)]^c \text{ occurs}\}\\
&\qquad + K_1[t \log t]^d\exp\;[K_2(t+\l)/\De_{r}]\exp\big[-\frac{x\th_r}{2
(12)^d} +\ep_r(t+\l)\log \Ga_r\big].
\teq(4.14)
\endalign
$$
For 
$$
x \ge \frac{4(12)^d}{\th_r}\big[\ep_{r}\log \Ga_r +
\frac{K_2}{\De_{r}}\big] (t+\l)
\teq(4.20)
$$
this yields
$$
\align
&P\{\Th(t) \text{ and for some } \wh \pi \in \Xi(J,\l,t),\;
\sum_{(\bold i,k)}{}^{(\wh \pi,r)} M(r, \bold i,k) I(\wh \pi, r, \bold
 i,k) \ge x\}\\
&\le P\{\text{for some }\wh \pi \in \Xi(J,\l,t),
\; \Th(t) \cap [H_2(\wh \pi,r)]^c\text{ occurs}\} 
+ K_1 [t\log t]^d\exp[-K_3x\th_r].
\teq(4.21)
\endalign
$$
\endproclaim
\demo{Proof}The first term in the right hand side of \equ(4.14) takes care 
of the event that $H_2(\wh \pi,r)$ fails for any $\wh \pi$.
It therefore suffices for \equ(4.14) to estimate
$$
\align
P\{\text{for some }\wh \pi \in \Xi(J,\l,t),\; &H_2(\wh\pi,r) 
\text{ occurs and }\\
&\sum_{(\bold i,k)}{}^{(\wh \pi,r)} M(r, \bold i,k) I(\wh \pi, r, \bold
 i,k) \ge x\}.
\teq(4.18)
\endalign
$$
If the event here occurs, then there is a $\wh \pi \in \Xi(\l,t)$  and
a subset, $\cS$ say, of the points $(\bold i,k)$ 
for which $\cB_{r}(\bold i,k)$ intersects $\wh \pi$, such that 
$I(\wh \pi,r,\bold i,k)=1$ for $(\bold i,k) \in \cS$ and 
$$
\sum_{(\bold i,k) \in \cS} |M(r,i,k)| \ge x.
$$
Moreover, $\cS$ contains at most $\ep_{r}(t+\l)$ points (because
$H_2(\wh \pi,r)$ occurs). 
We can split $\cS$ into the $2(12)^d$ subsets 
$$
\cS(\bold a, b) = \text{collection of $(\bold i,k)$ in $\cS$ 
with $(\bold i,k)
\equiv (\bold a,b)$},
$$
with $\bold a \in \{0,1,\dots,11\}^d, b = 0$ or $1$.
Then \equ(4.18) is bounded by the sum of
$$
P\Big\{\sum_{(\bold i,k) \in \cS(\bold a,b)} M(r,i,k) 
\ge \frac x{2(12)^d}\Big\}
\teq(4.19)
$$
over all possible $\cS(\bold a,b)$ corresponding to
some $\wh \pi \in \Xi(\l,t)$.

We know that any $\wh \pi \in
\Xi(\l,t)$ intersects at most 
$$
\la_r(\l) := 3^d\Big(\frac {t+\l}{\De_{r}} + 2\Big) 
\teq(4.19abc)
$$
$r$-blocks (see (6.30) in \cite {KSa} for $\nu = 1$ and with $r+1$
replaced by $r$). The set of $(\bold i,k)$ for which $\cB_r(\bold
i,k)$ intersects $\wh \pi$ has to be $\cL$-connected (see the lines
following \equ(3.2) for $\cL$).
Thus, as $\wh \pi$ varies over $\Xi(\l,t)$, and the starting point of
$\wh \pi$ varies over $\cC(t \log t)$, there are at most 
$[2t\log t +1]^d \exp[K_9\la_r(\l)]$
different possibilities for the collections 
$\{(\bold i,k):\cB_{r}(\bold i,k) \text{ intersects }\wh \pi\}$.
Here $K_9$ is some constant which depends
on $d$ only. 
Each $\cS$ has to be a subset of this collection, 
and once $\cS$ is given there are $2\cdot 12^d$ possibilities for
$(\bold a, b)$. Thus, there are at most
$$
2(12)^d2^{\la_r(\l)}[2t\log t+1]^d \exp[K_9\la_r(\l)]
\teq(4.19efg)
$$
possibilities for $\cS(\bold a,b)$.

Finally, for a fixed choice of $\cS(\bold a,b)$ we have by \equ(4.13) that
the probability in \equ(4.19) is bounded by
$$
\align
&P\Big\{\sum_{(\bold i,k) \in \cS(\bold a,b)} \wt M(r,\bold i,k) 
\ge \frac x{2(12)^d}\Big\}\\
&\le \exp\big[-\frac {x\th_r}{2(12)^d}\big]E\big\{\exp
\big[\th_r\sum_{(\bold i,k) \in
\cS(\bold a,b)} \wt M(r,\bold i,k)\big]\\
&\le \exp\big[-\frac {x\th_r}{2(12)^d}\big] \Ga_r^{\ep_r(t+\l)},
\teq(4.19cde)
\endalign
$$
where, for the last inequality, 
we used that $\cS(\bold a,b)$ has at most $\ep_r(t+\l)$ elements
on $H_2(\wh \pi,r)$ and \equ(4.13). This implies \equ(4.14) for
suitable $K_1,K_2$, because
\equ(4.18) is bounded by a sum of at most \equ(4.19efg) terms, 
each of which is bounded by \equ(4.19cde).

The inequality \equ(4.21) now follows from \equ(4.14) and 
\equ(4.20), because the latter implies that
$$
\frac {x\th_r}{2(12)^d} \ge \frac {x\th_r}{4(12)^d}
+ \big [\ep_r \log \Ga_r + \frac {K_2}{\De_r}\big](t+\l).
\tag "$\blacksquare$"
$$
\enddemo

Our first task is now to establish a representation for $j(t, \wh
\pi)$ of the form \equ(4.11), at least outside an event of small
probability.
Fix some $R_1 \ge 1$ and consider a sample point for which $\Th(t)$
occurs. If
$\wh \pi$ is a $J$-path which stays in $\cC(t \log t)$ during $[0,t]$,
then all $r$-blocks with $r \ge R$ which intersect $\wh \pi|_{[0,t]}$ 
must be good. Now recall that for each $r$, 
each point of space-time belongs to a
unique $r$-block $\cB_r(\bold i,k)$. 
We shall say that a jump in $\wh \pi$ from $x$ to
$y$ at time $s$ is {\it located at} $(x,s)$. 
For such a jump, either $(x,s)$ belongs to a good
$r$-block for all $r \ge R_1$, or there is a unique $r(x,s) \in (R_1,R]$
such that $(x,s)$ belongs to a good $r$-block for $r \ge r(x,s)$, but
belongs to a bad $[r(x,s)-1]$-block. In the former case we set $r(x,s)
= R_1$. Note that for any jump $(x,s)$, $r(x,s)$ is defined and the jump
lies in some $r(x,s)$-block. Moreover, this is a good block, by the
choice of $r(x,s)$.
Since we also know from \equ(4.17aa) that $J$-paths which start at
$\bold 0$ stay inside
$\cC(C_1t)$ during $[0,t]$, except for an exponentially 
small probability, we have for $t \ge$ some $t_1$ that 
outside the event in \equ(4.17aa) but in $\Th(t)$, it holds
$$
\align
J(t, \bold 0) = \sup_{\wh \pi(0) = \bold 0} 
\sum_{r = R_1}^{R(t)}&\text{[number of jumps $(x,s)$ of $\wh \pi$
with}\\ 
&\text{$(x,s)  \in \cC(t\log t) \times [0,t]$
and $r(x,s) = r$]}
\teq(4.23)
\endalign
$$
(the sup here is over the same set as in \equ(4.10)). The union of
the exceptional event in \equ (4.17aa) and $[\Th(t)]^c$ 
has probability at most $2/t^2$.

We now concentrate on estimating the summands appearing in the right
hand side in \equ(4.23). Let $(x,s) \in \cB_r(\bold i,k)$ be a jump
of $\wh \pi$. This jump is the jump of some
 particle $\rho$ at time $s$. We distinguish two kinds of jumps,
according as $\rho$ was outside or inside the pedestal $\cV_r(\bold
i,k)$ at time $(k-1)\De_r$. We define the corresponding quantity
$$
\align
M_{\text{out}}(r, \bold i, k) = [&\text{number of jumps inside 
$\cB_r(\bold i,k)$ by any particle $\rho$}\\
&\text{that was outside $\cV_r(\bold
i,k)$ at time $(k-1)\De_r$]},     
\endalign
$$
and its analogue $M_{\text{in}}(r,\bold i,k)$ with ``outside''
replaced by ``inside''. We further say that the block
$\cB_r(\bold i,k)$ is {\it contaminated} if it contains a jump of a
particle which was outside $\cV(\bold i,k)$ at time $(k-1)\De_r$ 
and take
$$
I_1(r,\bold i,k) := I[\cB_r(\bold i,k) \text{ is contaminated}].
$$
We point out that 
this definition of contaminated is somewhat stricter than the one
used in \cite {KSa} (just after (6.9)). 

We now start with a bound for 
$$
\sum_{(\bold i, k)}{}^{(\wh \pi,r)} 
M_{\text{out}}(r,\bold i,k) I_1(r, \bold i,k).
\teq(4.24)
$$
\proclaim{Lemma 16}
There exist constants $K_i, t_2$ and $R_2$ such that
for $t \ge t_2, R_2 \le r \le R(t)$ and $\l \ge 0$, 
$$
\align
&P\big\{\text{\rm there exists a }\wh \pi\in \Xi(J,\l,t) 
\text{ \rm such that }\\
&\phantom{MMMMMMMMMM}\sum_{(\bold i, k)}{}^{(\wh \pi,r)} 
M_{\text{out}}(r,\bold i,k) I_1(r, \bold i,k) \ge \frac{K_4 
(t+\l)}{ \De_r}\big\}\\
& \le K_5\exp\big[-K_6(t+ \l)/[\log t]^6\big].
\teq(4.29)
\endalign
$$
\endproclaim
\demo{Proof} We break the proof up into two steps.
\medskip
Step (i). In this step we reduce the calculations to some calculations
for discrete time random walks. This first step is standard weak
convergence theory and we leave many details to the reader.
We approximate the paths of the various particles by some
random walk paths which can jump only at times $j/n$ for some integer
$n \ge 1$ and $j = 0, 1, 2,\dots$
Specifically, we let $\{S_u^{(n)}\}_{u \ge 0}$ be a random walk starting 
at $\bold 0$ which can jump only at times $j/n$, with the
jump distribution
$$
q^{(n)}(y) = P\{S^{(n)}_{(j+1)/n} - S^{(n)}_{j/n} = y\} =
\cases 1-\frac Dn &\text{ if } y = \bold 0\\\frac D{2dn}
&\text{ if }y = \pm e_i, \; 1 \le i \le d,
\endcases
$$
($e_i$ is the $i$-th coordinate vector). For each particle $\rho$ we 
take $\{S^{(n)}_u(\rho)\}_u \ge 0$ as a copy of $\{S_u^{(n)}\}_{u \ge 0}$, and
we take the walks for the different $\rho$ as completely independent.
We then form what we shall call the $(n)$-system by letting 
$\rho$ move along the
path $t \mapsto \pi^{(n)}(t,\rho) :=\pi(0,\rho) + 
S^{(n)}_{\lfloor tn \rfloor}(\rho)$ for each of
the particles $\rho$.
Now it easy to see that for any finite collection of particles
$\big(\rho_{i_1}, \dots \rho_{i_K}\big)$, the $K$ dimensional process 
$t \mapsto \big(\pi^{(n)}(t, \rho_{i_1}, \dots, \pi^{(n)}(t, \rho_{i_K}\big)$
converges weakly (in the Skorokhod topology on the space
$D\big([0,\infty), (\Bbb Z^d)^K\big)$ to the process 
$t \mapsto \big(\pi(t, \rho_{i_1}, \dots, \pi(t, \rho_{i_K}\big)$.
This last process is the process of the true paths of $(\rho_{i_1},
\dots, \rho_{i_K})$. A simple way to prove this weak convergence 
is to apply Theorem 15.6 in \cite {B} (or rather the line following it
before the proof of Theorem 15.6). We then define the obvious analogue
of $N^*$, namely
$$
\align
N^{(n)}(x,t) &= \text{(number of particles at $(x,t)$ in 
the $(n)$-system)}\\
&= \text{(number of $\rho$ with $\pi^{(n)}(t, \rho) = x$)}.
\endalign
$$
Here we do not include the extra particle added at the origin at time
0; we only include the particles which were among the $N_A(x,0-)$ at
some $x$, just before the start of our system.
We also need an approximation to $N^{(n)}$ which only counts particles
which started in some finite cube. For this we fix some numbering of 
the particles $\rho_1, \rho_2, \dots$. Again, this excludes the extra 
particle added at $(\bold 0,0)$ if there is such a particle.
We then set 
$$
N^{(n,L)}(x,t) := \text{ (number of $i \le L$
with $\pi^{(n)}(t, \rho) = x$)}.
$$
It is convenient to set $N^{(\infty)}(x,t) = N^*(x,t)$ and similarly
$$
N^{(\infty,L)} =  \text{(number of $\rho$ among the first $L$ particles with
$\pi(t, \rho) = x)$}.
$$
Particles which start far out only have a small probability of 
reaching  $\cC(2t \log t)$ during $[0,2t]$. In fact, estimates like the
ones for (2.29)-(2.32) in \cite {KSb}
prove that for all $t > 0$ and $\eta > 0$ there exists an $L_0 = L_0(t,\eta)$
such that
$$
\align
P\{N^{(n,L)}(x,s) \ne N^{(n)}&(x,s) \text{ for any } (x,s) \in
\cC(2t\log t) \times [0,2t]\} \le \eta, \\
&\text{for all }L \ge L_0, 1 \le n \le \infty.
\teq(4.23aa)
\endalign
$$ 
Note the uniformity in $n$ here.
We can now replace $N^*$ by $N^{(n)}$ or $N^{(n,L)}$ in many of the
definitions. We indicate such a replacement by decorating the
appropriate quantity with a superscript $(n)$ or $(n,L)$ in a
self explanatory fashion, or by adding the qualification ``in the $(n)$
system or $(n,L)$-system.'' For instance,
$$
U^{(n,L)}_r(x,v) :=   \sum_{y \in \cQ_r(x)} N^{(n,L)}(x,v)
$$
and the block $\cB_r(\bold i,k)$ is good in the $(n,L)$-system if 
$$
\align
U^{(n,L)}_r(x,v) \le \ga_r\mu_AC_0^{dr} &\text{ for all $(x,v)$ for which}\\  
&\cQ_r(x) \subset V_r(\bold i) \text{ and } v \in [(k-1)\De_r,
(k+1)\De_r).
\endalign
$$
\equ(4.23aa) immediately implies that uniformly in $1 \le n \le \infty$
$$
P\Big\{\sum_{(\bold i, k)}{}^{(\wh \pi,r)} 
M_{\text{out}}^{(n,L)}(r,\bold i,k) I_1^{(n,L)}(r, \bold i,k) \ne
\sum_{(\bold i, k)}{}^{(\wh \pi,r)} 
M^{(n)}_{\text{out}}(r,\bold i,k) I^{(n)}_1(r, \bold i,k)\Big\} \le \eta
\teq(4.23bb)
$$
for $L \ge L_0(t,\eta)$, provided $t$ is so large that any $r$-block
which intersects $\cC(t\log t) \times [0,t]$ is contained in $\cC(2t
\log t) \times [0,2t]$. (It suffices for this last proviso that $3\De_r
= 3 C_0^{6r} \le t$.)

Next we claim that for each fixed finite $L, r$ and fixed finite set
$\cS$ of pairs $(\bold i,k)$, as $n \to \infty$,
$$
M^{(n)}_{\text{out}}(r,\bold i,k) 
\text{ converges weakly to }
M_{\text{out}}(r,\bold i,k) 
\teq(4.23ff)
$$
and
$$
\align
&\sum_{(\bold i, k)\in \cS}
M^{(n)}_{\text{out}}(r,\bold i,k) I^{(n)}_1(r, \bold i,k)
\text{ converges weakly}\\
&\phantom{MMMMMMMMM}\text{ to }\sum_{(\bold i, k) \in \cS}
M_{\text{out}}(r,\bold i,k) I_1(r, \bold i,k).
\teq(4.23cc)
\endalign
$$
This is an immediate application of the continuous mapping theorem (Theorem 5.1
or 5.2) in \cite {B}). Indeed, in any system of $L$ moving particles
with joint paths $s \mapsto (\pi_1(s), \dots ,\pi_L(s)) \in (\Bbb
Z^d)^L$ we can define $\overline U_r(x,v)$ for $v \le t$ 
as a functional of these paths by 
$$
\overline U_r(x,v) = \text{ (number of $i \in [1,L]$ with $\pi_i(v)
\in \cQ_r(x)$)}.
$$
We restrict ourselves to paths $\pi_i$ which are right continuous with
left limits, so that we view $\overline U_r(x,v)$ as a functional on 
the Skorokhod space $D([0,t], (\Bbb Z^d)^L)$, and we put the Skorokhod
topology on this space. Then $U^{(n,L)}(x,v)$ is just the value of 
$\overline U_r(x,v)$ at the point with $\pi_i(\cdot) = \pi^{(n,L)}(\cdot,
\rho_i)$. In a similar way we can view $I[\cB_r(\bold i,k)$ is good], 
$M_{\text{out}}(r,\bold i,k)$
and $I_1(r,\bold i,k)$ as the value at $\pi_i(\cdot) = \pi^{(n,L)}(\cdot,
\rho_i)$ of suitable functionals on $D([0,t], (\Bbb Z^d)^L)$. We
indicate these functionals on $D([0,t], (\Bbb Z^d)^L)$ by a bar over
the appropriate symbol.
Now it is not hard to see that $\overline U_r(x,v)$ and
$$
\align
&\overline I[\cB_r(\bold i,k) \text{ is good}]\\ 
&= I\big[\sup\{\overline U_r(x,v):
\cQ_r(x) \subset V_r(\bold i),(k-1)\De_r \le v < (k+1)\De_r\} \le
\ga_r\mu_A C_0^{dr}\big]
\endalign
$$
are continuous functionals on
$D([0,t], (\Bbb Z^d)^L)$ at all points 
$\big(\pi_1(\cdot), \dots, \pi_L(\cdot)\big)$
for which each $\pi_i$ is continuous at each
$\{j\De_r: j \in \Bbb Z\}$. 
In particular, this holds almost surely at the points with
$\pi_i(\cdot) = \pi(\cdot, \rho_i)$.
Similarly $\overline M_{\text{out}}(r, \bold i,k)$ is continuous 
at these same points.
Therefore, \equ(4.23cc) does indeed follow from the continuous mapping
theorem.

Finally we note that the event in the left hand side of \equ(4.29)
occurs if and only if
$$
\sum_{(\bold i,k) \in \cS}M_{\text{out}}(r, \bold i,k)  I_1(r, \bold i,k) \ge 
 \frac{K_4 (t+\l)}{ \De_r},
\teq(4.23ee)
$$  
for one of a number of possible collections $\cS$ of pairs $(\bold
 i,k)$. The possible 
collections $\cS$ are the collections of 
the form $\{(\bold i,k):\cB_r(\bold i,k)$ is good and intersects $\wh
\pi$\}, for some $\wh \pi \in \Xi(J,\l,t)$. The number of possibilities
for $\cS$ is finite, and whether $\cS$ is a possible collection depends
on the class $\Xi(J,\l,t)$ and on which blocks $\cB_r(\bold i,k)$ are
good. The indicator function of $\{\cS$ is
possible collection\} for a fixed collection $\cS$ is also a
continuous functional on $D([0,t], (\Bbb Z^d)^L)$. We can now
combine this observation with \equ(4.23bb) and \equ(4.23cc) 
to obtain the conclusion of
this step that the left hand side of \equ(4.29) is bounded by
 $$
\align
&\limsup_{L \to \infty} \limn 
P\{\text{there exists a }\wh \pi\in \Xi(J,\l,t) 
\text{ such that}\\
&\phantom{MMMM}\sum_{(\bold i, k)}{}^{(\wh \pi,r)} 
M^{(n,L)}_{\text{out}}(r,\bold i,k) I^{(n,L)}_1(r, \bold i,k) \ge \frac{K_4 
(t+\l)}{ \De_r}-\big\}.
\endalign
$$
In fact, since the collection of particles present in the
$(n,L)$-system increases to the collection of particles in the
$(n)$-system as $L \to \infty$, this expression is bounded by
$$
\align
&\limn 
P\{\text{there exists a }\wh \pi\in \Xi(J,\l,t) 
\text{ such that}\\
&\phantom{MM}\sum_{(\bold i, k)}{}^{(\wh \pi,r)} 
M^{(n)}_{\text{out}}(r,\bold i,k) I^{(n)}_1(r, \bold i,k) \ge \frac{K_4 
(t+\l)}{ \De_r}-\big\}.
\teq(4.23dd)
\endalign
$$

\medskip
Step (ii). In this step we derive a bound for \equ(4.23dd) in terms of
a large number of independent copies of the random walk
$\{S^{(n)}_u\}_{u \ge 0}$. We
follow the proof of Lemmas 10 and 11 in \cite {KSa} closely.

We take for 
$\{S^{(n)}_u(x,s,q)\}_{u \ge 0}$ 
a copy of $\{S^{(n)}_u\}_{u \ge 0}$ and take all these copies for different $x
\in \Bbb Z^d, s$ of the form $k/n$ and $q \ge 1$, completely 
independent. We further associate to each particle $\rho$ a 
uniform random variable
on $[0,1]$, $U(\rho)$ say, and all $U(\rho)$ and $\{S^{(n)}_u(x,s,q)\}$ are 
independent. Finally 
$$
\cW_r (\bold i) := \partial\prod_{s=1}^d [(i(s)
-3)\De_r,(i(s)+4)\De_r-1] = \partial V_r(\bold i),
\teq(4.29z)
$$
where $\partial$ denotes the topological boundary.
We now fix some $\bold a \in \{1, 2, \dots, 11\}^d$ and $b = 0$ or
$1$, and we want to look at the contribution to the sum \equ(4.24)
from the $(\bold i, k) \equiv (\bold a, b)$. For the sake of argument
let $b = 0$.
Assume that the paths of all $A$-particles till time 
$(k-1)\De_r$ with $k$ even have already been constructed in some way. 
In the case  $k=0$ this simply means that we begin 
with a mean $\mu_A$ Poisson system of 
$A$-particles at time $-\De_{r}$.
(The only change which is needed for the case $b=1$ is that we work
with odd $k$'s and start with a Poisson system at time $-2\De_{r}$
in that case.)
At each point $(x,(k-1)\De_{r})$
(in space-time) order all particles $\rho$ present so that their
associated uniform variables $U(\rho)$ are increasing. To the
$q$-th particle in this order associate the path $\{x+S^{(n)}_u(x,
(k-1)\De_{r},q)\}_{u \ge 0}$. This particle then moves to $x+S_1^n(x,
(k-1)\De_{r},q)$ at time $(k-1)\De_{r}+1/n$. We also associate to each
particle at each time an {\it index} $(y',v',q',g')$. A particle has
index $(y',v',q',g')$ at a certain time if its last associated random
walk is $S^{(n)}(y',v',q')$ and if the particle has moved $q'$ steps 
(or $g'/n$ time units) according to $S^{(n)}(y',v',q')$ since this
random walk was associated to the particle. Accordingly, the index
associated to the $q'$-th particle at $(x, (k-1)\De_r)$ at time 
$(k-1)\De_r + 1/n$ is $(x,(k-1)\De_r,q', 1)$.
\comment
The {\it index} associated  
with this particle is then taken to be $(x,(k-1)\De_{r},q, 1)$. The last
coordinate 1 here indicates that one step was taken since the last
choice of an associated random walk.
\endcomment
Assume we have constructed the paths of all particles
up to and including time $v
\in [(k-1)\De_{r}, (k+1)\De_{r})$ (with $v$ a multiple of $1/n$)
and that each particle has an index.
To construct the paths  $1/n$ time units 
further, we look for each $y \in \Bbb Z^d$ at all particles at $(y,v)$.
If $y$ does not belong to 
$$
\bigcup_{\bold j \equiv \bold a} \cW_r(\bold j),
\teq(6.50)
$$
and a particle at $(y,v)$ has index $(z,v',q,g)$, then this particle   
moves to $y +S^{(n)}_{g+1}(z,v',q)$ and its new index is $(z,v',q,g+1)$.
In other words it moves one step further in the random walk it is
presently associated with, and the last component of its index increases
by 1.
\comment
 to indicate the number of steps the particle has taken according
to its present associated random walk. (The first two components $z,v'$
indicate at which space-time point the particle started using this
random walk and the third component specifies which of the random
walks from the point $(z,v')$ the particle is using.)
\endcomment
If, on the other
hand, $y$ lies in the union \equ(6.50), then all particles at $y$
are again ranked according to increasing values of their uniform
random variables and a new random walk is associated to these particles. 
The particle with rank $q'$ will move to 
$y + S_1^n(y,v,q')$ at time $v+1/n$. Its index will then be
$(y,v,q',1)$. We continue this procedure till all positions at 
time $(k+1)\De_{r}$ have been determined. We then start anew with $k$
replaced by $k+1$. That is, we order all particles at one site
$(x,(k+1)\De_{r})$ and move the $q$-th particle at that site to $x+S_1^n(x,
(k+1)\De_{r}),q)$ and give it the index  
$(x,(k+1)\De_{r},q,1)$, and so on.

Basically, the above procedure switches each particle to a new random walk
every time the particle visits the set \equ(6.50).
It is clear that in the above construction all the $A$-particles perform
independent random walks with transition probability $q_A^n$.
\comment
The proof of Lemma 10 of \cite {KSa} shows that distribution of 
the paths of any given
finite collection of particles $\si_1, \dots \si_\l$ in this discrete system
converges (as $n \to \infty$) to the distribution  of the paths of 
these same particles in
the original system with continuous time random walks. Similarly, the
distribution in the discrete time system of the number of jumps of 
$\si_1, \dots, \si_\l$ in $\cB_r(\bold i,k)$ converges to the distribution
of these numbers of jumps in the continuous time system. 
\endcomment
Now, a particle $\rho$ whose jumps contribute to one of the sums \equ(4.23ee)
has to lie outside $\cV_r(\bold i,k)$ at time $(k-1)\De_r$,
but has some jump in $\cB_r(\bold i,k)$ during $[k\De_r,(k+1)\De_r)$.
Its space-time path in the discrete time system must contain a piece 
$(x_\xi,v),
(x_{\xi+1},v+1/n),\dots,(x_\ze, v+(\ze-\xi)/n)$ with $v$ a multiple of
$1/n$, which satisfies 
$$
\align
&x_\xi \in \cW_r(\bold i), x_\ze \in \partial
\prod_{s=1}^d[(i(s)-1)\De_r, (i(s)+2)\De_r -1]\text{ and }
x_\ka \text{ strictly}\\&\text{between }\cW_r(\bold i) \text{ and }
\partial \prod_{s=1}^d [(i(s)-1)\De_r, (i(s)+2)\De_r-1]
\text{ for }\xi < \ka < \ze,
\teq(4.25x)
\endalign
$$
and which is
traversed during $[(k-1)\De_r, (k+1)\De_r)$
(compare (6.17) in \cite {KSa}; note that the condition on $x_\ka$ is
misstated there). At the times $v+j/n, 1 \le j \le (\ze - \xi)/n$,
$\rho$ is at a position in the open cube $\prod_{s=1}^d \big(
(i-3)\De_r,(i+4)\De_r\big)$ and hence does not visit \equ(6.50).
Therefore, the random walk associated to $\rho$ remains the same at
the times $v+j/n, 0 \le j \le (\ze - \xi)/n$.
It follows that for $(x_\xi,v)$ to be the first point of 
such an excursion from $\cW_r(\bold i)$ to
$\partial \prod_{s=1}^d [(i(s)-1)\De_r, (i(s)+2)\De_r-1]$ it is
necessary that for an appropriate $q,\, S_1^{(n)}(x_\xi,v,q) \ne \bold 0$
and $\sup_{u \le 2\De_rn}
\|S_u^{(n)}(x_\xi,v,q)-S_1^{(n)}(x_\xi,v,q)\|
 \ge 2\De_r-1$. 
The last condition has to be satisfied because the minimal distance 
between  $\cW_r(\bold i)$ and 
$\partial \prod_{s=1}^d [(i(s)-1)\De_r, (i(s)+2)\De_r-1]$ is
$2\De_r$,and we are counting jumps in $\cB_r(\bold i,k)$ after time
$(k-1)\De_r$. These jumps must occur in the time interval
$[(k-1)\De_r, (k+1)\De_r]$, i.e., in at most $2\De_rn$ steps of
$S^{(n)}(x_\xi,v,q)$. (Note our terminology here: $S^{(n)}_u$ takes a
step each time $u$ increases by 1, but it has a jump only if
$S^{(n)}_{u+1} \ne S^{(n)}_u$.)
Suppose $S^{(n)}(x_\xi,v,q)$ indeed leaves  $\cW_r(\bold i)$ and reaches
$\partial \prod_{s=1}^d [(i(s)-1)\De_r, (i(s)+2)\De_r-1]$
before it returns to  $\cW_r(\bold i)$. In
this case, let $m = m(x_\xi,v,q)$ be the smallest integer for which 
$S_m^{(n)}(x_\xi,v,q) \in \partial \prod_{s=1}^d [(i(s)-1)\De_r,
(i(s)+2)\De_r-1]$. In the notation of \equ(4.25x), this is the 
number of steps it takes $S^{(n)}(x_\xi,v,q)$ to reach $x_\ze$.
The number of jumps of $\rho$ in $\cB_r(\bold i,k)$ between time $v$
and the next return to $\cW_r(\bold i)$ is then bounded by
the number of jumps 
of $\{S_u^{(n)}(x,v,q)\}$ for $m \le u \le m+2\De_r n$.
This number is independent of all random walks $S^{(n)}(y,w,s)$ with
$(y,w,s) \ne (x_\xi,v,q)$ and of the $S^{(n)}_u(x_\xi,v,q)$ for $u \le m$.
Moreover, if $\cB_r(\bold i,k)$ is good (in the $(n)$-system), 
then there are at most $\ga_r \mu_AC_0^{dr}+1$ particles at the 
space-time point $(x_\xi,v)$. Indeed 
$N^{(n)}(x_\xi,v) \le U^{(n)}_r(x_\xi,v) \le
\ga_0\mu_AC_0^{dr}$ by the definition of a good block, and the only possible
particle at 
$(x_\xi,v)$ possibly not counted by $N^{(n)}(x_\xi,v)$ is an extra
particle which was added at time 0 at the origin (see \equ(3.21cde)).
Therefore we only need to count jumps of some $\{S_u^{(n)}(x,v,q)\}$ with
$q \le \ga_0 \mu_A C_0^{dr}+1$. 
It follows that the total number of jumps
in $[(k-1)\De_r, (k+1)\De_r)$ in the good block $\cB_r(\bold i,k)$
of particles outside $\cV_r(\bold i)$ at time $(k-1)\De_r$ is stochastically
bounded by
$$
\align
\sum_{x \in \cW_r(\bold i)} \;&\sum_{v \in [(k-1)\De_r,(k+1)\De_r)}
\;\sum_{q \le \ga_0 \mu_A C_0^{dr}+1} 
I\big[S^{(n)}_1(x,v,q) \ne \bold 0,\\
&\phantom{MMMMM}
\sup_{u \le 2\De_rn} \|S_u^{(n)}(x,v,q)-S_1^{(n)}(x,v,q)\| \ge 2\De_r-2\big]\\
&\times [\text{number of jumps of $S^{(n)}_u(x,v,q) ,m \le
u \le m+2\De_rn]$}  
\teq(4.25)
\endalign
$$
($v$ is restricted to the multiples of $1/n$ in the second sum; the
bound here is valid in each $(n)$-system with $n < \infty$). 
Each of the random variables [number of jumps of $S^{(n)}_u(x,v,q) ,m \le
u \le m+2\De_rn]$ converges (as $n \to \infty$) in distribution to 
a Poisson variable, $X(x,v,q)$ say, with mean $2\De_rD$. Furthermore
$$
\align
\sum_{x \in \cW_r(\bold i)} \;&\sum_{v \in [(k-1)\De_r,(k+1)\De_r)}
\;\sum_{q \le \ga_0 \mu_A C_0^{dr}+1}
I\big[S^{(n)}_1(x,v,q) \ne \bold 0,\\
&\phantom{M}
\sup_{u \le 2\De_rn} \|S_u^{(n)}(x,v,q)-S_1^{(n)}(x,v,q)\| \ge 2\De_r-2]
\teq(4.25y)
\endalign
$$
converges (as $n \to \infty$) in distribution to a Poisson random
variable, $T= T(\bold i,k)$ say, of mean 
$$
\align
&K_7\De_r^dC_0^{dr}D \limn P\{\sup_{u \le 2\De_rn} 
\|S_u^{(n)}(x,v,q)-S_1^{(n)}(x,v,q)\| \ge 2\De_r-2]\}\\
&\le K_8 C_0^{7dr}\exp[-K_9\De_r]
\teq(4.25z)
\endalign
$$
for some constants $K_7-K_9$ which depend on $d, D$ and $\ga_0\mu_A$ only.
Moreover, $T$ and all $X(x,v,q)$ are independent. 
Thus (see \equ(4.23ff))
$P\{M_{\text{out}}(r,\bold i,k) > x\} \le P\{\wt M(r,\bold i,k)\ge  x\}$
for $\wt M(r,\bold i,k) =\sum_{j=1}^T X_j $
with Poisson variables $X_j$ with mean $2\De_rD$, independent of each
other and of $T$.

We finally show that \equ(4.13) holds for the $M_{\text{out}}(r, \bold
i,k), (\bold i,k) \in \cS(\bold a,b)$ for fixed $(\bold a,b)$, and with the
$\wt M(r, \bold i,k)$ as above and $\cS(\bold a,b)$ any collection of
pairs $(\bold i,k) \equiv (\bold a,b)$. Firstly, the 
sums in \equ(4.25) for different $(\bold i,k) \equiv (\bold a,b)$ use
different random walks $\{S^{(n)}_u\}$ and therefore are independent. 
From the argument in the last paragraph it then follows that 
the $M_{\text{out}}(r, \bold i,k), (\bold i,k) \in \cS(\bold a,b)$, are
dominated by an independent family of random variables $\wt M(r,\bold
i,k)$, each of which has the distribution of $\sum_{j=1}^T X_j$.
A straighforward calculation gives 
$$
E\{e^{\th \wt M(r,\bold i,k)}\} \le \exp\Big[K_8 C_0^{7dr}
\exp[-K_9\De_r]\big[e^{2\De_rD(e^\th-1)}\big]\Big].
\teq(4.26)
$$
Thus \equ(4.13) holds for any $r \ge 1$ with $\th_r =1$ and $\log
\Ga_r = K_{10}$ for any constant $K_{10} = K_{10}(d,D, \ga_0\mu_A)
\ge \sup_{r \ge 1} K_8C_0^{7dr}
\exp\big[-K_9\De_re^{2\De_rD(e-1)}\big]$.

In order to  apply Lemma 15, we have to have an estimate on 
$$
P\{\text{for some }\wh \pi \in \Xi(\l,t), \; [H_2(\wh \pi,r)]^c 
\text{ occurs}\}.
\teq(4.27)
$$
But this is trivial for $r \le R(t)$, if we take $\ep_r = 3^{d+1}/\De_r$.
Indeed, with this $\ep_r$ and $r \le R(t)$, $H_2(\wh \pi,r)$ never
fails, because the total number of $r$-blocks intersecting a given
$\wh \pi \in \Xi(\l,t)$ is at most $\la_r(\l) \le \ep_r(t+\l)$ (see
\equ(4.19abc) and recall that $\De_r = C_0^{6r} \le C_0^6[\log t]^{6/d}$ by 
\equ(4.19hij), and finally that we can take $t_2$ so that 
$C_0^6[\log t]^{6/d}) \le t$ for $t \ge t_2$). 

The lemma now follows from \equ(4.21) with $x$ equal to
the right hand side of \equ(4.20) with $\th_r, \ep_r$ and $\Ga_r$ as above.
\hfill $\blacksquare$
\enddemo

\comment
already in Lemma 11 of \cite {KSa}. Indeed, any
block $\cB_r(\bold i,k)$ which is contaminated as defined above, is
also contaminated in the sense of \cite {KSa}, 

============================================================

except that 
\cite {KSa} only counts blocks which are contaminated by one of the
particles which were present as $A$-particles at time $0-$, and not by
the added $B$-particle at time 0. Thus we may now have an additional
contribution to   .

. Therefore, 
$$
\sup_{\wh \pi \in \Xi(\l,t)}\sum_{(\bold i,k)}{}^{(\wh \pi,r)} 
I_1(r, \bold i,k) \le \Om^c_r(\l),
$$
with $\Om_r^c(\l)$ as defined before (6.10) in \cite {KSa}.
We now take 
$$
\ep_r = \ka\exp[-C_0^{(r-1)/2}]
\teq(4.28)
$$
 and apply
Lemma 11 in \cite {KSa}. This shows that there exist constants
$C_7=C_7(d), \ka_0,t_0$ (independent of $r, \l$) such that for $\ka \ge
\ka_0, t \ge t_0, 2 \le r \le R(t)$ and $\l \ge 0$ the probability \equ(4.27)
is at most 
$$
P\{\Om_r^c(\l) > \ep_r(t+\l)\} \le\exp[-C_7\ka(t+\l) \exp[-C_0^{(r-1)/2}].
\teq(4.27aa)
$$ 
(We note in passing that the proof of \equ(4.27aa) is essentially
given above. We merely have to leave out the factor 
[number of jumps of $S^{(n)}_u(x,v,q) ,0 \le
u \le 2\De_r]$ in \equ(4.25) to obtain a bound for $\Om_r(\l)$].)
Finally we take $\eta_r, \ep_r$ as in \equ(4.26) and \equ(4.28), $\ka
= \ka_0$ and 
$x = (t+\l)[\De_r]^{-1/2}$. We then see that we can choose $R_2$ and
$t_2$ so large that for $t \ge t_2$ and $R_2 \le r < R(t)$ 
$$
\align
&P\{\text{for some } \wh \pi \in \Xi(\l,t),\;
\sum_{(\bold i,k)}{}^{(\wh \pi,r)} M_{\text{out}}(r, \bold i,k) I_1(r, \bold
 i,k) \ge (t+\l) [\De_r]^{-1/2}\}\\
&\le \exp[-C_7\ka_0(t+\l) \exp[-C_0^{(r-1)/2}]
+ K_{15} [t\log t]^d\exp[-K_{18}(t+\l)\De_r^{-1/2}]\\
&\le \exp[-C_7(t+ \l)/\log t]
\teq(4.30)
\endalign
$$
(recall that $C_0^{R-1} \le [\log t]^{1/d}$).
Finally we can remove the $I_1$ in the sum in the left hand side,
because $I_1(r,\bold i,k)$ already equals 1 whenever 
$M_{\text{out}}(r,\bold i,k) \ne 0$.
\hfill $\blacksquare$
\enddemo

\noindent
{\bf Remark.} We shall later need a small variant of the estimate \equ(4.27aa).
Let 
$$
\wt \mu = 2\ga_0\mu_A.
\teq(4.30ab)
$$
$C_1$ is defined in Theorem 1 of \cite{KSb}. This constant depends
only on $\mu_A,d$ and $D$, since these are the only parameters appearing
in the model (now that the $A$ and $B$-particles have the same
jumprate $D$). Therefore, if $\mu_A$ is replaced by $\wt \mu$, then
Theorem 1 of \cite {KSb} again holds, but now with $C_1$ replaced by
some constant $\wt C_1$. Without loss of generality we take $\wt C_1$
to be a positive integer.

For this constant $\wt C_1$ we can consider 
$$
\align
\wt \Om_r(\l) := [&\text{number of good $r$-block $\cB(\bold i,k)$
 which intersect
$\wh \pi$ and for which}\\
&\text{some particle which is outside 
 $ \prod_{s=1}^d[(i(s)-3 \wt C_1)\De_r, (i(s)+1 + 3\wt C_1)\De_r)$}\\
&\text{at time $(k-1)\De_r$ enters 
 $ \prod_{s=1}^d[(i(s)-2 \wt C_1)\De_r, (i(s)+1 + 2\wt C_1)\De_r)$}\\
&\text{during $[(k-1)\De_r, (k+1)\De_r)$)}.
\endalign
$$
We can then use the proof of \equ(4.27aa) to show that also
$$
P\{\wt \Om_r(\l) \ge \ka(t+\l)\exp[-C_0^{(r-1)/2}]
(t+\l)\} \le\exp[-C_7\ka(t+\l) \exp[-C_0^{(r-1)/2}],
\teq(4.27bb)
$$
for $\ka \ge
\ka_0, t \ge t_0, 2 \le r \le R(t)$ and $\l \ge 0$,
with some changes in the constants $\ka_0, t_0$.
The proof proceeds as in Lemma 11 of \cite {KSa} or via a bound like
\equ(4.25). One has to redefine $\cW_r(\bold i)$ as
$$
\partial \prod_{s=1}^d [(i(s) -3\wt C_1)\De_r,(i(s)+1 + 3\wt C_1)\De_r-1],
\teq(4.27cc)
$$
and then needs to know that for a good block, there are at most  
$\ga_0\mu_AC_0^{dr}$ particles at a space-time point $(x,v)$ with $x
\in \cW_r(\bold a)\times \{(k-1)\De_r, (k+1)\De_r\}$, if $\cB_r(\bold
i,k)$ is good.
In particular we would
have to say that $\cB_r(\bold i,k)$ is good if  
$$
\align
&U_r(x,(k-1)\De_r) \le \ga_r\mu_AC_0^{dr} \text{ for all $x$ for which  }\\
&\phantom{MMMM}
\cQ_r(x) \subset \prod_{s=1}^d [(i(s)-3\wt C_1)\De_r, (i(s)+1+3\wt C_1)\De_r).
\endalign
$$
We then have to reprove the estimate \equ(4.5) for the new definition
of ``good''. One way to do this is to  repeat most of Section 6 of
\cite {KSa} with changes in many constants.
These changes are not difficult, but to spare the reader a great deal
of checking, we point out another possible approach. 
The block 
$\prod_{s=1}^d [(i(s) -3\wt C_1)\De_r,(i(s)+1 + 3\wt C_1)\De_r-1]$,
appearing in \equ(4.27cc) is the union of
translates of $\cB_r(\bold i,k)$ by a fixed finite number of vectors
$\bold z \in \Bbb Z^d$. Therefore the number of  bad
blocks of this form which intersect a space-time path $\wh \pi \in \Xi(\l,t)$
equals the sum of the number of bad blocks in the original
definition which intersect $\wh \pi - \bold z$ over the finite number
of possible $\bold z$. Since $\wh \pi - \bold z$ is again a space-time
path, we again have that it intersects at most $\Phi_r(\l)$
bad blocks of the original type. This does require that we extend the
definition of $\Xi(\l,t)$ to include paths in $\bold z + \cC(t\log t)$
instead of $\cC(t \log t)$. This, however, is much simpler change than
changing all constants, since $\bold z + \cC(t\log t)\subset \cC(2t \log
t)$ for large $t$ and replacing $\cC(t \log t)$ by $\cC(2t\log t)$ is
harmless. 
 \bigskip
\endcomment

Lemma 16 takes care of all contributions to \equ(4.23)
from jumps at some $(x,s)$ in some good $\cB_r(\bold i,k)$ with $r(x,s)= r
\in [R_1,R(t)]$, and
such that the particle which jumps at $(x,s)$ was outside $V_r(\bold
i,k)$ at time $(k-1)\De_r$. Next we consider the jumps 
at some $(x,s)$ in some good $\cB_r(\bold i,k)$ with $r(x,s)= r
\in [R_1+1,R(t)]$, and
such that the particle which jumps at $(x,s)$ was inside $V_r(\bold
i,k)$ at time $(k-1)\De_r$. Note that $r(x,s) = r > R_1$ implies that these
jumps lie in addition in a bad $(r-1)$-block.
We shall therefore estimate
$\sum_{(\bold i,k)}^{(\wh \pi,r)}
 M_{\text{in}}(r, \bold i,k)I_2(\wh \pi,r, \bold i,k)$
where 
$$
I_2(\wh \pi, r,\bold i,k) :=I[\cB_r(\bold i,k) \text{ is contained in a bad
$(r-1)$-block which intersects }\wh \pi].
$$
\proclaim{Lemma 17}
 There exist constants $R_3$ and $t_3$ such that
for $t \ge t_3, R_3 \le r \le R(t)$ and $\l \ge 0$,
$$
\align
&P\Big\{\text{$\Th(t)$ \rm and there exists a }\wh \pi\in \Xi(J,\l,t) 
\text{ \rm such that }\\
&\phantom{MMMMMMMMMM}\sum_{(\bold i, k)}{}^{(\wh \pi,r)} 
M_{\text{in}}(r,\bold i,k) I_2(\wh \pi,r, \bold i,k) \ge 
\frac{8(12)^dK_2 (t+\l)}{ \De_r}\Big\}\\
& \le 2\exp\big[- \sqrt{(t+\l)}\big].
\teq(4.31)
\endalign
$$
\endproclaim
\demo{Proof} Again this proof relies on \cite {KSa}. First we modify
the $M_{\text{in}}$ somewhat, 
so that we can verify \equ(4.13). If $\cB_r(\bold i,k)$ is
good, then $\cV_r(\bold i)$ contains at most
$[7\De_r]^d\ga_0\mu_A C_0^{dr}+1$ particles, so that $M_{\text{in}}(r,
\bold i, k)$ counts the number of jumps in $\cB_r(\bold i,k)$ during
$[k\De_r, (k+1)\De_r)$ of at most $[7\De_r]^d\ga_0\mu_A C_0^{dr}+1$
particles. (Again the one is added to the number of particles to take
into account the extra particle added at time 0.) 
$M_{\text{in}}(r,\bold i, k)$ is therefore bounded by the 
total number of jumps during
$[k\De_r, (k+1)\De_r)$ of the first $[7\De_r]^d\ga_0\mu_A C_0^{dr}+1$ particles
in $\cV_r(\bold i)$ in some arbitrary ordering of particles; if
there are fewer than $[7\De_r]^d\ga_0\mu_A C_0^{dr}+1$ particles in
$\cV_r(\bold i)$ we add artificial particles to raise the number to 
$[7\De_r]^d\ga_0\mu_A C_0^{dr}+1$ and count the jumps of the
artificial particles as well. Denote the resulting number of jumps 
by $\wt M_{\text{in}}(r,\bold i,k)$. By construction, each of the
$\wt M_{\text{in}}(r,\bold i,k)$ is a Poisson variable with mean 
$D\De_r\{[7\De_r]^{d+1}\ga_0\mu_A C_0^{dr}+1\}$. 
Moreover, if $\cB_r(\bold i,k)$
and $\cB_r(\bold {i'},k')$
have disjoint pedestals, then their corresponding $\wt M_{\text{in}}$
are independent since they count jumps of disjoint sets of  particles,
 and the cardinalities of the sets are non random. Thus \equ(4.13)
 holds for $\wt M(r,\bold i,k)$ Poisson variables with mean
$D\De_r\{ [7\De_r]^{d+1}\ga_0\mu_A C_0^{dr}+1\}$, and correspondingly
$$
\th_r = 1, \log \Ga_r = 2D\De_r\{[7\De_r]^{d+1}\ga_0\mu_AC_0^{dr}+1\}(e-1).
\teq(4.51wx)
$$

Next we check \equ(4.12). By definition of $I_2$,
$\sum_{(\bold i,k)}{}^{(\wh \pi,r)} I_2(\wh \pi,r, \bold i,k)$ is bounded by 
$\phi_{r-1}(\wh \pi)$, the
number of bad $(r-1)$-blocks which intersect $\wh \pi$. 
However, $\phi_{r-1}$
is already estimated in Lemma 15 of \cite {KSa}. 
\comment
We should increase the estimate of \cite {KSa} by 1, again to take
into account the possible influence of the extra particle added at
$(\bold 0,0)$, which does not appear in \cite {KSa}.
\endcomment
The proof of lemma 15 there (especially the one but last member 
of (6.47)) tells us that 
for suitable constants $K_i,C_i$, $t \ge$ some $t_3$ and $r-1 \ge d$
$$
\align
&P\{\text{$\Th(t)$ \rm and there exists a }\wh \pi\in \Xi(\l,t) 
\text{ \rm such that }\\
&\phantom{MMMMMMMM}\sum_{(\bold i, k)}{}^{(\wh \pi,r)} 
I_2(\wh \pi,r, \bold i,k) \ge  
K_{14}\ka_0(t+\l) \exp\big[-K_{13}C_0^{(r-1)/4}\big]\\
& \le \sum_{q=r-1}^{R-1}
\exp\big[-C_7\ka_0(t+\l)\exp[-C_0^{q/2}]\big]\\
&\phantom{MMMMMMMM} + \sum_{q=r-1}^{R-1}\exp\Big[-C_{10}\ka_0(t+\l)
\exp\big[-\frac 1{2(d+1)}\ga_0\mu_AC_0^{(d-\frac 34)q}\big]\\
&\le \exp\big[-\sqrt{(t+\l)}\big].
\endalign
$$
Thus, if we take
$$
\ep_r = K_{14}\ka_0\exp\big[-K_{13}C_0^{(r-1)/4}\big]
$$
then
$$
\align
&P\{\text{$\Th(t)$ \rm and there exists a }\wh \pi\in \Xi(\l,t) 
\text{ \rm such that }\sum_{(\bold i, k)}{}^{(\wh \pi,r)} 
I_2(r, \bold i,k) \ge \ep_r(t+\l)\}\\
&\le \exp\big[-\sqrt{(t+\l)}\big].
\endalign
$$
Finally, we apply \equ(4.21) with $x = 8(12)^dK_2(t+\l)/\De_r$
to obtain \equ(4.31) for $R_3 \le r \le R(t),t \ge t_3$,
 with suitable $R_3,t_3$.
\hfill $\blacksquare$
\enddemo

We go back to \equ(4.23).  Each jump at $(x,s)$ on some $J$-path 
is counted in some
$M_{\text{out}}(r,\bold i,k)$ or some $M_{\text{in}}(r,\bold
i,k)$. Lemma 16 takes care of all jumps of the former kind with
$R_2\le r(x,s) \le R(t)$, whereas Lemma 17 takes care of the jumps of the
latter kind, but only if $(R_1+1) \lor R_3 \le r(x,s) \le R(t)$. On $\Th(t)$
there are no jumps with $r > R(t)$ to consider. Without loss of
generality we can take $R_1 \ge R_2 \lor R_3$.
The only contributions to $J(t, \bold 0)$ which we still must estimate
are then counted in
$$
\sum_{(\bold i,k)}{}^{(\wh \pi,R_1)} M_{\text{in}} (R_1,\bold i,k).
\teq(4.33)
$$
This sum will be broken up into several subsums. But we must first
introduce a certain constant $\wt C_1$. Let 
$$
\wt \mu = 2\ga_0\mu_A.
\teq(4.30ab)
$$
$C_1$ is defined in Theorem 1 of \cite{KSb}. This constant depends
only on $\mu_A,d$ and $D$, since these are the only parameters appearing
in the model (now that the $A$ and $B$-particles have the same
jumprate $D$). Therefore, if $\mu_A$ is replaced by $\wt \mu$, then
Theorem 1 of \cite {KSb} again holds, but now with $C_1$ replaced by
some constant $\wt C_1$. Without loss of generality we take $\wt C_1$
to be a positive integer.

We now break the sum \equ(4.33) up into the two sums:
$$
\sum_{(\bold i,k)}{}^{(\wh \pi,R_1)} M_{\text{in}} (R_1,\bold i,k)
I_3(R_1,\bold i,k)
\text{ and }
\sum_{(\bold i,k)}{}^{(\wh \pi,R_1)} M_{\text{in}} (R_1,\bold i,k)
I_4(R_1,\bold i,k),
$$
where
$$
\align
I_3(r, \bold i,k)\;&:=I\Big[\text{$\cB_r(\bold i,k)$ is good and there
exists a J-path from}\\
&\text{some point }(x',s') \in \cB_r(\bold i,k), 
\text{ to a point $(x'',s'') \in$}\\
&\partial\prod_{s=1}^d
[(i(s)-1)\De_r, (i(s)+2)\De_r-1] \times
\big[s', \big(s' + \frac {\De_r}{2\wt C_1}\big) 
\land (k+1)\De_r\big)\Big]
\endalign
$$
and 
$$
I_4(r, \bold i,k)= 1-I_3(r,\bold i,k).
$$
It will turn out that the sum with $I_4$ can easily be reduced to the
sum with $I_3$ (see Lemma 23). However, the latter sum  will have to
be split up further. We define
$$ 
\align
I_{3,1}(r, \bold i,k) = I\Big[&\text{$\cB_r(\bold i,k)$ is a good
$r$-block, but some particle}\\
&\text{which is outside $V_r(\bold i)$ at time
$(k-1)\De_r$ enters }\\
&\text{$\prod_{s=1}^d [(i(s)-1)\De_r, (i(s)+2)\De_r-1]$
during $[(k-1)\De_r, (k+1)\De_r\Big]$};
\endalign
$$
$$ 
\align
I_{3,2}(r, \bold i,k) = I\Big[&\text{there
exists a J-path using only particles}\\
&\text{in $V_r(\bold i)$ at time $(k-1)\De_r$
and running from} \\
& \text{some point 
$(x',s') \in \cB_r(\bold i,k)$ to a point
$(x'',s'')\in$}\\
&\partial\prod_{s=1}^d
[(i(s)-1)\De_r, (i(s)+2)\De_r-1]\times\big[s',\big(s'
+ \frac{\De_r}{2\wt C_1}\big)\land (k+1)\De_r\big)\Big].
\endalign
$$
If $I_{3,1}(r, \bold i,k) = 0$, i.e., if no particles from the 
outside of $\cV_r(\bold i)$ enter 
$\prod_{s=1}^d [(i(s)-1)\De_r, (i(s)+2)\De_r-1]$
during $[(k-1)\De_r, (k+1)\De_r]$, but there is a J-path from
$\cB_r(\bold i,k)$ to $\partial \prod_{s=1}^d 
[(i(s)-1)\De_r, (i(s)+2)\De_r-1] \times [(k-1)\De_r, (k+1)\De_r]$, 
then this J-path cannot use any
particles which come from outside $\cV_r(\bold i)$. Consequently 
$$
I_3(r,\bold i,k) \le I_{3,1}(r,\bold i,k) + I_{3,2}(r,\bold i,k).
\teq(4.55)
$$
\proclaim{Lemma 18}
 There exist constants $R_4$ and $t_4$ such that
for $t \ge t_4, R_4 \le r \le R(t)$ and $\l \ge 0$,
$$
\align
&P\big\{\text{$\Th(t)$ \rm and there exists a }\wh \pi\in \Xi(J,\l,t) 
\text{ \rm such that }\\
&\phantom{MMMMMMM}\sum_{(\bold i, k)}{}^{(\wh \pi,r)} 
M_{\text{in}}(r,\bold i,k) I_{3,1}(r, \bold i,k) \ge
\frac{K_4(t+\l)}{\De_r}\big\}\\
&\le K_5\exp\big[-K_6(t+\l)/[\log t]^6.
\teq(4.55ab)
\endalign
$$
\endproclaim
\demo{Proof}
The sum $\sum_{(\bold i,k)}^{(\wh \pi,r)}I_{3,1}(r, \bold i,k)$ has
already been estimated on the event $\Th(t)$ in
the proof of Lemma 16 (or alternatively in Lemmas 10 and 11 of \cite
{KSa}). We mention a few of the steps, because essentially the same
estimate will be needed again right after Lemma 27 below (see
\equ(4.29z) for $\cW_r$).
Define
$$
\align
I_5(r, \bold i,k)\;&: = I\Big[\cB_r(\bold i,k) \text{ is good, but
there is a particle which is in $\cW_r(\bold i)$}\\
&\text{at some time
$u \in [(k-1)\De_r,(k+1)\De_r)$ and which visits}\\
&\text{$\prod_{s=1}^d [(i(s)-1)\De_r,
(i(s)+2)\De_r-1]$ at some later time in $[u, (k+1)\De_r)$}\Big].
\endalign
$$
Clearly $I_{3,1} (r,\bold i,k) \le I_5(r, \bold i,k)$.
One now uses the construction by means of
the discrete time random walks $\{S^{(n)}_u(x,s,q)\}$ as in the proof of
Lemma 16. The discrete time analogue of 
$\sum_{(\bold i,k)}^{(\wh \pi,r)}I_5(r, \bold i,k)$ 
is then the number of good
$r$-blocks $\cB_r(\bold i,k)$ which intersect $\wh \pi$ and for which there 
exists a particle whose path contains a piece with the properties  
in \equ(4.25x). The discrete time analogue of $I_5(r, \bold i,k)$
itself is stochastically bounded by the triple sum in \equ(4.25y).
As we saw, the weak limit of \equ(4.25y) is a Poisson variable
$T(\bold i,k)$ with mean bounded by \equ(4.25z). 
Moreover, the triple sums in \equ(4.25y) for different $(\bold i,k)
\equiv (\bold a,b)$ are independent, as observed just before \equ(4.26).
The number of summands in $\sum_{(\bold i,k)}^{(\wh \pi,r)}$ is at
most equal to the number of $r$-blocks which intersect $\wh \pi$, and
this is bounded by $\la_r(\l)$ (see \equ(4.19abc)).
Therefore, there exists some constant $K_{15}$, such that 
for large $t, r \le R(t)$ and for fixed $\wh \pi \in \Xi(J,\l,t)$,
$$
\align
&P\{\sum_{(\bold i,k)} {}^{(\wh \pi,r)} I_5(r, \bold i,k) \ge x\}\\
&\le \sum_{(\bold a, b)} P\{\text{a Poisson variable of mean }\la_r(\l)
K_8 C_0^{7dr}\exp[-K_9\De_r] \text{ is at least } \frac x{2(12)^d}\}\\
&\le 2(12)^d P\{\text{a Poisson variable of mean } K_{15}(t+\l)C_0^{(7d-6)r}
\exp[-K_9\De_r]\\
&\phantom{MMMMMMMMMMMMMMMMMMMMMMM} \text{ is at least } \frac x{2(12)^d}\}.
\teq(4.55z)
\endalign
$$
However, if $\mu \le 1$ and $y \ge 2e\mu$, then for $\th = \log \frac
y{2\mu} \ge 1$,
$$
\align
&P\{\text{a Poisson variable of mean }\mu \text{ is at least }y\}\\
&\le \exp[-\th y + \mu(e^\th -1)] \le \exp[-\th y/2] \le \exp[-\frac y2
\log \frac y2].
\teq(4.55zy)
\endalign
$$
We take $x = 4(12)^d [\De_r]^{-2d-3}(t+\l)$ in \equ(4.55z). Then
\equ(4.55zy) shows that for $R_4$ sufficiently large and $R_4 
\le r \le R(t)$ 
$$
\align
&P\big\{\sum_{(\bold i,k)} {}^{(\wh \pi,r)} I_{3,1}(r, \bold i,k) \ge
2[\De_r]^{-2d-3}(t+\l)\big\} \\
&\le P\big\{\sum_{(\bold i,k)}{}^{(\wh \pi,r)} I_5(r, \bold i,k) \ge
2[\De_r]^{-2d-3}(t+\l)\big\} \\
&\le 2(12)^d\exp\big[-[\De_r]^{-2d-3}(t+\l) \log 
\{[\De_r]^{-2d-3}(t+\l)\}\big].
\teq(4.55yx)
\endalign
$$
 
Now, by definition of $\sum_{(\bold i,k)}{}^{(\wh \pi,r)}$,
$$
\{\sum_{(\bold i,k)} {}^{(\wh \pi,r)} I_5(r, \bold i,k) \ge y\}
=\{\sum_{(\bold i,k)\in \cS}  I_5(r, \bold i,k) \ge y\},
$$
where 
$$
\cS= \{(\bold i,k):\cB_r(\bold i,k) \text{ is good and intersects }
\wh \pi|_{[0,t]}\}. 
\teq(4.55y)
$$
Consequently
$$
\align
&\{\text{for some }\wh \pi \in \Xi(J,\l,t), \sum_{(\bold i,k)} {}^
{(\wh \pi,r)} I_5(r, \bold i,k) \ge 2[\De_r]^{-2d-3}(t+\l)\}\\
&\phantom{MMMMM}\subset \bigcup_\cS \{\sum_{(\bold i,k)\in \cS}  
I_5(r, \bold i,k) 
\ge 2[\De_r]^{-2d-3}(t+\l)\}
\endalign
$$
and
$$
\align
&\le P\{\text{for some } \wh \pi \in \Xi(J,\l,t), \sum_{(\bold i,k)} {}^
{(\wh \pi,r)} I_5(r, \bold i,k) \ge 2[\De_r]^{-2d-3}(t+\l)\}\\
&\le \sum_\cS P\{\sum_{(\bold i,k)\in \cS}  I_5(r, \bold i,k) \ge 
2[\De_r]^{-2d-3}(t+\l)\}\\
&\le 2(12)^d
\sum_\cS \exp \big[-[\De_r]^{-2d-3}(t+\l)\log \{[\De_r]^{-2d-3}(t+\l)\}\big],
\teq(4.55x)
\endalign
$$
where the union and sum over $\cS$ runs over all collections $\cS$ of the form
\equ(4.55y) for some $\wh \pi \in \Xi(J,\l,t)$. Note that different
paths $\wh \pi$ may give rise to the same $\cS$, but that each
possible $\cS$ appears only once in the sum on the right hand side of
\equ(4.55x). As in the proof of
\equ(4.19efg), the number of possible collections $\cS$ is
bounded by $2^{\la_r(\l)}[2t\log t+1]^d \exp[K_9\la_r(\l)]$.
Thus for 
$\ep_r = 2[\De_r]^{-2d-3}$ and $t_4$ sufficiently large, we have for $t
\ge t_4, R_4 \le r \le R(t)$
$$
\align
&P\{\text{for some }\wh \pi \in \Xi(J,\l,t), 
\sum_{(\bold i,k)} {}^{(\wh \pi,r)} I_{3,1}(r, \bold i,k) \ge
\ep_r(t+\l)\}\\
&\le P\{\text{for some }\wh \pi \in \Xi(J,\l,t), 
\sum_{(\bold i,k)} {}^{(\wh \pi,r)} I_5(r, \bold i,k) \ge
\ep_r(t+\l)\}\\
&\le 2^{\la_r(\l)+1}[24t\log t+1]^d \exp\big[K_9\la_r(\l)
-[\De_r]^{-2d-3}(t+\l) \log \{[\De_r]^{-2d-3}(t+\l)\}\big]\\
&\le \exp  \big[-\frac 12[\De_r]^{-2d-3}(t+\l) \log \{[\De_r]
^{-2d-3}(t+\l)\}\big].
\teq(4.51uv)
\endalign
$$
This estimates the first term in the right hand side of \equ(4.14).

In addition we have already seen in the proof of Lemma 17
that  \equ(4.13) for the $M_{\text{in}}(r, \bold i,k)$ holds 
with $\wt M(r,\bold i,k)$ a Poisson variable and $\th_r, \log \Ga_r$
given in \equ(4.51wx).
Finally we apply Lemma 15 with $x = 8(12)^dK_2(t+\l)/\De_r$ once more to
obtain the lemma, but possibly with different values for the $K_i$
than in \equ(4.29).
\hfill $\blacksquare$
\enddemo

We now start on some technical preparations for estimating $\sum I_{3,2}$.
 For $q=1,2,\dots$ and $a \in \Bbb R$ we define
$$
[a]_q= a(a-1)\dots(a-q+1).
$$
\comment
We further define
$$
\wt V_r(\bold i) :=\prod_{s=1}^d[(i(s)-1)\De_r,(i(s)+2)\De_r).
$$
Note that $\wt V_r(\bold i,k) \subset V_r(\bold i)$
\endcomment
We also need the following $\si$-fields and random variables.
$$
\align
\cJ_r\big(\bold i,(k-1)\De_r\big) :=\; &\text{$\si$-field generated by 
the $N_A(x,0-),  x \in
\Bbb Z^d$, and all paths}\\
&\text{during $[0,(k-1)\De_r]$, as well as the paths on $[(k-1)\De_r, 
\infty)$}\\
&\text{of all particles outside $V_r(\bold i)$ at time $(k-1)\De_r$},
\teq(4.51zy)
\endalign
$$
$$
\align
\wt \cJ_r\big(\bold i,(k-1)\De_r\big) := \;&\text{$\si$-field generated 
by the locations and
numbers of}\\
&\text{particles in $V_r(\bold i)$ at time $(k-1)\De_r$},
\teq(4.51yx)
\endalign
$$
$$
\align
L(x,u) = L_r(x,\bold i,k-1,u) = [&\text{number of particles at $x$ at
time $(k-1)\De_r+u$
which}\\
&\text{were in $V_r(\bold i)$ at time $(k-1)\De_r$}].
\endalign
$$
\proclaim{Lemma 19} There exists an $R_5$, such that if $r \ge R_5$ 
 and $\De_r/2 \le u
\le 3\De_r$, then for distinct $a_1,\dots, a_\l \in \Bbb Z^d$ and
$q_1, \dots,q_\l \in \{1,2,\dots\}$ it holds on the event 
$\{\cV_r(\bold i,k)$ is good\},
$$
\align
&E\big\{\prod_{i=1}^\l [L(a_i,u)]_{q_i} \big|
\cJ_r\big(\bold i,(k-1)\De_r\big)\big\} \\
&=E\big\{\prod_{i=1}^\l [L(a_i,u)]_{q_i} \big|
\wt \cJ_r\big(\bold i,(k-1)\De_r\big)\big\} \\
&\le [2\ga_0 \mu_A]^{\sum_{i=1}^\l q_i}.
\teq(4.30abc)
\endalign
$$
\endproclaim
\demo{Proof} Once we know the numbers and locations of the particles
in $\cV_r(\bold i,k)$, the $L(x,u)$ are
determined by the increments after $(k-1)\De_r$ in the paths of the
particles in $\cV_r(\bold i,k)$. These increments are independent of
$\cJ_r\big(\bold i,(k-1)\De_r\big)$. The conditioning
on $\cJ_r\big(\bold i,(k-1)\De_r\big)$ only effects the $L(x,u)$
through the determination of which
particles are in $\cV_r(\bold i,k)$, because these are the only
particles to be counted in $L(x,u)$. The equality in \equ(4.30abc) is
immediate from this.

We now first prove \equ(4.30abc) in the special case $\l=1, q_1=1$. For
brevity we write $a$ instead of $a_1$. 
The conditional expectations in \equ(4.30abc) are now at most
$$
\sum_{y \in V_r(\bold i)} N^*(y,(k-1)\De_r) P\{y+ S_u = a\}
+\sup_yP\{y+S_u = a\}.
\teq(4.34)  
$$
The last term has to be included because the extra particle added at
time 0 is not counted in the $N^*$, even though it may be in
$V_r(\bold i)$ at time $(k-1)\De_r$.
We have to show that \equ(4.34) is at most $2\ga_0 \mu_A$.
This part of the proof 
is very similar to the proof of Lemma 8, with $p$
replaced by $r$. In fact it is somewhat easier. We take 
$\cM(\pmb \l)$ as in Lemma 8 (with $p$ replaced by $r$) 
but this time define $\La$ by
$$
\La = \La(\bold i, r)=\{\pmb \l \in \Bbb Z^d: \cM(\pmb \l) \subset
V_r(\bold i)\},
$$
and for each $\pmb \l \in \La$ we take $y_{\pmb \l} \in \cM(\pmb \l)$
such that
$$
P\{y_{\pmb \l} + S_u = a\} = \max_{y \in \cM(\pmb \l)}P\{y + S_u =
a\}.
$$ 

From here on one can follow the proof of Lemma 8, merely reversing
some inequality signs and making use of 
$$
\sum_{y \in \cM(\pmb \l)} N^*(y,(k-1)\De_r) \le \ga_0 \mu_A
C_0^{dr}
$$
for all $\pmb \l \in \La(\bold i,p)$, which holds because $\cV_p(\bold
i,k)$ is good.  Note that the analogue of \equ(3.21) now becomes
$$
\align
&\sum_{\pmb\l \in \La} \sum_{y \in \cM(\pmb\l)}
\ga_0\mu_A P\big\{y_{\pmb\l}+S_u =a \}\\ 
&\le  \sum_{\pmb\l \in \La} \sum_{y \in \cM(\pmb\l)}
\ga_0\mu_A P\big\{y+S_u =a\}\\
&\phantom{MM}+\sum_{\pmb\l \in \La} \sum_{y \in \cM(\pmb\l)}
\ga_0\mu_A \big|P\big\{y_{\pmb\l}+S_u =a\}
-P\big\{y+S_u = a\}\big|.
\teq(3.21x)
\endalign
$$
Clearly the first double sum in the right hand side here is
at most $\ga_0\mu_A$ for all $a$.
Moreover, the last double sum is at
most $K_2 \ga_0\mu_A C_0^r[\log u]^du^{-1/2}$ by (5.26) and (6.37) 
in \cite {KSa}. Also, $\sup_y P\{y+S_u = a\} = O(u^{-1/2})$ by the
local central limit theorem. The desired bound $2\ga_0\mu_A$ for
\equ(4.34) for $r \ge $ some $R_5$ now follows.

We now turn to the general case of \equ(4.30abc). Write $Q$ for
$\sum_{i=1}^\l q_i$ and let 
$$
I_i(\rho) =I[\rho \text{ moves from $\cV_r(\bold i,k)$ to $a_i$ at time $u$}].
$$
Note that $\prod_{i=1}^\l [L(a_i,u)]_{q_i}$ equals the number of
distinct ordered $Q$-tuples of particles, with
$q_i$ of these particles located at $a_i$, at time $u$, and which were in
$V_r(\bold i)$ at time $(k-1)\De_r$. Set $Q_0 = 0$ and for $j \ge 1$
set $Q_j=\sum_{i=1}^j q_i$. Then this number of $Q$-tuples can be written as
$$
\sum{}^{(Q)} \prod_{j=1}^\l \prod_{i= Q_j+1}^{Q_{j+1}} I_j(\rho_i),
\teq(4.30bc)
$$
where $\sum{}^{(Q)}$ denotes the sum over all ordered $Q$-tuples
of {\it distinct} particles $\rho_1, \dots, \rho_Q$ which are 
in $V_r(\bold i)$ at time $(k-1)\De_r$. Let us write 
$y_i$ for the position of $\rho_i$ at time $(k-1)\De_r$.
If we take the conditional
expectation of \equ(4.30bc), given  $\wt \cJ_r\big(\bold i,(k-1)\De_r\big)$,
we find that the middle member of \equ(4.30abc) equals
$$
\align
&E\big\{\sum{}^{(Q)} \prod_{j=1}^\l \prod_{i= Q_j+1}^{Q_{j+1}} I_j(\rho_i)\big|
\wt \cJ_r\big(\bold i,(k-1)\De_r\big)\big\}\\
&=\sum{}^{(Q)}  \prod_{j=1}^\l \prod_{i=Q_j+1}^{Q_{j+1}} 
P\{S_u = a_j- y_i\}\\
&\le \prod_{j=1}^\l \prod_{i=Q_j+1}^{Q_{j+1}} \Big [\sum_{\rho_i \in
V_r(\bold i)} P\{S_u = a_j- y_i\}\Big]\\
&\le \prod_{j=1}^\l \prod_{i=Q_j+1}^{Q_{j+1}} [2\ga_0\mu_A] = [2\ga_0\mu_a]^Q.
\endalign
$$
The first equality here holds because the $\rho_i$ are distinct, and their
paths are therefore independent. The first inequality holds, because all
products which appear in the left hand side also appear in the right
hand side after expanding the right hand side. The second inequality
is true by virtue of the bound $2\ga_0\mu_A$ for \equ(4.34).
\hfill $\blacksquare$
\enddemo
\proclaim{Lemma 20} Without loss of generality we can take $R_5$ so
large that for $r \ge R_5$
$$
P\{I_{3,2}(r,\bold i,k)=1|\cJ_r\big(\bold i,(k-1)\De_r\big)\} 
\le \De_r^{-2(d+1)^2}
\teq(4.35)
$$
on the event 
$$
\{\cV_r(\bold i,k) \text{ \rm is good}\}.
\teq(4.36)
$$
\endproclaim
\demo{Proof} The event $\{I_{3,2}(r, \bold i,k)=1\}$ 
is determined by the location
of the particles in $\cV_r(\bold i,k)$ and by the paths of these
particles during $[(k-1)\De_r, (k+1)\De_r)$. From this one easily sees that
the conditional probability in the left hand side of \equ(4.35) equals
$$
P\{I_{3,2}(r, \bold i,k)=1|\wt \cJ_r\big(\bold i,(k-1)\De_r\big)\}.
\teq(4.37)
$$
To estimate this probability on the event \equ(4.36),
we note that on this event there are at most $7^d\De_r^d+1$
particles in $V_r(\bold i)$ at time $(k-1)\De_r$. The
probability that any given one of these particles, $\rho$ say, 
has two jumps within $1/n$ time units from each other during
$[(k-1)\De_r, (k+1)\De_r)$ is at most
$$
\align
&\sum_{j \ge 1}P\{\text{$j$-th jump of $\rho$ after $(k-1)\De_r$ occurs
before $(k+1)\De_r$ and the}\\
&\phantom{MMMMM}\text{next jump occurs $\le 1/n$ time units later 
$\big|\wt \cJ_r\big(\bold i,(k-1)\De_r\big)\}$}\\
&\le \frac Dn 
\sum_{j \ge 1}P\{\text{$j$-th jump of $\rho$ after $(k-1)\De_r$ occurs
before }\\
&\phantom{MMMMM}(k+1)\De_r\big|
\wt \cJ_r\big(\bold i,(k-1)\De_r)\}\\
&\le \frac {2\De_rD^2}n.
\endalign
$$
Therefore, on the event \equ(4.36),
$$
\align
&P\{\text{some particle in $\cV_r(\bold i,k)$ has two jumps within $1/n$ time
units}\\
&\phantom{MM}\text{of each other during }[(k-1)\De_r, (k+1)\De_r)\big| 
\wt \cJ_r\big(\bold i,(k-1)\De_r\big)\}\\
&\le (7^d \De_r^d+1) \frac{2\De_r D^2}n.
\teq(4.39)
\endalign
$$

For the remainder of this proof we take $n = \De_r^{3(d+1)^2}$.
Assume that the event in the left hand side of \equ(4.39) does not
occur.
Now if $\{I_{3,2}(r, \bold i,k) = 1\}$ occurs, and the J-path in this event
starts at $(x',s')$ and $j/n \le s' < (j+1)/n$, then each particle  
at $x'$ at time $s'$ is also at $x'$ at one or both of the times $j/n,
(j+1)/n$. We can therefore let the J-path begin at $(x',j/n)$ or 
$(x',(j+1)/n)$. Consequently, after raising $R_5$ and hence $n$, if
necessary, to make $\De_r - 1/n \ge \De_r/2$, 
the left hand side of \equ(4.35) is at most
$$
\align
&\frac {(7^d \De_r^d+1)2\De_rD^2}{\De_r^{3(d+2)^2}}\\  
&\quad +\sum_{\De_r-1/n \le j/n < 2\De_r+1/n}\; 
\sum \Sb i(s)\De_r \le x(s) < (i(s)+1)\De_r\\ 1 \le s \le d \endSb
P\big\{\text{there exists
a J-path}\\   
&\qquad\text{from $\big(x,(k-1)\De_r +j/n\big)$ to
$ \partial \prod_{s=1}^d 
[(i(s) - 1)\De_r, (i(s)+2)\De_r-1]$}\\
&\qquad \text{of time duration $\le \De_r/(2\wt
C_1)+ 1/n \le \De_r/\wt C_1$ and which uses}\\
&\qquad\text{only particles which are in $V_r(\bold i)$ at 
time $(k-1)\De_r$}\big|\cJ_r\big(\bold i,(k-1)\De_r\big)\}. 
\teq(4.41)
\endalign
$$

We next prepare for the estimation of the probability in the right
hand side here. Fix $j/n$ and $x$ for the time being. We shall 
condition on the numbers and locations of the particles at time
$(k-1)\De_r+ j/n$ which were in $V_r(\bold i)$ at 
time $(k-1)\De_r$. Recall that the number of such particles at $x$ is
denoted by 
$L(x,j/n)=L_r(x,\bold i,k-1,j/n)$. 
We are going to apply Proposition 4 and Remark 2 after it and 
(the proof of) Theorem 1 in \cite {KSb}. 
To this end we bring in the following process.
First we start the $A$-system by choosing the $N_A(x,0-)$ as i.i.d. mean
$\mu_A$ Poisson variables and add an extra $A$-particle at $(\bold 0,0)$. 
We let this $A$-system run till time $(k-1)\De_r$. We then
continue from time $(k-1)\De_r$ 
with only the $A$-particles in $\cV_r(\bold i,k)$. At time
$(k-1) \De_r + j/n$ we add one further $B$-particle at $x$.
We then let this process with the extra $B$-particle 
continue from  time $(k-1)\De_r +
j/n$, using the same rules as for the $\{Y_t\}$ process, 
that is, $A$-particles turn into $B$-particles when they
coincide with a $B$-particle, but particles cannot recuperate from type
$B$ to type $A$. Let $\nu(x, j/n)$ denote the 
number of $B$-particles outside $x+\cC(\De_r)$ at time $(k-1)\De_r +
j/n + \De_r/\wt C_1$  in the resulting process.
Then
$$
\align
&E\{\nu(x,j/n)|\cJ_r\big(\bold i,(k-1)\De_r\big)\}=
E\{\nu(x,j/n)|\wt \cJ_r\big(\bold i,(k-1)\De_r\big)\}\\
&= E\Big\{E\big\{\nu(x,j/n)\big|L(z,j/n), \; z \in \Bbb Z^d\big\}\Big|
\wt \cJ\big(\bold i,(k-1)\De_r\big)\Big\}.
\teq(4.38)
\endalign
$$
\equ(4.30abc) says that the conditional distribution of the
$L(z,j/n), \; z \in \Bbb Z^d$, given 
$\wt \cJ\big(\bold i,(k-1)\De_r\big)$ satisfies 
condition (2.51) of \cite {KSb} with $\mu_A$ replaced
by 
$$
\wt \mu := 2 \ga_0\mu_A.
$$
We think of the $L(\cdot, j/n)$ as the random initial condition for
the process from time $(k-1)\De_r + j/n$ on of the particles in
$\cV_r(\bold i,k)$ plus the one extra $B$-particle inserted at $x$ at time
$(k-1)\De_r + j/n$.  
Proposition 4 and Remark 2 after it and Theorem 1 
in \cite {KSb} then show that on the
event \equ(4.36) the right hand side of \equ(4.38) is at most equal to the 
$$
\align
&E\{\text{number of $B$-particles outside $\cC(\De_r)$ at time
$\De_r/\wt C_1$}\\
&\qquad \text{in the $\{Y_t\}$-process with the initial number of}\\
&\qquad \text{$A$-particles distributed as i.i.d., mean $\wt \mu$ Poisson}\\
&\qquad \text{variables plus one $B$-particle at $\bold 0\}$}\\
&\le 2e^{-\De_r/\wt C_1} \text{ (see (1.3) in \cite {KSb})}.
\teq(4.42)
\endalign
$$

We now return to the estimation of \equ(4.41). By virtue of Lemma 14, 
the probability in \equ(4.41) is at most
$$ 
\align
&P\{\text{some J-path starting at $\big (x,(k-1)\De_r+j/n\big)$
and using only particles}\\
&\phantom{MMM}\text{from $\cV_r(\bold i,k)$ 
leaves $x+\cC(\De_r)$ during }[0,\De_r/\wt C_1]\big|
\cJ_r\big(\bold i,(k-1)\De_r\big)\}\\
&\le 2E\{\text{number of $B$-particles outside $\cC(\De_r)$ at time
$\De_r/\wt C_1$, in the process}\\
&\phantom{MMM}\text{using only particles from $\cV_r(\bold i,k)$ plus
one $B$-particle at $(x,(k-1)\De_r+j/n)$},\\
&\phantom {MMM}\text{as described above}
\big|\cJ_r\big(\bold i,(k-1)\De_r\big)\}\\
&\le 4e^{-\De_r/\wt C_1} \text{ (by \equ(4.38) and \equ(4.42))}.
\endalign
$$

To conclude we substitute the last estimate into \equ(4.41) to obtain
$$
\align
&P\{I_{3,2}(r, \bold i,k)=1|\cJ_r\big(\bold i,(k-1)\De_r\big)\}\\
&\le \frac {(7^d \De_r^d+1)2\De_rD^2}{\De_r^{3(d+1)^2}} + 
\sum \Sb \De_r-1/n \le j/n \\< 2\De_r+1/n \endSb\; 
\sum \Sb i(s)\De_r \le x(s) < (i(s)+1)\De_r\\
1 \le s \le d \endSb
4e^{-\De_r/\wt C_1}\\
&\le \De_r^{-2(d+1)^2},
\endalign
$$
on the event \equ(4.36), provided $R_5$ is taken large enough.
\hfill $\blacksquare$
\enddemo

\proclaim{Lemma 21} Take $R_1 = \max\{R_2, R_3,R_4, R_5\}$. Then
there exist constants $K_{16}$-$K_{18}$, depending on $d$ only,
 and $t_5$ such that for $t \ge t_5$ and $R_1 \le r \le R(t)$ 
$$
\align
P\big\{\sup_{\wh \pi \in \Xi(J,\l,t)} &\sum_{(\bold i,k)}{}^{(\wh \pi,r)}
I_{3,2} (r, \bold i,k) \ge K_{16}\De_r^{-2d-3}(t+\l)\big\}\\
&\le K_{17}\exp \big[-K_{18}\De_r^{-2d-3}(t+\l)\big].
\teq(4.56)
\endalign
$$
\endproclaim
\demo{Proof} This proof follows the standard outline of Lemmas 10 and
11 in \cite {KSa}. Define
$$
Y(\bold i,k) = I[\text{$\cV_r(\bold i,k)$  is good, but $I_{3,2}(r,\bold
i,k)=1$}].
$$
Also, let $\{Z(\bold i,k)\}$ be a system of independent random
variables with
$$
P\{Z(\bold i,k) = 1\} = 1-P\{Z(\bold i,k) = 0\} = \De_r^{-2(d+1)^2}.
$$
We claim that for fixed $\bold a \in \{0, \dots,11\}^d, b = 0$ or 1
$$
\align
\{Y(\bold i,k): \;&(\bold i,k) \text{ such that $(\bold i,k) \equiv
(\bold a,b)$ and $\cB_r(\bold i,k)$}\\
&\text{intersects }\cC(t\log t)\}
\endalign
$$
lies stochastically below the collection 
$$
\align
\{Z(\bold i,k): \;&(\bold i,k) \text{ such that $(\bold i,k) \equiv
(\bold a,b)$ and $\cB_r(\bold i,k)$}\\
&\text{intersects }\cC(t\log t)\}.
\endalign
$$
This claim follows immediately from \equ(4.35). Indeed, 
the event $\{\cV_r(\bold i,k)$ is good\} lies in $\cJ_r\big(\bold i, 
(k-1)\De_r\big)$. Also 
the events $\{Y(\bold {i'},k')=1\}$ for $k' \le k, 
(\bold {i'},k') \ne (\bold i, k),
(\bold {i'},k') \equiv (\bold a, b)$, belong to $\cJ_r\big(\bold i, 
(k-1)\De_r\big)$. Finally, $P\{Y(\bold i,k) = 1|\cJ_r
\big(\bold i,(k-1)\De_r\big) \le \De_r^{-2(d+1)^2}$, by Lemma 20.
(Note that $Y(\bold i,k)= 0$ on the complement of the event
\equ(4.36).) 
With our claim established, it follows that the left
hand side of \equ(4.56) is at most
$$
\sum_{(\bold a,b)} 
P\big\{\sup_{\wh \pi \in \Xi(\l,t)} \sum_{(\bold i,k)\equiv (\bold
a,b)}{}^{(\wh \pi,r)}
Z( \bold i,k) \ge [2\cdot 12^d]^{-1} K_{16}\De_r^{-2d-3}(t+\l)\big\}.
\teq(4.52)
$$
We shall not give further steps in the proof of \equ(4.56), because
from \equ(4.52) on it is the same as for Lemma 11 in \cite {KSa},
with $\chi_{r+1}$ and $r+1$ there replaced by
$\De_r^{-2(d+1)^2}$ and $r$, respectively (see also the proof of
Theorem 9 in \cite {L}).
\hfill $\blacksquare$
\enddemo

\proclaim{Lemma 22} There exist a constant $t_6$
such that for $t \ge t_6$ and $R(t) \ge r \ge \max\{R_2,R_3,R_4,R_5\}$,
$$
\align
P\{ \text{for some } \wh \pi \in \Xi(J,\l,t),\;
&\sum_{(\bold i,k)}{}^{(\wh \pi,r)} M_{\text{in}}(r, \bold i,k) 
I_{3,2}(r, \bold
 i,k) \ge \frac{K_4(t+\l)}{\De_r}\\
&\le  K_5\exp\big[-K_6(t+\l)/[\log t]^{12d+18}\big].
\teq(4.62)
\endalign
$$
\endproclaim 
\demo{Proof} This is now a familiar application of Lemma 15. We use
that \equ(4.13) for the $M_{\text{in}}(r, \bold i,k)$ holds 
with $\wt M(r,\bold i,k)$ a Poisson variable and $\th_r, \log
\Ga_r$ as in \equ(4.51wx). Further, \equ(4.56) gives us
an estimate for the first term in the right hand side of \equ(4.21),
with $I$ replaced by $I_{3,2}$ and 
$$
\ep_r = K_{16}\De_r^{-2d-3}
$$ 
in the definition of $H_2(\wh \pi, r)$. The lemma follows from Lemma
15 with $x = 8(12)^dK_2$ $(t+\l)/\De_r$.
\hfill $\blacksquare$
\enddemo

The next lemma will deal with $\sum_{(\bold i,k)}^{\wh \pi,r)} 
M_{\text{in}}I_4$, but only for $r=R_1$.
\proclaim{Lemma 23} There exist some constants $C_{13}, R_6, K_{17},
K_{18}$ and $t_7$ (all
independent of $\l$) such that for
$t \ge t_7, R_1 \ge \max\{R_j:2\le j \le 6\}$  and $\l \ge C_{13}t$
it holds
$$
\align
&P\{\Th(t) \text{ and for some } \wh \pi \in \Xi(J,\l,t),\;
\sum_{(\bold i,k)}{}^{(\wh \pi,R_1)} M_{\text{in}}(R_1, \bold i,k) 
I_4(R_1, \bold  i,k) \ge \l/4\}\\
&\le 2K_{17}\exp\big[-K_{18}(t+\l)/[\log t]^{12d+18}\big]
+ K_1 [t\log t]^d\exp[-K_3\l/4].
\teq(4.67)
\endalign
$$
\endproclaim
\demo{Proof} We begin with proving the deterministic inequality
$$
\sum_{(\bold i,k)}{}^{(\wh \pi, r)} I_4(r, \bold i,k) \le 
2\cdot 3^d \wt C_1\sum_{(\bold i,k)}{}^{(\wh \pi, r)} I_3(r, \bold i,k) + 4
\cdot 3^d\wt C_1t/\De_r.
\teq(4.50)
$$
This inequality holds for each $r$ and each $\wh \pi$.
To see this, fix $\wh \pi$
and consider the time
intervals $\chi_j :=\big[k\De_r+ j\De_r/(2\wt C_1), 
k\De_r+ (j+1)\De_r/(2\wt C_1)\big)$.
An $r$-block $\cB_r(\bold i,k)$ can intersect $\wh \pi|_{[0,t]}$ only
during a $\chi_j$ with $0 \le j < 2\wt C_1 t/\De_r$.
Fix such a $j$ and assume that for this $j$, 
$\wh \pi|_{\chi_j}$ intersects exactly $\si_j$ distinct good $r$-blocks. 
 There is then a subcollection of at least $\al_j :=
\lceil 3^{-d} \si_j\rceil$ of these blocks such that no two of them 
are adjacent
on $\cL$. Denote this subcollection of  good
$r$-blocks by $\cB_r({\bold i_1},k), \dots, \cB_r(\bold i_{\al_j},
k)$, where $k$ is such that $\chi_j \subset [k\De_r, (k+1)\De_r)$ (only
$r$-blocks with this value of $k$ can intersect $ \wh \pi$ during $\chi_j$).
Without loss
of generality assume these blocks are ordered in the order in which 
$\wh \pi|_{\chi_j}$ first visits them. Then, for each $u < \al_j$ let 
$(x',s')$ be the earliest point in 
$\cB_r(\bold i_u,k)\cap \wh \pi|_{\chi_j}$, and $(x'',s'')$ the
earliest point in
$\cB_r(\bold i_{u+1},k)\cap \wh \pi|_{\chi_j}$.
By our choice of the blocks $\cB_r(\bold i_\l,k)$ we then have that
$\|\bold i_{u+1}-\bold i_u\| > 1$ and $(x'',s'') \in \cB_r(\bold
i_{u+1},k)$. Hence $x''$ lies outside $\prod_{s=1}^d \big[(i_u(s)-1)\De_r,
(i_u+2)\De_r\big)$, so that the piece of $\wh \pi$ from $(x',s')$
to $(x'',s'')$ is a $J$-path from $(x',s') \in \cB_r(\bold i_u,k)$
to $(x'',s'')$ with $s',s'' \in \chi_j$ and 
$x''$ outside of $\prod_{s=1}^d \big[(i_u(s)-1)\De_r,
(i_u+2)\De_r\big)$ and  $s'' \in \big[s',\big(s'+ \De_r/(2\wt
C_1)\big)$.
Thus there have to be at least $\al_j - 1 \ge
3^{-d}\si_j -1$ good $r$-blocks $\cB_r(\bold i,k)$ which are counted
in $\sum_{(\bold i,k)}{}^{(\wh \pi,r)}I_3(r, \bold i,k)$. 
The blocks so obtained for
one given value of $j$ are distinct by construction. However, we may
obtain the same block $\cB_r(\bold i,k)$ a number of times for different
values of $j$. However, we already saw that this can happen only for
$\chi_j \subset [k\De_r, (k+1)\De_r)$, and hence only for $2\wt C_1$
values of $j$.  Consequently,
$$
\align
&\sum_{(\bold i,k)}{}^{(\wh \pi,r)}[I_3(r, \bold i,k) + I_4(r, \bold
i,k)] \\
&= [\text{number of good $r$-blocks which intersect $\wh \pi$}] \\
&\le  \sum_{0 \le j < 2\wt C_1 t/\De_r}\si_j,
\endalign
$$
and
$$
\align
&\sum_{(\bold i,k)}{}^{(\wh \pi,r)}I_3(r, \bold i,k) \ge 
\frac 1{2\wt C_1}\sum_{0 \le j
< 2\wt C_1 t/\De_r} [\al_j - 1] \ge  
\frac 1{2\wt C_1} \sum_{0 \le j < 2\wt C_1 t/\De_r}
[3^{-d}\si_j -1] \\
&\ge \frac 1{2 \cdot 3^d \wt C_1}\sum_{(\bold i,k)}{}^{(\wh \pi,r)}
[I_3(r, \bold i,k) + I_4(r, \bold i,k)] - 2t/\De_r.
\endalign
$$
\equ(4.50) is now immediate.

To prove \equ(4.67)
we shall apply Lemma 15 once more. As in Lemma 17 we have \equ(4.13) for 
$M_{\text{in}}(R_1, \bold i,k) \le \wt M(R_1,\bold i,k)$ with the
$\wt M(R_1,\bold i,k)$ Poisson variables with $\th_{R_1}, \log
\Ga_{R_1}$ as in \equ(4.51wx) with $r$ replaced by $R_1$.
To get a bound for $P\{\Th(t)$ and for some $\wh \pi \in \Xi(J, \l,t),
\; [H_2(\wh \pi,R_1)]^c \text{ occurs}\}$ with $I(R_1, \bold i,k)$ replaced by
$I_4(R_1, \bold i, k)$, we take 
$$
\ep_{R_1} = 4\cdot 3^d \wt C_1 (K_{16}+1)[\De_{R_1}]^{-2d-3}
+\frac{4\cdot 3^d}{1+C_{13}}[\De_{R_1}]^{-1},
$$
and then $R_6$ and $C_{13}$ so large that for $R_1 \ge R_6$
$$ 
\ep_{R_1} \log \Ga_{R_1} \le \frac{K_2}{\De_{R_1}} \quad \text{and}
\quad \frac{8(12)^dK_2}{\De_{R_1}} \le \frac 18.
\teq(4.66)
$$
and
$$
1+ C_{13} \ge \frac{2 [\De_{R_1}]^{2d+2}}{K_{16}} \text{ and } C_{13} \ge 1.
\teq(4.64)
$$  
Note that these requirements can indeed be satisfied for $\th_{R_1} =
1$ and $\log \Ga_{R_1}$ as in \equ(4.51wx) by taking $R_6$ and
$C_{13}$ large in the the indicated order 
(but with $C_{13}$ dependent on $\De_{R_1}$.
We further restrict ourselves to 
$$
\l \ge C_{13}t.
\teq(4.63)
$$

Now, by \equ(4.55), \equ(4.51uv) and \equ(4.56)
$$
\sup_{\wh \pi \in \Xi(J,\l,t)} \sum_{(\bold i,k)} {}^{(\wh \pi,R_1)}
I_3(R_1, \bold i,k) \le 2 (K_{16}+1)[\De_{R_1}]^{-2d-3}(t+\l)
$$
outside an event of probability at most
$$
\align
& \exp  \big[-\frac 12[\De_r]^{-2d-3}(t+\l) \log \{[\De_r]
^{-2d-3}(t+\l)\}\big]\\
&\phantom{MMMMMMMMMMMM}+ K_{17}\exp \big[-K_{18}\De_r^{-2d-3}(t+\l)\big]\\
&\le 2K_{17}\exp \big[-K_{18}(t+\l)/[\log t]^{12d+18}\big],
\teq(4.65zy)
\endalign
$$
provided $t \ge t_7$ and $r \ge R_6$, with $t_7$ and $R_6$ are large
enough (but with the conditions on $R_6$ and $t_7$ still independent
of $\l$). 
We then have from 
\equ(4.50) with $R_1$ for $r$ and \equ(4.63) that
$$
\align
&\sup_{\wh \pi \in \Xi(J,\l,t)} 
\sum_{(\bold i,k)}{}^{(\wh \pi, R_1)} I_4(R_1, \bold i,k)\\
 &\le 
\Big[4\cdot 3^d \wt C_1 (K_{16}+1)[\De_{R_1}]^{-2d-3} + 
\frac{4\cdot 3^d\wt C_1}{1+C_{13}}[\De_{R_1}]^{-1}\Big](t+\l)\\
&=  \ep_{R_1}(t+\l),
\teq(4.65)
\endalign
$$
outside an event of probability at most \equ(4.65).
We therefore have shown that for suitable $t_7, R_6$ (independent of
$\l$ and for $t \ge t_7, R_1 \ge R_6$ and $\l \ge C_{13}t$
$$
\align
P\{\Th(t)\text{ and for some }\wh \pi \in \Xi(J,\l,t),
\;& \sum_{(\bold i,k)}{}^{(\wh \pi,R_1)} I_4(R_1, \bold i,k) > 
\ep_{R_1}(t+\l)\}\\
&\le 2K_{17}\exp \big[-K_{18}(t+\l)/[\log t]^{12d+18}\big].
\endalign
$$
Finally, we apply \equ(4.21) with $x = \l/4$. The values in
\equ(4.66)-\equ(4.64) have been chosen that this satisfies \equ(4.20)
with $r = R_1$. This results in \equ(4.67).
\hfill $\blacksquare$
\enddemo

\demo{Proof of Proposition 13} The definitions of $M_{\text{out}}$ and
$M_{\text{in}}$,  and the lines just before \equ(4.33)
show that for $\wh \pi \in \Xi(J, \l,t)$, on the event $\Th(t)$,
$$
\align
&j(t,\wh \pi) \le \sum_{r=R_1}^{R(t)}\sum_{(\bold i,k)}{}^{(\wh \pi, r)} 
M_{\text{out}}(r, \bold i,k)I_1(r, \bold i,k) 
+ \sum_{r=R_1+1}^{R(t)}\sum_{(\bold i,k)}{}^{(\wh \pi, r)} 
M_{\text{in}}(r, \bold i,k)I_2(r, \bold i,k)\\
&\phantom{iMM}+\sum_{(\bold i,k)}{}^{(\wh \pi, R_1)}
M_{\text{in}}(r, \bold i,k)
I_3(r, \bold i,k)+ \sum_{(\bold i,k)}{}^{(\wh \pi, R_1)}
M_{\text{in}}(r, \bold i,k)I_4(r, \bold i,k).
\teq(4.57)
\endalign
$$
Now any $\wh \pi \in \Xi(\l,t)$ has $\l$ jumps during $[0,t]$ and
therefore, if $\Xi(J,\l,t)$ is nonempty, then for some $\wh \pi$
 one of the four sums in the right hand side 
here must be at least $\l/4$. Consequently,
$$
\align
P\{\Th(t) &\text{ and }\Xi(J,\l,t) \ne \emptyset\} \\
&\le  P\{\Th(t) \text{ and }\sup_{\wh \pi \in \Xi(J,\l,t)}
\sum_{r=R_1}^{R(t)}
\sum_{(\bold i,k)}{}^{(\wh \pi, r)} M_{\text{out}}(r, \bold i,k)
I_1(r, \bold i,k) \ge \l/4\}\\ 
&+P\{\Th(t) \text{ and }\sup_{\wh \pi \in \Xi(J,\l,t)} 
\sum_{r=R_1+1}^{R(t)}
\sum_{(\bold i,k)}{}^{(\wh \pi, r)} M_{\text{in}}(r, \bold i,k)
I_2(r, \bold i,k)\ge \l/4\}\\
&+P\{\Th(t) \text{ and }\sup_{\wh \pi \in \Xi(J,\l,t)} 
\sum_{(\bold i,k)}{}^{(\wh \pi, R_1)}
M_{\text{in}}(r, \bold i,k) I_3(r, \bold i,k) \ge \l/4\}\\
&+P\{\Th(t) \text{ and }\sup_{\wh \pi \in \Xi(J,\l,t)} 
\sum_{(\bold i,k)}{}^{(\wh \pi, R_1)}
M_{\text{in}}(r, \bold i,k) I_4(r, \bold i,k) \ge \l/4\}.
\teq(4.60)
\endalign
$$
We now restrict ourselves to $\l \ge C_{13}t$
and take $R_1 \ge \max\{R_j: 2\le j \le 6\}$ such that also 
$$
\De_{R_1} \ge \frac {C_0^6}{C_0^6-1}\big[16K_4 \lor 64(12)^dK_2\big].
\teq(4.60aa)
$$
Finally we take $t \ge \max\{t_j:1 \le j \le 7\}$ and
large enough for some further inequalities below. We stress that all
these requirements do not depend on the value of $\l$.

Now, to estimate the first term in the right hand side of \equ(4.60),
assume that 
$$
\sum_{(\bold i, k)}{}^{(\wh \pi,r)} 
M_{\text{out}}(r,\bold i,k) I_1(r, \bold i,k) \ge \frac{K_4 
(t+\l)}{ \De_r}\big\}
$$
for all $\wh \pi \in \Xi(J,\l,t)$. Then also for all such $\wh \pi$
$$
\align
&\sum_{r=R_1}^{R(t)}
\sum_{(\bold i, k)}{}^{(\wh \pi,r)} 
M_{\text{out}}(r,\bold i,k) I_1(r, \bold i,k) \ge \frac{K_4 
(t+\l)}{ \De_r}\big\}\\
&< \sum_{r=R_1}^\infty \frac {K_4(t+\l)}{\De_r} = \frac
{K_4(t+\l)}{\De_{R_1} }\frac {C_0^6}{C_0^6-1} \text{ (since $\De_r =
C_0^{6r}$)} \\
&< \frac\l 4 \text{ (by \equ(4.60aa)}. 
\endalign
$$
It therefore follows from \equ(4.29), and the fact that
$R(t) \sim [d\log C_0]^{-1}\log \log t$ (see \equ(4.19hij)),
that the first term in the right
hand side of \equ(4.60) is at most
$$
\sum_{R_1 \le r \le R(t)} K_5\exp\big[-K_6(t+ \l)/[\log t]^6\big] 
\le K_{19} (\log \log t)\exp\big[-K_6(t+ \l)/[\log t]^6\big] 
$$
for some constant $K_{19}$.

In the same way, but using \equ(4.31) instead
of \equ(4.29), we obtain that the second term in the right hand side of
\equ(4.60) is at most
$$
K_{19}( \log \log t)\sum_{\l \ge C_{13}t} \exp\big[- \sqrt{(t+\l)}\big]
$$

The third term in the right hand side contributes at most
$$ 
K_{19}(\log \log t)\exp\big[-K_6(t+\l)/[\log t]^{12d+18}\big],
$$
by virtue of \equ(4.55ab) and \equ(4.62).
Finally, under condition \equ(4.63), by \equ(4.67), 
the fourth term in the right hand side is at most
$$
2K_{17}\exp\big[-K_{18}(t+\l)/[\log t]^{12d+18}\big]
+ K_1 [t\log t]^d\exp[-K_3\l/4].
$$

 We now substitute these estimates in the right hand side of
 \equ(4.60) 
and sum over $\l \ge C_{13}t$. This yields
$$
P\{\Th(t) \text{ and } \Xi(J,\l,t) \ne \emptyset \text{ for some }
\l \ge C_{13}t\}
\le K_{20} \exp [-K_{21}\sqrt t]
\teq(4.70)
$$
for suitable constants $K_{20}, K_{21}$ (which depend on $\De_{R_1}$,
but that is fixed now) and all large $t$. We add $P\{[\Th(t)]^c\}$ (see
\equ(4.13bb))
to obtain that $P\{\Xi(J,\l,t) \ne \emptyset\} \le 2t^{-2}$ for
large $t$. Hence, by Borel-Cantelli,
$\cap_{\l \ge C_{13}t} \Xi(J, \l, t) = \emptyset$ for all 
large integers $t$ a.s.
In view of \equ(4.17aa) and the lines
following it, this implies that $J(t, \bold 0) \le C_{13}t$ for all
large integers $t$ a.s. Since $J(t, \bold 0)$ is nondecreasing in $t$ this
implies Proposition 13 with $C_{12} = 2C_{13}$.
\hfill $\blacksquare$
\enddemo
\medskip \noindent
{\bf Remark 6.} Note that \equ (4.70) proves the explicit estimate
$$
P\{\Th(t) \text{ and }J(t,x) \ge C_{12}t\} \le K_{20}\exp[-K_{21}\sqrt t]
\teq(4.69)
$$
for each fixed $x \in \Bbb Z^d$, for all large $t$.

\subhead
5. Extinction for large $\la$
\endsubhead
\numsec=5
\numfor=1
In this section we show that $\la_c < \infty$. We shall use the
$r$-blocks $\cB_r(\bold i,k)$ and their pedestals $\cV_r(\bold i,k) =
V(\bold i,k) \times \{(k-1)\De_r\}$ as defined in Section 4. $\wt C_1$
is defined just after \equ(4.33). We shall work with the
$\{Y_t(\la)\}$-process in this section. This process 
starts with independent mean
$\mu_A$ Poisson variables $N_A(x,0-)$ for the number of $A$-particles
``just before time 0'' and one additional $B$-particle at $\bold 0$ at time 0,
as explained in the abstract. The $B$-particles turn back into
$A$-particles at rate $\la$, independent of everything else; $\la$ is
called the recuperation rate. A particle $\rho'$ which recuperated at
time $s'$ turns into a $B$-particle again at time $s'' := \inf \{s >
s': \rho'$ jumps onto another $B$-particle $\rho''$ or vice versa at
time $s$\}.

If there is a $B$-particle at the
space-time point $(x,t)$, then there is a genealogical path which
starts at $(\bold 0,0)$ and reaches $(x,t)$. In particular, 
this means that, for some $\l$,
there exist times $s_0 =0 < s_1 < \cdots  < s_\l < s_{\l+1} = t$ 
and particles $\rho_0,\rho_1, \cdots, \rho_\l$ such that $\rho_0$ is
a $B$-particle at $(\bold 0,0)$, $\rho_\l$ is the given  $B$-particle
at $(x,t)$, and $\rho_i$ jumps onto the position of $\rho_{i-1}$
or vice versa at time $s_i, \, 1 \le i \le \l$; moreover, $\rho_i$
is of type $B$ during $[s_i, s_{i+1}],\, 0 \le i \le \l$.
 (See the construction in
the proof of Proposition 5 in \cite {KSa} as well as the comments in
the paragraph following \equ(2.0bb) above.)  
Note that it is not necessary that all particles
$\rho_i, \, 0 \le i \le \l$, are distinct; it is possible that $\rho_i
= \rho_j$ if $|i-j| > 1$. This is due to the possibility of
recuperation, and cannot be ruled out, as was done in the
case without recuperation studied in \cite {KSb}. 
We shall extend our definition of J-path somewhat, so that
a genealogical path such as just discussed is also a J-path. In
section 4 we considered only $A$-particles. But the paths of the
particles are not influenced at all by the types under our basic
assumption that the $A$ and $B$-particles perform the same random
walk. We can therefore define a J-path to be any path which coincides
at all times with some particle, irrespective of type. Otherwise these
paths are exactly as discussed in the beginning of Section 4. All
arguments of the preceding section, and in particular, its
principal result, Proposition 13, remain valid. To see this, one
simply has to ignore the types of all particles. 
A genealogical 
path for a particle at time $t$ coincides at each time in $[0,t]$ 
with some $B$-particle, and is 
therefore a J-path on $[0,t]$. Moreover, in our model,
it has to start at
$\bold 0$, because that is the only site with $B$-particles at time 0.

In this section we want to prove the following result:
\proclaim{Proposition 24} For sufficiently large $\la$
there a.s. exists a (random) time $\tau < \infty$ such that there are
no $B$-particles in $\{Y_t(\la)\}$ after $\tau$.
\endproclaim

The idea of the proof is as follows. Assume that there is a
$B$-particle at $(x,t)$ and let $\wh \pi: [0,t] \to \Bbb Z^d\times
[0,t]$ be its 
genealogical path. For a fixed large $r$ we consider all 
$r$-blocks which intersect $\wh \pi$. Of course there have to be at
least $\lfloor t/\De_r \rfloor$ such blocks, since each $r$-block only
extends over an interval of length $\De_r$ in the time direction. 
The next lemma is the principal
one. It states that for each of these $r$-blocks at least one of
four events $G(j)$ has to occur. We shall then show in a series of lemmas
that there are (with high probability for large $t$) for each $j$  
at most $t/(10\De_r)$ $r$-blocks which intersect $\wh \pi$ and have $G(j)$
occurring. Actually the next lemma leaves one exceptional case. At the
end of the section we show that with high probability this exceptional case
contributes at most a bounded number of $r$-blocks which intersect
$\wh \pi$. In total that gives at most $5t/(10\De_r) < \lfloor t/\De_r \rfloor$
$r$-blocks which intersect $\wh \pi$. This contradiction shows that
with high probability there are for large $t$ 
no points $(x,t)$ with a $B$-particle.
\medskip
For integral $z \ge -1$ we shall use the abbreviation
$$
v_r(k,z) := [k+z/(4\wt C_1)]\De_r.
$$

\proclaim{Lemma 25} Let there be a $B$-particle at $(x,t)$ and let
$\wh \pi:[0,t] \to \Bbb Z^d \times [0,t]$ be a genealogical path from
$(\bold 0,0)$ to $(x,t)$ and let $\cB_r(\bold i,k)$ be an $r$-block
which intersects $\wh \pi$ in a point $(y'',s'')$
with $k \ge 1$ and $v_r(k,z) \le s'' < v_r(k,z+1)\le t$ for some 
integer $z \in[0,4\wt C_1)$.
Then one of the following four events
must occur:
$$
G(1) = G_r(1, \bold i,k) :=\{\cB_r(\bold i,k) \text{ is bad}\};
$$
$$
\align
&G(2,z) = G_r(2,z, \bold i,k) := \{\cB_r(\bold i,k) \text{ is good,
and in the system}\\
&\quad\text{which continues from time $v_r(k,z-1)$ with 
the particles}\\
&\quad\text{in $V_r(\bold i)$ only, there are still some 
$B$-particles at time $v_r(k,z)$}\};
\endalign
$$
$$
\align
&G(3,z) = G_r(3,z, \bold i,k) := \{\cB_r(\bold i,k) 
\text{ is good, but
there is a particle which is}\\
&\quad \text{outside $V_r(\bold i)$ at some time
$u \in [v_r(k,z-1),(k+1)\De_r)$ and which visits}\\
&\quad\text{$\prod_{s=1}^d [(i(s)-1)\De_r,
(i(s)+2)\De_r-1]$ at some later time in $[u, (k+1)\De_r)$}\};
\endalign
$$
$$
\align
&G(4,z) = G_r(4,z,\bold i,k) := \{\cB_r(\bold i,k) 
\text{ is good, but
there is a J-path from}\\
&\quad\text{some $(y',s')$ to $(y'',s'')$ with $y' \in \partial 
\prod_{s=1}^d [(i(s)-1)\De_r,(i(s)+2)\De_r-1]$, }\\
&\quad y'' \in \prod_{s=1}^d[i(s)\De_r, (i(s)+1)\De_r)
\text{ and $v_r(k,z-1)\le s' \le s'' \le v_r(k,z+1)$,}\\
&\quad \text{and this J-path uses only
particles which were in $V_r(\bold i)$ at time $v_r(k,z-1)$}\}.
\endalign
$$
\endproclaim
\demo{Proof}
If $\cB_r(\bold i,k)$ intersects $\wh \pi$, then they
must have some point $(y'',s'')$ with
$y'' \in
\prod_{s=1}^d [(i(s)\De_r, (i(s)+1)\De_r)$ and $s'' \in [k\De_r, (k+1)\De_r)$
in common (by the definition of $\cB_r(\bold i,k$)). Clearly there must then
exist an integer $z \in [0, 4\wt C_1)$ such that
$v_r(k,z) \le s'' < v_r(k,z+1)$. We fix such a $(y'',s'')$ for the remainder
of this proof. Since $(y'',s'')$ is on
the genealogical path for $(x,t)$ there must be a $B$-particle present
at $(y'',s'')$, as we already pointed out. Let $\rho^*$ be such a $B$-particle.

If $G(1)$ fails, then $\cB_r(\bold i,k)$ 
is good, so that this
may be assumed to be the case in $G(2,z)$-$G(4,z)$. Now assume that none of
$G(1), G(2,z)$ or $G(3,z)$ occur.
Since $G(3,z)$ fails, any particle at
$(y'',s'')$ is one of the particles which was in $V_r(\bold i)$ at
time $v_r(k,z-1)$. In particular, this must be true for $\rho^*$.
Consider the genealogical path for $\rho^*$. 
Let $\wh \pi_0$ be the piece 
of this last genealogical path over the time interval $[v_r(k,z-1),s'']$.
This is a genealogical path for $\rho^*$ in a system which starts with all the
particles at time $v_r(k,z-1)$. Assume $\wh \pi_0$ arises
from particles $\rho_i$, such that $\rho_i$ jumps to the
position of $\rho_{i-1}$ or vice versa at time $s_i,\, 1 \le i \le \l$,
and that $\rho _i$ has type $B$ during  $[s_i, s_{i+1}],\, 0 \le i
\le \l$, with $s_0 =v_r(k,z-1), s_{\l+1} = s''$, and $\rho_\l = \rho^*$, 
as explained in the second paragraph of this section.
Let $s_{i_0} \le v_r(k,z) < s_{i_0+1}$. We claim that one of the 
$\rho_i$ with $i \le i_0$ must have been outside
$V_r(\bold i)$ at time $v_r(k,z-1)$. Indeed, if this is not the
case, then the system starting with only the particles in $V_r(\bold i)$
at time $v_r(k,z-1)$ has at least one $B$-particle at the time $v_r(k,z)$. 
To see this, observe that if $\rho_0. \dots, \rho_{i_0}$ all came
came from $V_r(\bold i)$ at time 
$v_r(k,z-1)$, then, by induction on $i$, each of these $\rho_i$ would be
a particle of type $B$ during $[s_i, s_{i+1}]$ in the system of
particles which were in $V_r(\bold i)$ at time $v_r(k,z-1)$. 
In particular $\rho_{i_0}$ would be of type $B$ at time $v_r(k,z)$
in this sytem.
This would contradict the assumption that $G(2,z)$ does not occur.
Our claim follows. 
 
In particular, there is a maximal index $m\le \l$ for which 
$\rho_m$ was outside $V_r(\bold i)$ at some time in $[v_r(k,z-1),s'']$. 
In fact
this maximal $m$ is less than $\l$, since $\rho_\l= \rho^*$ is a particle in 
$\prod_{s=1}^d [(i(s)\De_r, (i(s)+1)\De_r)$, and $G(3,z)$ fails.  
Since $\rho_m$ is outside $V_r(\bold i)$ at some time $u \in
[v_r(k,z-1),s'']$, it does
not enter $\prod_{s=1}^d [(i(s)-1)\De_r,
(i(s)+2)\De_r-1]$ during $[u, (k+1)\De_r)$ (because $G(3,z)$ fails).
This means that at time $s_{m+1}$,\, $\rho_m$ and  hence also
$\rho_{m+1}$, are outside $\prod_{s=1}^d [(i(s)-1)\De_r,
(i(s)+2)\De_r-1]$. The path $\wh \pi_0$ therefore must intersect
$\partial \prod_{s=1}^d [(i(s)-1)\De_r,(i(s)+2)\De_r-1]$ sometime 
during $[s_{m+1}, s''] \subset [v_r(k,z-1), v_r(k,z+1)]$, 
because its endpoint at time $s''$ lies in
$\cB_r(\bold i,k)$. Let the latest intersection of $\wh
\pi_0$ with $\partial \prod_{s=1}^d
[(i(s)-1)\De_r,(i(s)+2)\De_r-1]$ occur at time $s' \in [s_{m+1},
s'')$ and position $y'$. 
Then the piece of $\wh \pi_0$ over the time interval $[s',s'']$
is a J-path which uses at most the 
particles $\rho_{m+1}, \cdots, \rho_\l$,
all of which were in $V_r(\bold i)$ at time $v_r(k,z)$
(by our choice of $m$). Thus $G(4,z)$ occurs with
this J-path, while $(y',s')$ and $(y'',s'')$ have all the required properties.
\hfill $\blacksquare$
\enddemo
We now start on showing that each $G(j)$ occurs on relatively few
blocks. We remind the reader of the definitions of $\Phi_r(\l)$ and
$\Xi(J,\l,t)$ in \equ(3.75), \equ(3.75b), \equ(4.13bcb) (see also (6.1), (6.8)
and (6.9) in \cite {KSa}; note that ``good'' is now defined as in \cite {KSa} 
and not as in \cite {KSb}). In \equ(4.17aa) and the lines following it
we showed that
$$
P\{\text{there exists a J-path starting at $\bold 0$ which is not in
$\cup_{\l \ge 0}\Xi(J,\l,t)$}\} \le 2e^{-t} 
$$
for large $t$. In addition, Proposition 13 (or rather \equ(4.70))
says that for suitable constants $C_{13}, K_{20}, K_{21}$ 
$$
P\{\Th(t) \text{ and }\cup_{\l \ge C_{13}t} \Xi(J,\l, t) \ne \emptyset\}
\le K_{20}\exp[-K_{21}\sqrt t].
$$
Finally, \equ(4.13bb) says that $P\{[\Th(t)]^c\} \le t^{-2}$ for large $t$.
It follows from these that
$$
\align
&P\{\text{there is a $B$-particle at time $t$}\}\\
&\le P\{\text{there exists some
genealogical path $\wh \pi$ leading to a $B$-particle at time $t$}\}\\
&\le 2e^{-t} + t^{-2}+ K_{20}\exp[-K_{21}\sqrt t]\\
&\qquad +P\{\text{$\Th(t)$ and there exists a path in 
$\cup_{\l < C_{13}t} \Xi(J,\l,t)$ which starts at }\bold 0\}\\
&\le 2e^{-t} + t^{-2}+ K_{20}\exp[-K_{21}\sqrt t] \\
&\qquad + \sum_{j=1}^4
P\{\Th(t)\text{ and there exists a path in } \cup_{\l < C_{13}t} \Xi(J,\l,t)
\text{ which starts}\\
&\phantom{MMMMMM} \text{at $\bold 0$ and intersects more than $t/(10 \De_r)$
$r$-blocks $\cB_r(\bold i,k)$}\\
&\phantom{MMMMMM} \text{with $k \ge 1$ for which $G(j)$ occurs}\}\\
&\qquad+P\{\Th(t)\text{ and there exists a path in } \cup_{\l <
C_{13}t} \Xi(J,\l,t)
\text{ which starts at $\bold 0$}\\
&\phantom{MMMMMM} \text{and which intersects more than $t/(10 \De_r)$
$r$-blocks $\cB_r(\bold i,0)$}\}.
\endalign
$$
It therefore suffices for Proposition 24 to prove for
some fixed $r$ and $1 \le j \le 4$
$$
\align
P\{&\Th(t)\text{ and there exists a path in } \cup_{\l < C_{13}t} \Xi(J,\l,t)
\text{ which}\\
&\text{starts at $\bold 0$
and intersects more than $t/(10 \De_r)$
$r$-blocks}\\
&\text{$\cB_r(\bold i,k)$ with $k \ge 1$ for which $G(j)$ occurs}\} \to 0,
\teq(5.1)
\endalign
$$
as well as 
$$
\align
&P\{\Th(t)\text{ and there exists a path in } \cup_{\l <
C_{13}t} \Xi(J,\l,t)
\text{ which starts at $\bold 0$}\\
&\phantom{M} \text{and which intersects more than $t/(10 \De_r)$
$r$-blocks $\cB_r(\bold i,0)$}\} \to 0,
\teq(5.1z)
\endalign
$$
as $t \to \infty$.

For $G(1)$ \equ(5.1) is contained in Lemma 15 of \cite {KSa}.
Indeed, Lemma 15 in \cite {KSa} proves that for suitable constants
$K_{13}, K_{14}, \ka_0$
$$
\align
&P\{\text{there exists some path in $\Xi(\l,t)$ which
intersects more than}\\
&\phantom{M}\text{$K_{14}\ka_0 (t+\l) \exp[-K_{13}C_0^{r/4}]$ 
bad $r$-blocks for some $r \ge d, \l \ge 0$}\} \\
&\le \frac 2{t^2}
\teq(5.2)
\endalign
$$
for all large $t$.
We merely have to take $r_1 \ge d$ so large that 
$$
K_{14}\ka_0 (1+C_{13})
\exp[-K_{13}C_0^{r_1/4}] \le \frac 1{10 \De_{r_1}} 
$$
to obtain for any $r \ge r_1$
$$
\align
P\{&\text{there exists some path in $\cup_{\l < C_{13}t} \Xi(\l,t)$ 
which intersects}\\
&\text{more than $t/(10 \De_r)$ bad blocks}\} \le \frac 2{t^2},
\teq(5.3)
\endalign
$$
which gives \equ(5.1) for $j=1$.

\comment
In order to obtain (5.1) for $j=1$ we also need to show that
$$
\align
P\{&\text{there exists some path in $\cup_{\l < C_{13}t} \Xi(\l,t)$ 
which intersects more than}\\
&\text{$t/(20 \De_r)$ blocks $\cB_r(\bold i,k)$ for which $k \ge 1$ and  
$\cB_r(\bold i,k-1)$ is bad}\} \le \frac 2{t^2}.
\teq(5.3z)
\endalign
$$
But this follows by applying \equ(5.3) after shifting our path by
$\De_r$ in the time direction. More precisely, if $\wh \pi(s) =
(\overline \pi(s),s)$ define a new space-time path $\pi^*$ on
$[0,t]$ by
$$
\pi^*(s) = \cases (\overline \pi(s+ \De_r),s) &\text{ if } 
0 \le s \le t-\De_r\\ (\overline \pi(t),s) &\text{ if } t-\De_r \le s
\le t.
\endcases
$$
If $\wh \pi \in \Xi(\l,t)$, then $\pi^* \in \cup_{m \le \l} \Xi(m,t)$
and if $\cB_r(\bold i,k)$ intersects $\wh \pi$, then $\cB_r(\bold
i,k-1)$ intersects $\pi^*$. Thus if the event in braces in \equ(5.3z)
occurs for a path $\wh \pi$, then the event in braces in \equ(5.3)
occurs for $\pi^*$. 
\endcomment

The next two lemmas will imply (5.1) for $j=2$. It is the only place where the
recuperation rate $\la$ plays a role. For simplicity we formulate this lemma
only in the form in which we use it, even though there is a more
general version. We generalize the definition \equ(4.51zy) to
$$
\align
\cJ_r(\bold i,u) :=\; &\text{$\si$-field generated by 
the $N_A(x,0-),  x \in
\Bbb Z^d$, and}\\
&\text{all paths during $[0,u]$, as well as the paths}\\
&\text{on $[u,\infty)$ of all particles outside $V_r(\bold i)$ at time $u$},
\teq(4.51zzy)
\endalign
$$
Similarly we generalize the definition of a good
pedestal. Specifically, we say that
$V_r(\bold i) \times \{u\}$ is good, if
$$
U_r(x,u) \le \ga_r\mu_AC_0^{dr} \text{ for all $x$ for which }  
\cQ_r(x) \subset V_r(\bold i).
$$
The definition of a good $r$-block then shows that if  
$\cB_r(\bold i,k)$ is good, then so is $V_r(\bold i)\times \{u\}$
for any $u \in [(k-1)\De_r, (k+1)\De_r)$.

\proclaim{Lemma 26} For each $r \ge 1, T \ge 0$ and $0 < \ep
\le 1$ there exists a $\la_1(r,T,\ep)$ such that for all $\la \ge
\la_1$, all $(\bold i,k)$ and all $u \in [(k-1)\De_r, (k+1)\De_r)$
$$
\align
&P\big\{\text{$V_r(\bold i) \times\{u\}$ is good and 
the system which consists at time}\\ 
&\phantom{Mi}\text{$u \in [(k-1)\De_r, (k+1)\De_r)$ of the
particles in $V_r(\bold i)$ only, and}\\
&\phantom{Mi}\text{which develops from time $u$ on 
according to the rules for}\\
&\phantom{Mi}\text{$\{Y_t(\la)\}$, has some $B$-particles at time
$u+T\,\big|\cJ_r(u, \bold i)\big\}$}\\
&\le \ep.
\teq(5.5)
\endalign
$$
\endproclaim
\demo{Proof} Let $\rho_1, \rho_2, \dots, \rho_f$ be all the particles
in $V_r(\bold i)$ at time $u$. If $V_r(\bold i)\times \{u\}$ is good, then
there are at most $\nu = \nu_r := [7\De_r]^d\ga_0\mu_AC_0^{dr}+1$ particles in
$V_r(\bold i)$ at time $u$, so
that $f \le \nu$ (see the beginning of the proof of Lemma 17 
for an explanation of the extra term 1 here).
Let $N$ be a large integer and set $u_m = u+ mT/N,\, 0
\le k \le N$. It suffices to show that with probability at least
$1-\ep$ all particles $\rho_1, \dots, \rho_f$ have type $A$ at some
$u_m, \, 0 \le m \le N$. Indeed if this happens at time $u_m$, then
the particles $\rho_1, \dots, \rho_f$ will all have type $A$ at all
times after $u_m$ (since we are ignoring interactions with all other
particles in the system of this lemma). But whatever types and
locations $\rho_1,
\dots \rho_f$ have at time $u_m$, there is a conditional probability
of at least $\exp[-fDT/N][1-e^{-\la T/N}]^f$ that none of the particles
$\rho_i, 1 \le i \le f$, has a jump during $[u_m,u_{m+1}]$, but that
all of them have a recuperation event during $[u_m,u_{m+1}]$. If this
happens, then all $\rho_i$ will be of type $A$ at time
$u_{m+1}$. (Note that here we use our rule that a jump is needed
before a recuperated particle can become reinfected.) 
It follows from this that the left hand side of \equ(5.5) is at most
$$
\Big[ 1- \exp[-fDT/N][1-e^{-\la T/N}]^f\Big]^N.
$$
Now take $N_0(\nu_r,T,\ep)$ such that
$$
\Big[ 1- \frac 12\exp[-fDT/N_0]\Big]^{N_0} \le \ep \text{ for all }f \le \nu_r,
$$
and then $\la_1 = \la_1(r,T,\ep)$ such that $[1-e^{-\la_1
T/N_0}]^{\nu_r} \ge 1/2$.
\equ(5.5) holds for this value of $\la_1$.
\hfill $\blacksquare$
\enddemo
\proclaim{Lemma 27} For each $r \ge 1$ there exists a $\la_0(r)$ such
that for $\la \ge \la_0$, \equ(5.1) with $G(j)$ replaced by
$\cup_{0 \le z < 4\wt C_1} G_r(2,z,\bold i,k)$ holds.
\endproclaim
\demo{Proof} Fix $r \ge 1$ and $0 \le z < 4\wt C_1$. Define
$Y(z,\bold i,k) = I[G_r(2,z,i,k)]$ and let $\wt Y(z, \bold i,k)$
be the indicator function of
$$
\align
\big\{&\text{$V_r(\bold i) \times\{v_r(k,z-1)\}$ is good and 
the system which consists}\\
&\phantom{Mi}\text{only of the particles in $V_r(\bold i)$ at time 
$v_r(k,z-1)$ and which}\\
&\phantom{Mi}\text{develops from time $v_r(k,z-1)$ on 
according to the rules}\\
&\phantom{Mi}\text{for $\{Y_t(\la)\}$ has some $B$-particles at time
$v_r(k,z)$}\big\}
\endalign
$$
It is immediate from the definitions that $Y(z,\bold i,k) \le \wt
Y(z,\bold i,k)$. Moreover,
by applying the preceding lemma with $u =
v_r(k-1,z)$ and $T = \De_r/(4\wt C_1)$, 
we see that we can find a $\la_0$ such that for all $\la \ge \la_0, \;
0\le z < 4\wt C_1$, and $0 < \ep \le (4\wt C_1)^{-1}$
$$
P\{\wt Y(z,\bold i,k) = 1|\cJ_r(v_r(k,z-1),\bold i)\} \le \ep.
\teq(5.5z)
$$
A fortiori, the same inequality holds if $\wt Y$ is replaced by $Y$.
In fact we have more. Let $\bold a \in \{0, 1, \cdots,11\}^d$ and $b =
0$ or 1. Let further $Z(z,\bold i,k)$ be a family of independent
random variables with
$$
P\{Z(z,\bold i,k) = 1\} = 1-P\{Z(z,\bold i,k)=0\} = \ep.
$$
\equ(5.5z) shows that the conditional probability of $\{\wt Y(z,\bold
i,k)=1\}$, given all the $\wt Y(z,\bold j, \l)$ with $\bold j \equiv \bold
i \mod(\bold a), \l \equiv k \mod (b)$ and $(\bold j,\l)$ preceding
$(\bold i,k)$ in the lexicographic order, is at most $\ep$. 
Just as in the proof of Lemma 21, this shows that for fixed $z$,
the family $\{\wt Y(z,\bold i,k): (\bold i,k) \equiv (\bold a,
b)\}$ lies stochastically below the family $\{Z(\bold i,k): 
(\bold i,k) \equiv (\bold a,b)\}$. Again this statement remains true
if $\wt Y$ is replaced by then smaller $Y$.

We can now continue exactly as in Lemma 11 of \cite {KSa}
or Lemma 21 in Section 4. 
Note that if there exists a $\wh \pi \in \Xi(J,\l,t)$ which
intersects more than $K_1\ep^{1/(d+1)}(t+\l)/\De_r$ 
blocks $\cB_r(\bold i,k)$ for which 
$\cup_{0 \le z < 4\wt C_1}G_r(2,z, \bold i,k)$ occurs, then there
exists a $0 \le z < 4 \wt C_1$ such that $\wh \pi$ intersects 
more than $K_1\ep^{1/(d+1)}(t+\l)/[4\wt C_1\De_r]$ 
blocks $\cB_r(\bold i,k)$ for which 
$G_r(2,z, \bold i,k)$ occurs.
We therefore have for $t \ge 1$ and
for some constants $K_1$-$K_3$, which do not depend on $\ep,\l$ or $r$,
$$
\align
&P\{\text{there exists an $\l < C_{13}t$ and 
a path $\wh \pi \in \Xi(J,\l,t)$
such that  $\wh \pi$ intersects }\\
&\phantom{iiMMMMMMM}\text{more than $K_1\ep^{1/(d+1)}(t+\l)/\De_r$ 
blocks $\cB_r(\bold i,k)$ for}\\
&\phantom{iiMMMMMMM}\text{which 
$\cup_{0 \le z < 4\wt C_1}G_r(2,z, \bold i,k)$ occurs}\}\\
&\le \sum_{0 \le z < \wt 4C_1} P\{\text{there exists an $\l < C_{13}t$ and 
a path $\wh \pi \in \Xi(J,\l,t)$
such that  $\wh \pi$}\\
&\phantom{iiMMMMMMM}\text{intersects more than 
$K_1\ep^{1/(d+1)}(t+\l)/[4\wt C_1\De_r]$ 
blocks}\\
&\phantom{iiMMMMMMM}\text{$\cB_r(\bold i,k)$ for which 
$G_r(2,z, \bold i,k)$ occurs}\}\\
&\le  \sum_{0 \le z < \wt 4C_1} \sum_{(\bold a,b)} \sum_{0 \le \l < C_{13}t}
P\{\text{there exists a path $\wh \pi \in \Xi(J,\l,t)$
such that  $\wh \pi$}\\
&\phantom{iiMMMMMMM}\text{intersects more than
$K_1\ep^{1/(d+1)}(t+\l)/[8(12)^d\wt C_1\De_r]$ blocks}\\
&\phantom{iiMMMMMMM}\text{$\cB_r(\bold i,k)$ with 
$(\bold i,k) \equiv (\bold a,b)$
for which $G_r(2,z,\bold i,k)$ 
occurs}\}\\
&\le \sum_{(\bold a,b)} \sum_{0 \le \l < C_{13}t}K_2\exp\big[-K_3
\frac{(t+\l)}{\De_r} \ep^{1/(d+1)}\big]. 
\endalign
$$
For $\ep$ so small that $K_1(1+C_{13})\ep^{1/(d+1)} 
< 1/10$, this gives \equ(5.1)
for $\cup_{0\le z < 4\wt C_1} G(2,z)$.

\hfill $\blacksquare$
\enddemo

The case $j=3$ of \equ(5.1) has already been handled in the proof of
Lemma 18, where we introduced $I_5$. Indeed, since simple random walk
cannot jump across $\cW_r(\bold i)$,\;
$I[\cup_{0 \le z < 4\wt C_1} G(3,z,\bold i,k)] \le I_5(r, \bold
i,k)$. 
Thus, the number of 
good blocks $\cB_r(\bold i,k)$ which intersect any given $J$-path $\wh
\pi|_{[0,t]} \in \Xi(\l,t)$ and for which the event $\cup_{0 \le z < 4\wt C_1}
G(3,z,\bold i,k)$ occurs is
bounded by $\sum_{(\bold i,k)}^{(\wh \pi,r)} I_5(r,\bold i,k)$. 
The inequality \equ(4.51uv) therefore applies,
and \equ(5.1) now follows with $\cup_{0 \le z < \wt 4C_1}G(3,z)$ 
in the place of $G(j)$
and for any fixed $r$ which is large enough that $2[\De_r]^{-2d-3}(1+C_{13})
< 1/(10\De_r)$, say for $r \ge r_3$.

Finally we turn to \equ(5.1) with $j=4$. As we shall show now, all 
the steps for this
estimate are already given in the estimates for $\sum_{(\bold
i,k)}{}^{(\wh \pi r)} I_{3,2}(r, \bold i,k)$ in the preceding
section. Define
$$
\align
&\wt G(4,z)= \wt G_r(4,z, \bold i,k)\\ 
&= \big\{\cV_r(\bold i,k) \text{ is good, but
there is a J-path from some }\\
&\phantom{Mi}\text{$(y',s')$ to $(y'',s'')$ with $y' \in \partial 
\prod_{s=1}^d [(i(s)-1)\De_r,(i(s)+2)\De_r-1]$, }\\
&\phantom{Mi}y'' \in \prod_{s=1}^d[i(s)\De_r, (i(s)+1)\De_r)
\text{ and $v_r(k,z-1)\le s' \le s'' \le v_r(k,z+1)$,}\\
&\phantom{Mi}\text{and this J-path uses only
particles which were in $V_r(\bold i)$ at time $v_r(k,z-1)$}\big\}.
\endalign
$$
Then, by the definition of a good block $\cB_r(\bold i,k)$ and its 
pedestal $\cV_r(\bold i,k) = V_r(\bold i) \times \{(k-1)\De_r\}$,
we have $G(4,z) \subset \wt G(4,z)$.
Moreover, we have the following analogue of \equ(4.41)
for $n = \De_r^{3(d+1)^2}$:
$$
\align
&P\{\wt G_r(4,z,\bold i,k) \text{ occurs }\big|
\cJ_r\big(\bold i,(k-1)\De_r\big)\}\\
&\le \frac {(7^d \De_r^d+1)2\De_rD^2}{\De_r^{3(d+2)^2}}\\
&+\sum \Sb v_r(k,z)-1/n \le j/n \\ \le v_r(k,z+1)+1/n \endSb \;\sum_{x \in
\partial \prod_{s=1}^d 
[(i(s) - 1)\De_r, (i(s)+2)\De_r-1]}
P\big\{\cV_r(\bold i,k) \text{ is good and}\\
&\phantom{MMMMM}\text{there exists a J-path from 
$\big(x,j/n)$ to $\prod_{s=1}^d[i(s)\De_r, 
(i(s)+1)\De_r)$}\\
&\phantom{MMMMM}\text{of time duration $\le \De_r/(2\wt C_1)+1/n \le 
\De_r/\wt C_1$ and
which uses only}\\
&\phantom{MMMMM}\text{particles which are in $V_r(\bold i)$ at 
time $v_r(k,z-1)$}\big|
\cJ_r\big(\bold i,(k-1)\De_r\big)\big\}.
\teq(5.7)
\endalign
$$
We can then follow the proof of Lemma 20 from \equ(4.41) on to obtain that
the left hand side of \equ(5.7) is at most
$$
\align
&\frac {(7^d \De_r^d+1)2\De_r D^2}{\De_r^{3(d+1)^2}}\\
&\phantom{MMM} + \sum \Sb v_r(k,z)-1/n \le j/n\\ \le v_r(k,z+1)+1/n
\endSb 
\;\sum_{x \in
\partial \prod_{s=1}^d 
[(i(s) - 1)\De_r, (i(s)+2)\De_r-1]}
4e^{ -\De_r/\wt C_1}. 
\teq(5.7aa)
\endalign
$$
In turn, it is easy to see that there exists
an $r_3$ such that for $r \ge r_3$ the expression \equ(5.7aa) is 
for each $0 \le z < 4\wt C_1$ 
at most $\De_r^{-2(d+1)^2}$ (as in Lemma 20).
We then also obtain for $r \ge r_4$
$$
\align
&P\big \{\cB_r(\bold i,k) \text{ is good and $G_r(4,z,\bold i,k)$
occurs for some $0 \le z < 4\wt C_1$}
\big| \cJ_r\big(\bold i, (k-1)\De_r\big)\}\\ 
&\le \sum_{0 \le z < 4\wt C_1}P\big \{\wt G_r(4,z,\bold i,k)
\big| \cJ_r\big(\bold i, (k-1)\De_r\big)\} \le 4\wt C_1 \De_r^{-(2(d+1)^2}. 
\endalign
$$
As in the last lemma the collection of random variables 
$$
\wt Y_r(\bold i,k) :=I[\wt G_r(4,z,\bold i,k) \text{ occurs for some 
$0 \le z < 4\wt C_1$}]
$$
with $(\bold i,k) \equiv (\bold a, b)$ lies stochastically below a
family of independent random variables $Z_r(\bold i,k)$ satisfying
$$
P\{Z_r(\bold i,k) = 1\} = 1- P\{Z_r(\bold i,k)=0\} =
4\wt C_1\De_r^{-2(d+1)^2}.
$$
Again we can now follow the proof of Lemma 11 or (5.43) in \cite {KSa}
to conclude that (for $r \ge r_4$)
$$
P\big\{\sup_{\wh \pi \in \Xi(\l,t)} \sum_{(\bold i,k)}{}^{(\wh \pi,r)}
Y_r(\bold i,k) \ge K_{15} \frac{(t+\l)}{ \De_r^{2d+3}}\big\} \le K_{16} 
\exp\big[-K_{17}\frac{(t+\l)}{ \De_r^{2d+3}}\big]
$$
($\sum_{(\bold i,k)}^{(\wh \pi,r)}$ is as in \equ(4.11)). 
If we take $r_4$ such that $K_{15}(1+C_{13})[\De_{r_4}]^{-2d+3} \le
1[10\De_{r_4}]^{-1}$, then
\equ(5.1) for any $r \ge r_5$ and 
with $G(j)$ replaced by $\cup_{0 \le z < 4\wt C_1} G_r(4,z,\bold i,k))$
is an immediate consequence.

Because the case $k=0$, was excluded in Lemma 25 we still need an
estimate for the sup over $\wh \pi \in\Xi(J,\l,t), \wh \pi(0) = \bold 0$, of 
the number of blocks $\cB_r(\bold i,0)$ which intersect $\wh
\pi$. There are at most $t/(10\De_r)$ blocks $\cB_r(\bold i,0)$ with
$\|\bold i\|\le K_{22}t^{1/d}$, with $K_{22}$ some constant which
depends on $d$ and $\De_r$ only. If there is a block $\cB_r(\bold i,0)$ with
$\|i\| >   K_{22}t^{1/d}$ which intersects $\wh \pi$, then some initial piece
of $\wh \pi$ forms a $J$-path from $\bold 0$ 
to the outside of $\cC( K_{22}t^{1/d})$. Since all points in blocks
$\cB_r(\bold i,0)$ have time coordinate less than $\De_r \le 
[2 C_1]^{-1}K_{22}t^{1/d}$ (for large $t$), we obtain
by means of \equ(4.17), \equ(4.17z) and by (1.3) in \cite{KSb}
$$
\align
&P\{\text{there exists some $\wh \pi \in \Xi(J,\l,t)$ with $\wh \pi(0) =
\bold 0$ such that $\wh \pi$ intersects more}\\
&\phantom{MMMMMMMMMMMMMMMM}\text{than $t/(10\De_r)$ $r$-blocks}\}\\
&\le P\{\text{there exists a $J$-path from $\bold 0$ 
to the outside of $\cC( K_{22}t^{1/d})$ of time}\\
&\phantom{MMMMMMMMMMMMMMMM}\text{duration less than 
$[2C_1]^{-1}K_{22}t^{1/d}$}\\
&\le 4\exp \big[-[2 C_1]^{-1} K_{22}t^{1/d}\big].
\teq(5.20)
\endalign
$$
Thus, also \equ(5.1z) holds.

We now take $r = \max\{r_i:1 \le i \le 5\}$ and $\la \ge \la_0(r)$. 
Then \equ(5.1) holds for $1 \le j \le 4$ and also \equ(5.1z) holds.
As discussed right after the statement of Proposition 24, these
properties imply Proposition 24.

\enddocument